\documentclass[a4paper, 10pt]{scrartcl}

\usepackage[top=2.5cm, bottom=2.5cm, left=2.5cm, right=2.5cm]{geometry}

\usepackage{moreverb}
\usepackage{cite}
\usepackage{amsmath, amsfonts, mathtools}

\usepackage[colorlinks=false,bookmarksopen,bookmarksnumbered,citecolor=black,urlcolor=black]{hyperref}

\newcommand\BibTeX{{\rmfamily B\kern-.05em \textsc{i\kern-.025em b}\kern-.08em
T\kern-.1667em\lower.7ex\hbox{E}\kern-.125emX}}

\makeatletter
\let\save@mathaccent\mathaccent
\newcommand*\if@single[3]{%
 \setbox0\hbox{${\mathaccent"0362{#1}}^H$}%
 \setbox2\hbox{${\mathaccent"0362{\kern0pt#1}}^H$}%
 \ifdim\ht0=\ht2 #3\else #2\fi
 }
\newcommand*\rel@kern[1]{\kern#1\dimexpr\macc@kerna}
\newcommand*\widebar[1]{\@ifnextchar^{{\wide@bar{#1}{0}}}{\wide@bar{#1}{1}}}
\newcommand*\wide@bar[2]{\if@single{#1}{\wide@bar@{#1}{#2}{1}}{\wide@bar@{#1}{#2}{2}}}
\newcommand*\wide@bar@[3]{%
 \begingroup
 \def\mathaccent##1##2{%
	\let\mathaccent\save@mathaccent
	\if#32 \let\macc@nucleus\first@char \fi
	\setbox\z@\hbox{$\macc@style{\macc@nucleus}_{}$}%
	\setbox\tw@\hbox{$\macc@style{\macc@nucleus}{}_{}$}%
	\dimen@\wd\tw@
	\advance\dimen@-\wd\z@
	\divide\dimen@ 3
	\@tempdima\wd\tw@
	\advance\@tempdima-\scriptspace
	\divide\@tempdima 10
	\advance\dimen@-\@tempdima
	\ifdim\dimen@>\z@ \dimen@0pt\fi
	\rel@kern{0.6}\kern-\dimen@
	\if#31
	  \overline{\rel@kern{-0.6}\kern\dimen@\macc@nucleus\rel@kern{0.4}\kern\dimen@}%
	  \advance\dimen@0.4\dimexpr\macc@kerna
	  \let\final@kern#2%
	  \ifdim\dimen@<\z@ \let\final@kern1\fi
	  \if\final@kern1 \kern-\dimen@\fi
	\else
	  \overline{\rel@kern{-0.6}\kern\dimen@#1}%
	\fi
 }%
 \macc@depth\@ne
 \let\math@bgroup\@empty \let\math@egroup\macc@set@skewchar
 \mathsurround\z@ \frozen@everymath{\mathgroup\macc@group\relax}%
 \macc@set@skewchar\relax
 \let\mathaccentV\macc@nested@a
 \if#31
	\macc@nested@a\relax111{#1}%
 \else
	\def\gobble@till@marker##1\endmarker{}%
	\futurelet\first@char\gobble@till@marker#1\endmarker
	\ifcat\noexpand\first@char A\else
	  \def\first@char{}%
	\fi
	\macc@nested@a\relax111{\first@char}%
 \fi
 \endgroup
}
\makeatother

\newcommand{\ie}{i.e.\ }
\newcommand{\eg}{e.g.\ }
\newcommand{\eq}[1]{\begin{align}#1\end{align}}
\newcommand{\bs}[1]{{\boldsymbol{#1}}}
\newcommand{\mc}[1]{{\mathcal{#1}}}
\newcommand{\mf}[1]{{\mathfrak{#1}}}
\newcommand{\mrm}[1]{{\mathrm{#1}}}
\newcommand{\mbb}[1]{{\mathbb{#1}}}
\newcommand{\tx}[1]{{\text{#1}}}
\newcommand{\field}[1]{\mbb{#1}}
\newcommand{\total}{\mrm{d}}
\newcommand{\pd}[2]{\frac{\partial #1}{\partial #2}}

\newcommand{\cotanspace}[2]{ T^*_{#1} \mc{#2} }
\newcommand{\tanbndl}[1]{ T \mc{#1} }
\newcommand{\cotanbndl}[1]{ T^* \mc{#1} }

\newcommand{\cobasis}[2]{ \bs{\total} {#1}^{#2} }
\newcommand{\extd}{\bs{\total}}

\newcommand{\ital}{\normalfont \itshape}
\newcommand{\mg}{\bs{g}}

\newcommand{\simplx}[1]{\widebar{#1}}
\newcommand{\clssimplx}[2]{\widebar{#1}_{\bs{#2}}}
\newcommand{\ssimplx}[2]{#1_{\bs{#2}}}
\newcommand{\refsimplx}[1]{\hat{#1}}
\newcommand{\Dsimplx}[1]{\Delta\widebar{{#1}}}
\newcommand{\Dksimplx}[2]{\Delta_{#2}\widebar{{#1}}}
\newcommand{\freud}[1]{\mathcal{F}(#1)}
\newcommand{\Dfreud}[1]{\Delta\mathcal{F}(#1)}
\newcommand{\Dkfreud}[2]{\Delta_{#2}\mathcal{F}\left(#1\right)}
\newcommand{\incfreud}[1]{\widetilde{\mathcal{F}}(#1)}
\newcommand{\Dkincfreud}[2]{\mc{S}_{#2}\widetilde{\mathcal{F}}(#1)}
\newcommand{\hchildren}[1]{\mathcal{C}(#1)}
\newcommand{\leaf}[1]{\mathfrak{I}(#1)}
\newcommand{\upadjac}[1]{\mrm{anc}\left(#1\right)}
\newcommand{\upadjacstar}[1]{\mrm{anc}_\star\left(#1\right)}

\newcommand{\extensplus}[4]{E^{#1+}_{\ssimplx{#2}{#3},\simplx{#4}}}
\newcommand{\extensminus}[4]{E^{#1-}_{\ssimplx{#2}{#3},\simplx{#4}}}
\newcommand{\PLplus}[3]{\mc{P}_{#1}\Lambda^{#2}(\simplx{#3})}
\newcommand{\PLminus}[3]{\mc{P}^-_{#1}\Lambda^{#2}(\simplx{#3})}
\newcommand{\Tr}[2]{\tx{Tr}_{\ssimplx{#1}{#2}}}
\newcommand{\fcspcminus}[5]{V_{#1}^{#2-}(\simplx{#3},\ssimplx{#4}{#5})}
\newcommand{\fcspcplus}[5]{V_{#1}^{#2+}(\simplx{#3},\ssimplx{#4}{#5})}

\newcommand{\trspcminus}[4]{W_{#1}^{#2-}(\ssimplx{#3}{#4})}
\newcommand{\trspcplus}[4]{W_{#1}^{#2+}(\ssimplx{#3}{#4})}

\usepackage{scalerel}

\DeclareMathOperator*{\assembly}{\scalerel*{\mathsf{A}}{\bigoplus}}

\usepackage{pgfplots}
\pgfplotsset{compat=1.17}
\usetikzlibrary{positioning}
\usetikzlibrary{calc}
\usepackage{mathbbol}
\usepackage{enumerate}
\usepackage{tensor}
\usepackage{tikz}
\usetikzlibrary{decorations.pathmorphing,arrows,intersections}
\tikzset{zigzag/.style={decorate, decoration=zigzag}}

\usepackage{contour}

\usepgfplotslibrary{colormaps}
\pgfplotsset{
	 colormap={unityjet}{
		  rgb255=(107,2,40)
		  rgb255=(163,80,44)
		  rgb255=(169,160,98)
		  rgb255=(132,155,98)
		  rgb255=(70,126,102)
		  rgb255=(62,54,99)
	 },
	 colormap={arnold}{
rgb255=(39.015,99.96,24.99)
rgb255=(77.009995,146.115,32.894997)
rgb255=(126.99,187.935,65.025)
rgb255=(184.11,224.91,133.875)
rgb255=(230.01001,245.05501,208.08)
rgb255=(247.095,247.095,247.095)
rgb255=(252.95999,223.89,238.935)
rgb255=(240.97499,182.06999,218.02501)
rgb255=(222.105,119.085,173.90999)
rgb255=(197.115,27.029999,124.950005)
rgb255=(142.03499,1.0200001,82.11)
	 },
		  colormap={BuRd}{
rgb255=(5.1,47.94,96.9)
rgb255=(32.894997,102.0,172.125)
rgb255=(67.065,146.87999,195.075)
rgb255=(146.115,197.115,222.105)
rgb255=(209.09999,228.99,239.955)
rgb255=(247.095,247.095,247.095)
rgb255=(252.95999,219.04501,198.9)
rgb255=(244.035,164.985,130.05)
rgb255=(213.94499,95.88,77.009995)
rgb255=(177.99,23.97,43.095)
rgb255=(103.020004,0.0,31.11)
	 },
				colormap={BkOr}{
rgb255=(0.0,0.0,0.0)
rgb255=(174.0,112.0,0.0)
	 },
					 colormap={Purples}{
rgb255=(251.93999,250.92001,252.95999)
rgb255=(238.935,236.895,245.05501)
rgb255=(218.02501,218.02501,235.11)
rgb255=(187.935,188.955,220.06499)
rgb255=(158.1,154.01999,199.92)
rgb255=(128.01,124.950005,185.89499)
rgb255=(106.08,81.09,162.94499)
rgb255=(83.895,39.015,143.055)
rgb255=(62.984997,0.0,124.950005)
	 },
}

\definecolor{mylightred}{RGB}{211,79,73}
\definecolor{mydarkred}{RGB}{199,44,38}
\definecolor{mylightgreen}{RGB}{78,153,67}
\definecolor{mydarkgreen}{RGB}{43,129,33}
\definecolor{mylightpurple}{RGB}{150,107,178}
\definecolor{mydarkpurple}{RGB}{126,78,160}
\definecolor{mylightblue}{RGB}{49,101,205}
\definecolor{mydarkblue}{RGB}{20,92,205}

\tikzset{
  juliadot/.style args={#1,#2}{shape=circle,line width=0.03ex,minimum width=0.4ex,fill=#1,draw=#2}
}

\begin{document}

\title{An hp-hierarchical framework for the \\finite element exterior calculus\thanks{The authors acknowledge the generous support of Udo Nackenhorst at the Institut f\"{u}r Baumechanik und Numerische Mechanik, Leibniz Universit\"{a}t Hannover, Appelstr.\ 9A, 30167 Hannover, Germany for this project. This paper is dedicated to the memory of Erwin Stein, whom the authors thank for the many fruitful discussions on this project. Finally, the authors thank Christoph Schwab for pointing out an earlier reference~\cite{Braess2000} to the polynomial eigenvalue problem on simplices which had originally been attributed to a later work~\cite{Braess2005} in an earlier version of this article.} }

\author{Robert L.\ Gates\thanks{Correspondence to: Robert Gates (robert.gates@ibnm.uni-hannover.de)} \ and Maximilian Bittens\thanks{Maximilian Bittens (maximilian.bittens@bgr.de)}}

\date{December 31, 2020}
\maketitle

\begin{abstract}
\noindent The problem of solving partial differential equations (PDEs) on manifolds can be considered to be one of the most general problem formulations encountered in computational multi-physics. The required covariant forms of balance laws as well as the corresponding covariant forms of the constitutive closing relations are naturally expressed using the bundle-valued exterior calculus of differential forms or related algebraic concepts. It can be argued that the appropriate solution method to such PDE problems is given by the finite element exterior calculus (FEEC). The aim of this essay is the exposition of a simple, efficiently-implementable framework for general $hp$-adaptivity applicable to the FEEC on higher-dimensional manifolds. A problem-independent spectral error-indicator is developed which estimates the error and the spectral decay of polynomial coefficients. The spectral decay rate is taken as an admissibility indicator on the polynomial order distribution. Finally, by elementary computational examples, it is attempted to demonstrate the power of the method as an engineering tool.
\par\vskip\baselineskip\noindent
\textbf{Keywords: Mixed finite element methods, Finite element exterior calculus, Hierarchical basis, $hp$-adaptivity, Riemannian manifolds, Incompressible elasticity.}
\end{abstract}

\vspace{-6pt}

\section{Introduction}
\vspace{-2pt}

Partial differential equations (PDEs) are ubiquitous in the mathematical modeling of physical, environmental, and financial phenomena. A number of methods have been developed for the numerical approximation of such equations, among them the finite difference, finite volume, and finite element methods. These methods essentially are concerned with means for representing discretely the quantities and differential operators involved in the PDE. Unfortunately, for many types of problems, such numerical methods have been plagued by various shortcoming which essentially are related to the issues of discrete consistency and stability~\cite{Arnold2006}. These shortcomings can lead to numerical solutions which are polluted by artifacts or are entirely incorrect. Depending on the experience of the user and the complexity of the problem, some of these artifacts may be difficult to detect by inspection. While numerous so-called stabilizations~\cite{Simo1990, Tezduyar1991, Tobitzka1991} have been developed and applied successfully to special cases of PDEs and applications, a general remedy had remained elusive throughout the past century. In addition to the stability of the numerical scheme, the physical accuracy of numerical solutions to PDEs is governed by the ability to ensure local conservation in the discrete case. In the past, the attainment of discrete local conservation had also been difficult, as it involved not only the correct mathematical modeling, but also the choice of numerical method.

Beginning in the early second half of the last century in electromagnetism with Kron's method of tearing~\cite{Kron1953} and in mathematics with Whitney's geometric integration theory~\cite{Whitney1957}, a confluence of research~\cite{Branin1966, Tonti1972, Tonti1974, Dodziuk1974, Dodziuk1976, Bossavit1988} in the field of numerical analysis of PDEs and the field of differential geometry led to the realization that the algebraic structure of the partial differential equation must be preserved in the discrete case in order to guarantee stability and conservation of the numerical scheme. In the field of finite elements this led to the Raviart-Thomas and Brezzi-Douglas-Marini elements, as well as variational time integrators. More recent advancements such as the discrete exterior calculus~\cite{Hirani2003}, higher-order Whitney forms~\cite{Hiptmair2001a} and the famed finite element exterior calculus~\cite{Arnold2006} generalized these developments. Other extensions include those of subdivision Whitney forms~\cite{Wang2008} and de Rham-compatible isogeometric analysis~\cite{Buffa2011}. During the same period, in finite element methods, the efficiency of $hp$-discretizations became widely accepted in engineering practice. However, to this date, fully general $hp$-methods for structure-preserving discretization have not been considered. Most developments confine themselves to the discretization of the de Rham complex in up to three dimensions. The present work aims to outline a method of going the extra mile, generalizing to arbitrary dimension and including the elasticity complex.

The motivation for this work is to obtain a simple and unified computational framework for $hp$-adaptivity on piecewise smooth Riemannian manifolds, general enough to be applied to multi-physics problems in arbitrary dimensions. The appropriate unifying framework selected as a basis for the following developments is that of the finite element exterior calculus. The generalization to arbitrary dimension of course precludes ``refinement rules'' such as red/green-refinement in two dimensions, which leads us to the requirement of a hierarchical basis. Thus, building upon the works by Yserentant (1985)~\cite{Yserentant1985} and Grinspun (2003)~\cite{Grinspun2003}, this work develops curvilinear cell complexes with local hierarchical subdivisions and non-uniform polynomial order distribution. By construction, such a method entirely avoids the introduction of hanging mesh entities (such as hanging vertices, lines, faces, etc.). In order to obtain a linearly-independent basis compatible with a non-uniform hierarchical subdivision of the cell complex, we use the refinement equation to introduce linear constraints on Whitney forms which are entirely local to the subdivided subcells. 

While the presented method is also in principle applicable to polygonal cells, this work restricts itself to simplicial complexes subdivided using an appropriately-modified Freudenthal simplicial splitting~\cite{Freudenthal1942, Edelsbrunner2000, Goncalves2006} for reasons of simplicity in implementation. A $p$-hierarchical decomposition of the finite element exterior calculus bases on simplices amenable to implementation is given along with a simple and general algorithm for its computation. This decomposition is used to define finite element spaces of non-uniform polynomial order on the simplicial complex. The authors have recently discovered that non-uniform polynomial order distributions in the finite element exterior calculus were analyzed by Licht (2017)~\cite{Licht2017} in parallel to this project, where a single-level minimum rule was employed. This work extends the classical minimum rule to a variant which is compatible with the non-uniform hierarchical subdivision.

Building upon the $p$-hierarchical decomposition of the bases, we develop a problem-independent $n$-simplex version of the one-dimensional spectral error indicator given by Mavriplis (1994)~\cite{Mavripilis1994} by utilizing a polynomial eigenvalue problem given by Braess (2005)~\cite{Braess2005}. Apart from the solution coefficients and the element mapping, all operations required in evaluating the proposed a posteriori error indicator are entirely local to the reference element. As in the original version of the spectral error indicator, the spectral decay rate is a by-product of the error-estimation process and gives a measure for the solution regularity. This regularity measure is used to determine whether to refine in $h$ or in $p$. The emphasis here is placed on generality, we thus make no reference to improved efficiency attainable for specific problems by other methods.

Throughout this work we give an outline of our implementation of the method for arbitrary dimensions in the Julia language. Finally, we present numerical examples in up to three dimensions, including a newly-proposed reference solution of the Cook's membrane problem in incompressible, linearized plane-strain elasticity.

\section{Fundamentals}
We begin with some fundamentals on the exterior calculus of differential forms, simplices as well as on the assembly and the inner product of finite element basis functions. To this end, let a {\ital $k$-vector} be defined as a ${0\choose k}$-tensor of an $n$-dimensional vector space $V$ over another vector space $W$. Its tensor-representation is antisymmetric with respect to all permutations of the indices and is defined as the multi-linear map
\eq{
\underbrace{V^* \times \ldots \times V^*}_{k-\tx{fold}} \rightarrow W \,.
}
A $k$-vector is an element of a vector space, the {\ital $k$-th exterior power} $\Lambda^k(V; W)$. Writing $\Lambda^k(V)$ in place of $\Lambda^k(V; \field{R})$, the exterior power $\Lambda^k(V; W)$ can be identified with the tensor product space $W \otimes \Lambda^k(V)$. Exterior powers are subspaces of the {\ital exterior algebra} $\Lambda(V)$ generated by the {\ital exterior product} $\wedge:\Lambda^k\times \Lambda^p \rightarrow \Lambda^{k+p}$. The exterior product has the properties
\eq{
	\bs{\omega} \wedge \bs{\nu} = (-1)^{kp}\,  \bs{\nu} \wedge \bs{\omega} \,,
}
and $\bs{\omega} \wedge \bs{\omega} = 0$ for $k$ odd. Since the exterior algebra follows the rules of the determinant, $k$-vectors are useful in representing oriented physical quantities, such as flux or circulation. A $k$-vector inherits properties of tensors. For vector-valued forms, we will distinguish wedge products by the operator used to treat the vector values, \eg for an appropriate inner product $\cdot : W \times W \rightarrow \field{R}$ we would write $\overset{\cdot}{\wedge} : \Lambda^k(V; W) \times \Lambda^k(V; W) \rightarrow \Lambda^{k+p}(V; \field{R})$. For $\otimes : W \times W \rightarrow W \otimes W$ we would then have $\overset{\otimes}{\wedge} : \Lambda^k(V; W) \times \Lambda^k(V; W) \rightarrow \Lambda^{k+p}(V; W\otimes W)$.

On a Riemannian manifold $\mc{M}$, we work with an exterior algebra over cotangent spaces and bundles. The basis of the $k$-th exterior power $\Lambda^k(\cotanspace{P}{M})$ also denoted succinctly as $\Lambda^k(\mc{M})$ is given locally in a chart about $P \in \mc{M}$ by wedge-multiples of the tangent basis covectors $\cobasis{x}{i}$:
\eq{
\{\cobasis{x}{\bs{\sigma}} \;|\; \bs{\sigma} \in \Sigma(k, n)\}, \quad \cobasis{x}{\bs{\sigma}} = \cobasis{x}{\sigma(1)} \wedge \ldots \wedge \cobasis{x}{\sigma(k)}\,,
}
where $\Sigma(k, n)$ is the set of all ${n \choose k}$ combinations of $1,\ldots,n$. Hence, with coefficients $\omega_{\bs{\sigma}}$, the local form is
\eq{
	\bs{\omega} = \sum\limits_{\bs{\sigma} \in \Sigma(k,n)} \omega_{\bs{\sigma}} \, \cobasis{x}{\bs{\sigma}} \, .
}
A {\ital differential $k$-form} on $\mc{M}$ is obtained by letting $P$ vary and thus is an element of the exterior power $\Lambda^k(\cotanbndl{M})$ on the cotangent bundle.

Given a metric $\bs{g}$ at point $P$ on $\mc{M}$, we define the {\ital Hodge star} point-wisely through the local inner-product of differential forms
\eq{
	\langle\bs\omega, \bs\nu \rangle_{\Lambda^k(\cotanspace{P}{M})} = \sum\limits_{\bs{\sigma} \in \Sigma(k,n)}\omega_{\bs\sigma} \nu_{\bs\sigma}\gamma_{\bs\sigma} = \bs\omega \wedge\star_{\mg}  \bs\nu \,,
}
where $\gamma_{\bs\sigma}$ is a geometric coefficient equal to one in the Euclidean case. The Hodge star is thus defined as the linear map $\star_{\mg} : \Lambda^{k}(\cotanspace{P}{M}) \rightarrow \Lambda^{n-k}(\cotanspace{P}{M})$ such that with $\bs \nu = (1,\ldots,n)\setminus \bs \sigma$, in coordinates,
\eq{
	{\star_{\mg}\cobasis{x}{\bs \sigma}} =  \sqrt{|\det(g_{ij})|}\;\mrm{sign}((\bs \sigma, \bs \nu)) \sum\limits_{\bs\kappa \in S(\bs\sigma)} \mrm{sign}(\bs \kappa) \prod\limits_{i = 1}^k g^{\bs\sigma(i)\bs\kappa(i)} \cobasis{x}{\bs \nu} \label{eq:hodgeStar}
}
where $S(\bs\sigma)$ is the set of all permutations of $\bs{\sigma}$. The Hodge star transforms naturally under general linear group transformations, as long as the metric is also transformed. 

The {\ital exterior derivative} is defined as the map $\extd : \Lambda^{k}(\cotanbndl{M}) \rightarrow \Lambda^{k+1}(\cotanbndl{M})$ between differential forms, which for $\omega \in \Lambda^0(\cotanbndl{M})$, is simply the total differential:
\eq{
	\extd \omega = \pd{\omega}{x^1} \,\cobasis{x}{1} + \ldots + \pd{\omega}{x^n} \,\cobasis{x}{n} \, .
}
For $\bs{\omega} \in \Lambda^k$, the definition 
\eq{
	\extd\bs{\omega} = \extd \omega_{\bs{\sigma}} \wedge \cobasis{x}{{\bs{\sigma}}}\, 
}
then follows from the wedge product.
It is nilpotent: $\extd \circ \extd = 0$ and satisfies the Leibniz rule, \ie for $\bs{\omega} \in \Lambda^{k}$, $\bs{\nu} \in \Lambda^{p}$,
\eq{
	\extd  \left( \bs{\omega} \wedge \bs{\nu} \right) = \extd  \bs{\omega} \wedge \bs{\nu} + (-1)^k   \bs{\omega} \wedge \extd \bs{\nu} \, . \label{eq:leibnizrule}
}
Given the Hodge star and the exterior derivative, the exterior coderivative is defined as $\bs{\delta} = (-1)^{n(k+1)+1+s} {\star\extd\star}$ where $s$ is the metric signature.

\subsection{Geometric decomposition of the simplex}
\label{sec:geometric_decomposition_simplex}

Let the {\ital geometric $n$-simplex} $\simplx{T}$ be defined as the convex hull of $n+1$ points $\{\bs{\xi}_i\}_{i=0}^{n}$ in $\mathbb{R}^n$ termed the vertices of $\simplx{T}$:
\begin{align}
\simplx{T} = \{ \lambda_0 \bs{\xi}_0 + \ldots + \lambda_n \bs{\xi}_n \; : \; \sum\limits_{i=0}^n \lambda_i = 1 \; \tx{and} \; \lambda_i \ge 0 \} \, .
\end{align}
For the {\ital reference simplex} $\refsimplx{T}$, the vertices are $\bs{\xi}_0 = (0, \ldots, 0)$, $\bs{\xi}_1 = (1, 0, \ldots, 0)$, \ldots, $\bs{\xi}_n = (0, \ldots, 0, 1)$. %
With the notation $\Sigma_0(k,n)$ as the set of increasing maps $\{0,...,k\} \rightarrow \{0,...,n\}$ for $0\leq k \leq n$, we define the {\ital geometric $k$-subsimplex} as the open set of points
\begin{align}
\ssimplx{f}{\sigma} = \{ \lambda_{\sigma(0)} \bs{\xi}_{\sigma(0)} + \ldots + \lambda_{\sigma(k)} \bs{\xi}_{\sigma(k)} \; : \; \sum\limits_{i=0}^{k} \lambda_{\sigma(i)} = 1 \; \tx{and} \; \lambda_{\sigma(i)} > 0 \}\,,
\end{align}
where $\bs{\sigma} \in \Sigma_0(k,n)$ and $\{\bs{\xi}_i\}_{i \in \bs{\sigma}}$ are vertices of $\simplx{T}$. Hence, the closure $\mrm{cl}(f_{(0,\ldots,n)}) = \simplx{T}$.
 The negative counterpart of a $k$-subsimplex $\ssimplx{f}{\sigma}$ is defined as the subsimplex $-\ssimplx{f}{\sigma} = f_{\pi(\bs\sigma)}$ with the vertices $\bs{\sigma}$ permuted by an odd permutation $\pi$. %
The geometric $n$-simplex $\simplx{T}$ can thus be decomposed into its subsimplices as
\begin{align}
\simplx{T} = \bigcup\limits_{f \in \Dsimplx{T}} f \quad \text{with} \quad \Dsimplx{T} = \left\{ \ssimplx{f}{\sigma}  \, :\, \bs{\sigma} \in \Sigma_0(k,n), k \in (0,\ldots,n) \right\}.
\end{align}
We refer to the set $\Dsimplx{T}$ as the {\ital geometric descendants} of simplex $\simplx{T}$. 
\newcommand\IMAGESWIDTHBD{.22}

\begin{figure}[htb]
	\centering
	\begin{minipage}[t]{\IMAGESWIDTHBD\linewidth}
        		\begin{tikzpicture}
			\node[inner sep=0pt]{\includegraphics[trim={35cm 20cm 37cm 14cm},clip,width=\linewidth]{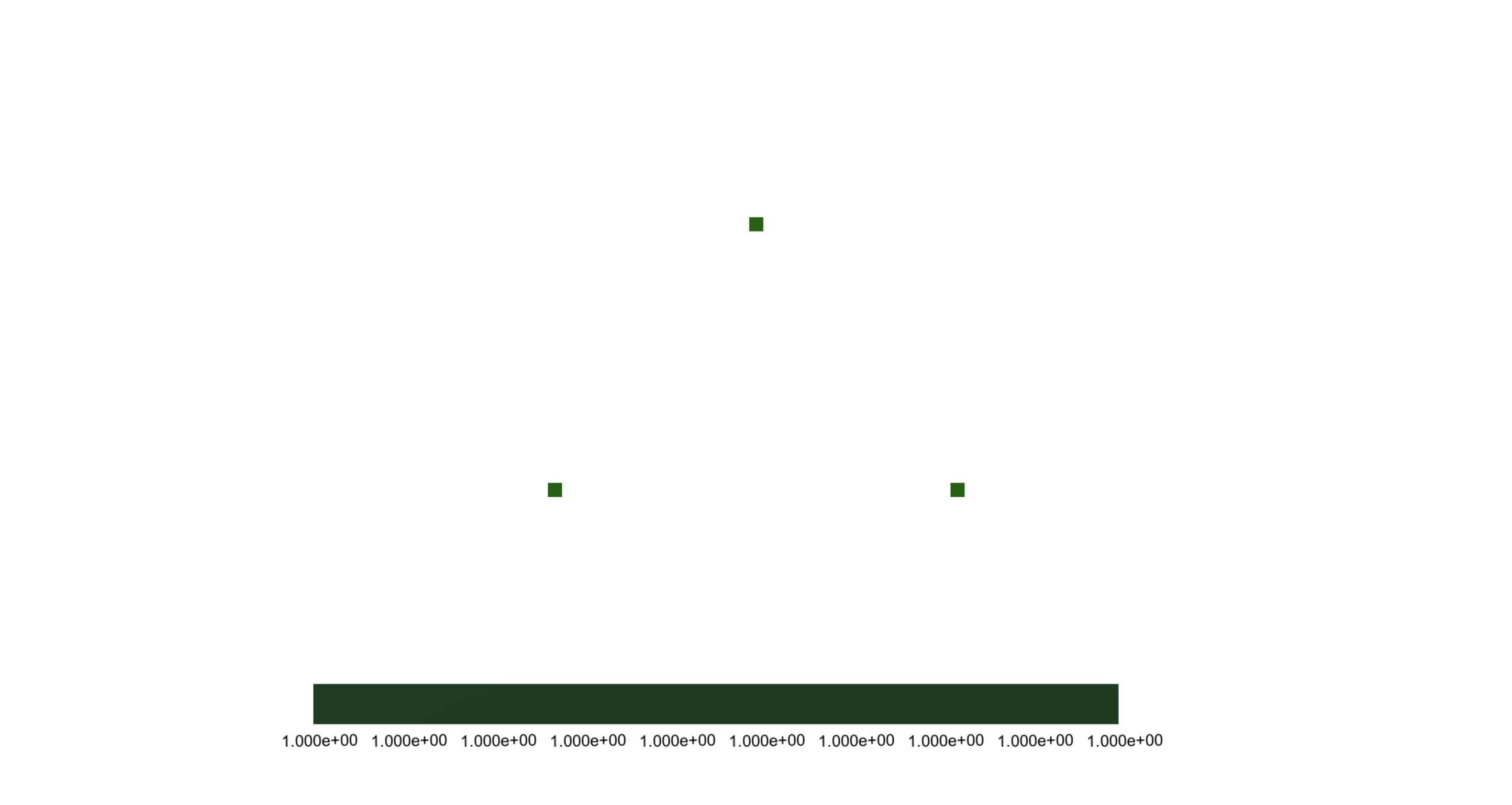}};
			\node at (-.92,-.8) {$f_{(1)}$};
			\node at (1.17,-.8) {$f_{(2)}$};
			\node at (.5,.8) {$f_{(3)}$};
		\end{tikzpicture}
		\begin{center}
        		(a) $\Dksimplx{T}{0}$
		\end{center}
    	\end{minipage}%
    	\hfill
    	\begin{minipage}[t]{\IMAGESWIDTHBD\linewidth}
		\begin{tikzpicture}
			\node[inner sep=0pt]{\includegraphics[trim={35cm 20cm 37cm 14cm},clip,width=\linewidth]{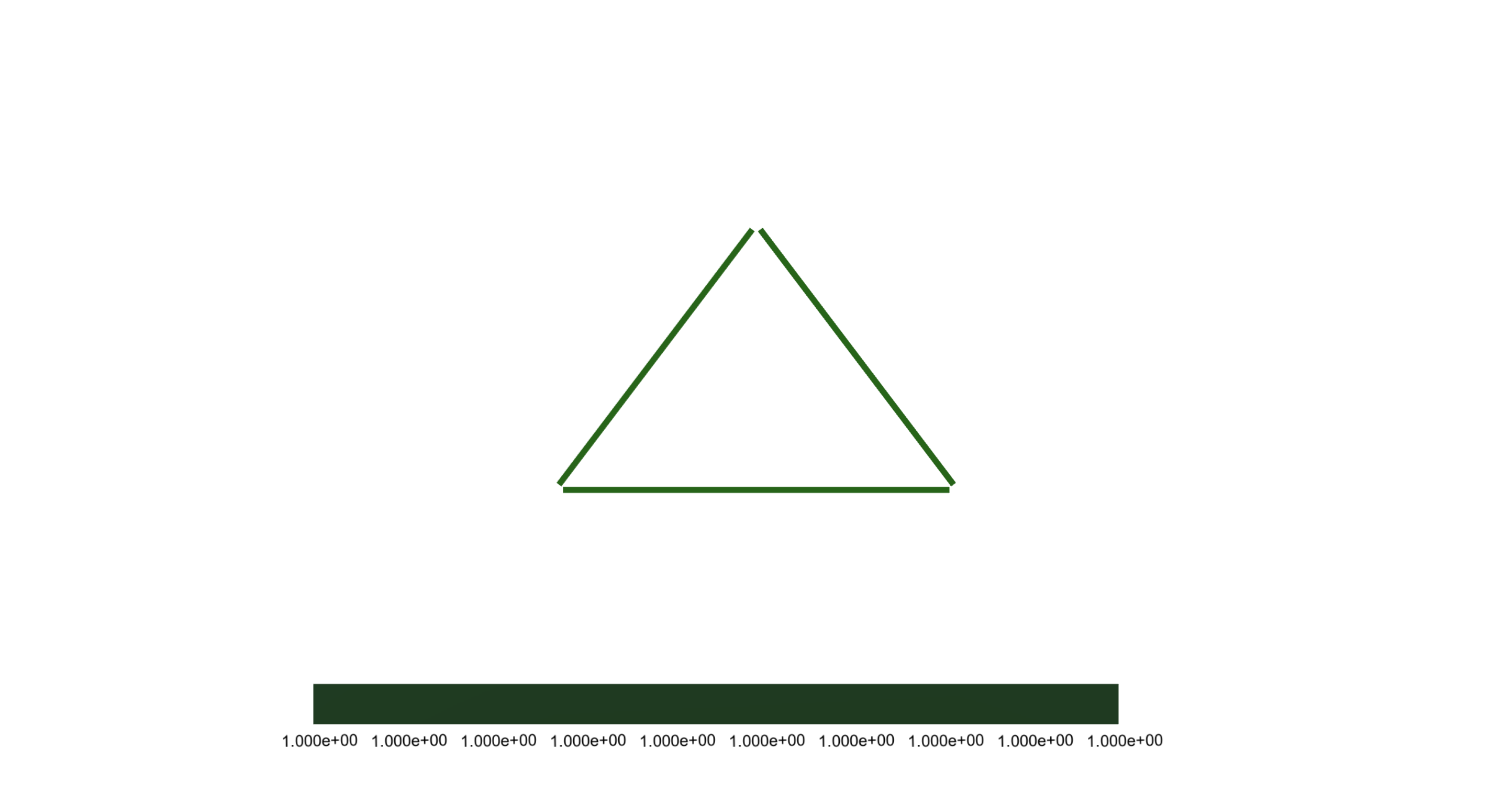}};
			\node at (-1.,.2) {$f_{(1,3)}$};
			\node at (0.2,-0.65) {$f_{(1,2)}$};
			\node at (1.3,.2) {$f_{(2,3)}$};
		\end{tikzpicture}
        		\begin{center}
        		(b) $\Dksimplx{T}{1}$		
		\end{center}
	\end{minipage}%
	\hfill
    	\begin{minipage}[t]{\IMAGESWIDTHBD\linewidth}
		\begin{tikzpicture}
			\node[inner sep=0pt]{\includegraphics[trim={35cm 20cm 37cm 14cm},clip,width=\linewidth]{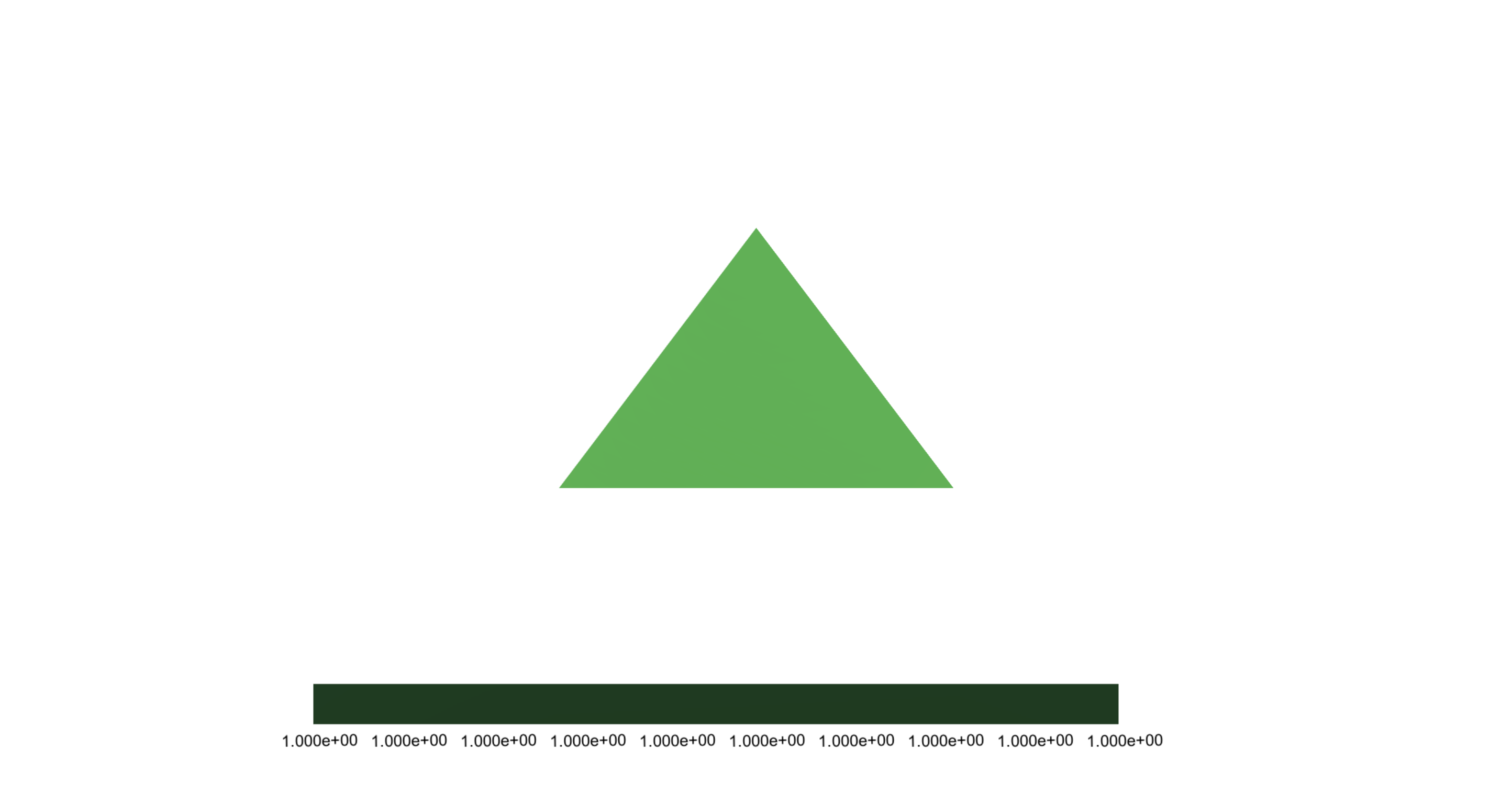}};
			\node at (0.15,-0.3) {$f_{(1,2,3)}$};
		\end{tikzpicture}
        		\begin{center}
        		(c) $\Dksimplx{T}{2}$
		\end{center}
	\end{minipage}%
	\hfill
    	\begin{minipage}[t]{\IMAGESWIDTHBD\linewidth}
		\begin{tikzpicture}
			\node[inner sep=0pt]{\includegraphics[trim={35cm 20cm 37cm 14cm},clip,width=\linewidth]{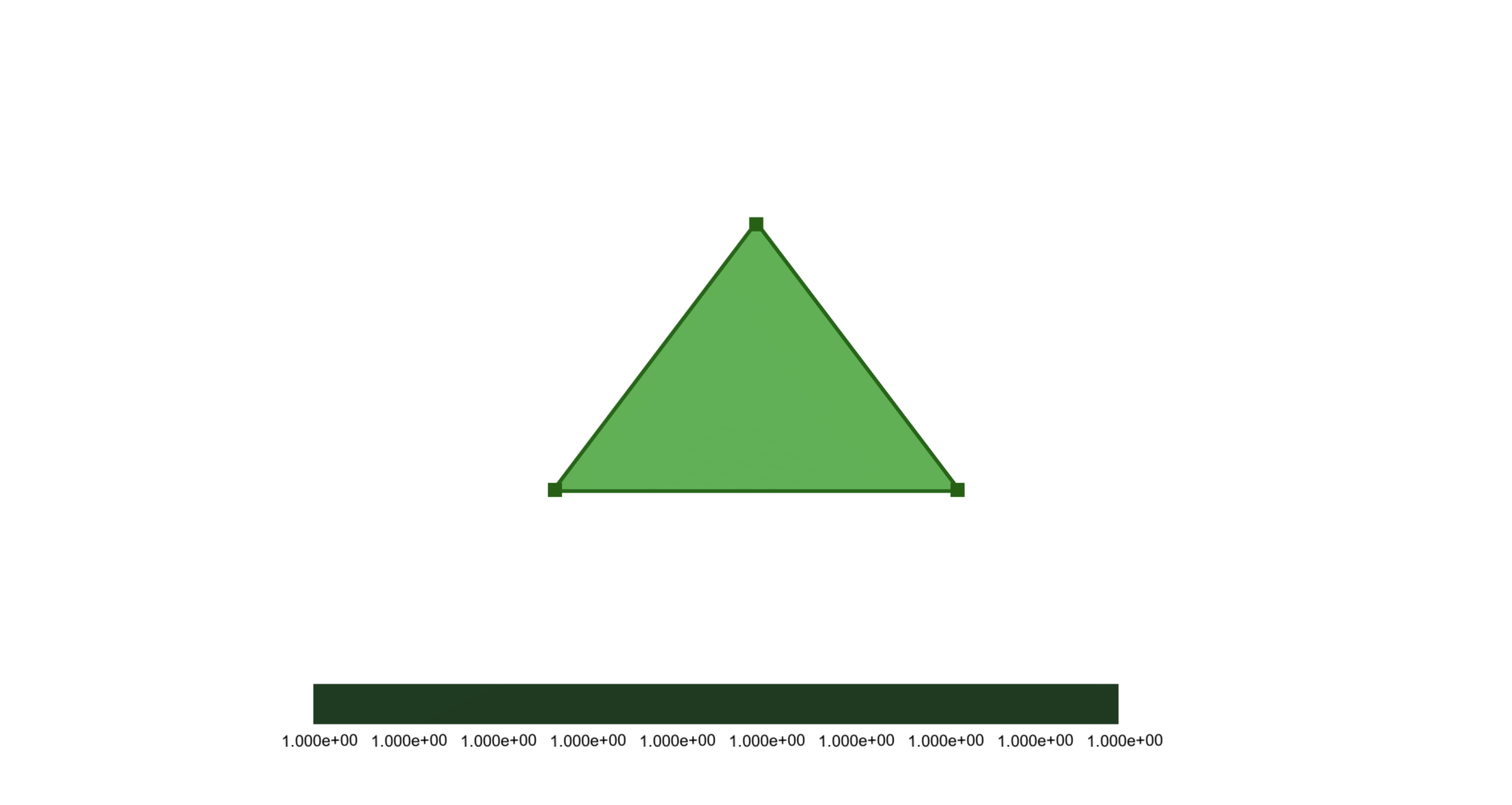}};
		\end{tikzpicture}
        		\begin{center}
        		(d) $\simplx{T} = \mrm{cl}(f_{(1,2,3)})$
		\end{center}
	\end{minipage}%
	\caption{Geometric decomposition of 2-simplex}
\end{figure}
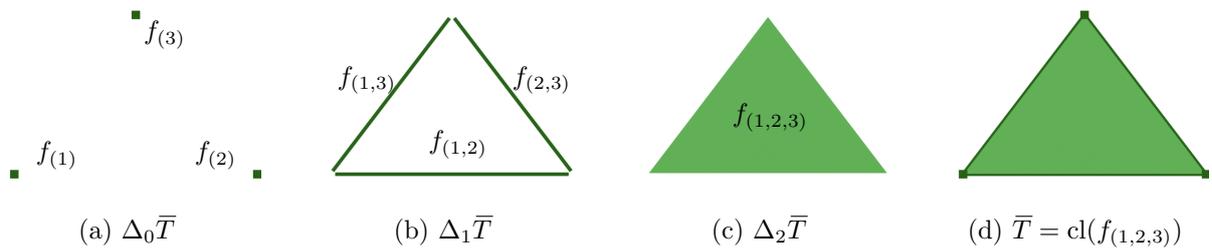

 In addition, the selector
$\mc{S}_k(\Dsimplx{T}) = \left\{ \ssimplx{f}{} \, : \, \mrm{dim}f = k\,,\, \ssimplx{f}{} \in \Dsimplx{T} \right\}$
gathers all $k$-dimensional subsimplices, with the abbreviated form $\Dksimplx{T}{k} = \mc{S}_k({\Dsimplx{T}})$. Two simplices $\simplx{U},\simplx{V}$ are termed adjacent, if $\simplx{U}\cap\simplx{V}\neq\emptyset$. %
We define as {\ital geometric ancestors} to the $k$-subsimplex $f \in \Dsimplx{T}$ the subsimplices of dimension $d > k$
\begin{align}
	\upadjac{\Dsimplx{T},f} = \{ s \; : \; f \subset \mrm{cl}(s) \wedge s \not\subseteq f, \, s \in \Dsimplx{T} \} \,,
\end{align}
where the condition $\ssimplx{f}{\sigma} \subset \mrm{cl}(\ssimplx{s}{\mu})$ is equivalent to the combinatorial condition $\bs{\sigma} \subseteq \bs{\mu}$. For notational convenience, we write
\begin{align}
	\upadjacstar{\Dsimplx{T},f} = \upadjac{\Dsimplx{T},f} \cup \{f\} \,,
\end{align}
as the set of all subsimplices of $\Dsimplx{T}$ which are geometric ancestors to $f$, including $f$ itself. Finally, we define the set-theoretic boundary of an $n$-simplex $\simplx{T}$ as the union of simplices in the set
\eq{
	\partial \simplx{T} = \{ (-1)^i \, \mrm{cl}(f_{(0,\ldots,i-1,i+1,\ldots, n)}) \; : \; i \in (0,\ldots, n) \} \,.
}

\subsection{The reference element mapping}
\begin{figure}[htb]
	\centering
	\includegraphics[width=.35\linewidth]{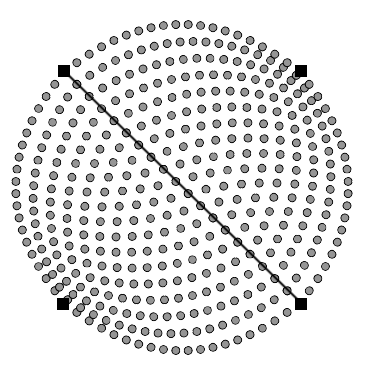}\hspace{.1\linewidth}
	\includegraphics[width=.35\linewidth]{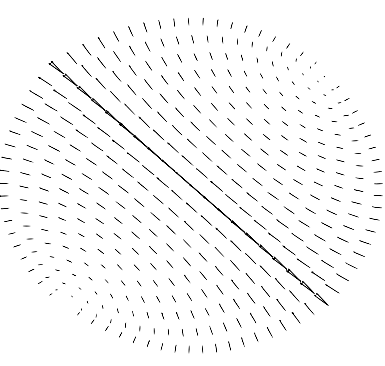}
	\caption{Cartesian grid of points on the reference element mapped onto two curvilinear triangles (vertices marked with squares, common edge marked as a black line) (left); Whitney 1-form basis function for the common edge mapped to the geometry (right)}
	\label{fig:elementmapping}
\end{figure}

In the following, we discuss the mappings required for the evaluating the inner product of differential forms expressed in terms of a finite element basis which is defined on a reference element. It is remarked that the finite element exterior calculus builds upon affine mappings to the reference element; from a practical standpoint however, it is possible to employ more elaborate mappings, such as those represented by polynomials or rational functions. While we do not provide mathematical analysis of such element mappings, we do not observe any unexpected deterioration of accuracy when using a sufficiently-high number of integration points in the element integrals. This suggests that, at least from a practical standpoint, curvilinear elements do not present a problem. A polynomial mapping of a Nédélec edge basis function is depicted in Fig.~\ref{fig:elementmapping}. We do not make reference to refinement in $h$ or $p$ in this section, as it is not necessary for the definition of the reference element mapping. We close the section by discussing a simple and efficient assembly procedure for higher-order finite elements using orientation-preserving and -reversing mappings to a single reference element.

The reference element mapping is the map from an open region of the computational manifold to an open region of Euclidean space; we thus consider the reference element mapping in the context of charts. We restrict ourselves to a piecewise-smooth approximation to an $n$-dimensional orientable Riemannian manifold $(\Omega, \bs{g})$ with boundary for which an atlas can be defined by choosing suitably many subsets $\{U^\varepsilon_i\}_{i\in A}$ of $\Omega$, including regular and boundary charts. A regular chart is defined such that each $\eta_i : U^\varepsilon_i \rightarrow \refsimplx{T}^{\varepsilon}$ is a diffeomorphism, where $\refsimplx{T}^{\varepsilon}$ denotes the open $\varepsilon$-neighborhood of the reference $m$-simplex $\refsimplx{T}$. We restrict the choice of charts such that when $U^\varepsilon_i \cap U^\varepsilon_j \neq \emptyset$, as $\varepsilon \rightarrow 0$, $\eta_i(U^\varepsilon_i \cap U^\varepsilon_j) \rightarrow \mrm{cl}(f)$ and $\eta_j(U^\varepsilon_i \cap U^\varepsilon_j) \rightarrow \mrm{cl}(g)$ where $f, g$ are potentially different subsimplices of $\refsimplx{T}$.

\subsubsection{Integration over a positively-oriented chart}
Take $(U^\varepsilon_i, \eta'_i)$ as a positively-oriented chart and let $\bar{\eta}'_i = (\eta'_i)^{-1}$ be a map of degree 1, then the inner product of $\omega, \nu \in \Lambda^k(U^\varepsilon_i)$ is given by
\eq{
	\int\limits_{U^\varepsilon_i} \omega \wedge \star_{\bs{g}} \nu = \int\limits_{\eta'_i(U^\varepsilon_i)} (\bar{\eta}'_i)^* \left[ \omega \wedge \star_{\bs{g}} \nu \right] \, ,
}
where $\star_{\bs{g}} : \Lambda^k \rightarrow \Lambda^{m-k}$ is the Hodge star with metric $\bs{g}$. Since the pullback is compatible with the exterior product and the Hodge star commutes with the pullback when also pulling back the metric, we have
\eq{
	(\bar{\eta}'_i)^*\left[ \omega \wedge  \star_{\bs{g}} \nu \right] = (\bar{\eta}'_i)^*  \omega \wedge  \star_{(\bar{\eta}'_i)^*\bs{g}} (\bar{\eta}'_i)^* \nu \, .
}
We wish to define the finite element basis $k$-forms on the reference simplex $\refsimplx{T}$ for reasons of efficient implementation and integration. To this end, we identify $(\bar{\eta}'_i)^* \omega$ with the linear combination $\omega_j \phi_j \in \Lambda^k(\eta'_i(U^\varepsilon_i))$ and $(\bar{\eta}'_i)^* \nu$ with $\nu_j \phi_j \in \Lambda^k(\eta'_i(U^\varepsilon_i))$ such that
\eq{
	(\bar{\eta}'_i)^*\left[ \omega \wedge  \star_{\bs{g}} \nu \right] = \omega_k \phi_k \wedge  \star_{(\bar{\eta}'_i)^*\bs{g}} \phi_l \nu_l \, .
}
For an operator $A$ which commutes with the pullback (such as the exterior derivative $\extd$), integration over the reference element is enabled by the identity
\eq{
	\int\limits_{U^\varepsilon_i} A\omega \wedge \star_{\bs{g}} A\nu = \omega_k\nu_l \int\limits_{U^\varepsilon_i}  A(\eta_i')^*\phi_k  \wedge \star_{\bs{g}}  A(\eta_i')^*\phi_l  = \omega_k\nu_l \int\limits_{\eta'_i(U^\varepsilon_i)}  A\phi_k \wedge  \star_{(\bar{\eta}'_i)^*\bs{g}} A\phi_l \, . \label{eq:refintident}
}
This expression generalizes the Piola-transform maps used in vector elements to differential forms of arbitrary degree $k$ on an $n$-dimensional manifold. For the assembly procedure employed in this work, the coordinate chart map $\eta$ need not be orientation-preserving such that more special attention must be paid to the reference element integral, as outlined in section~\ref{sec:basis_continuity}. To keep computational costs low in implementation, the $p$-hierarchical decomposition of the bases is exploited by numerical integration of the linear and bilinear forms using the $p$-hierarchical Grundmann-Möller rule~\cite{Grundmann1978} on the $n$-simplex.

\subsubsection{Integration over the boundary}
\label{sec:int_over_bnd}

Let $(U^\varepsilon, \eta')$ denote a positively-oriented boundary chart of $\Omega^\varepsilon$ with $\eta' : U^\varepsilon \rightarrow \refsimplx{T}^\varepsilon_F$ where $\refsimplx{T}^\varepsilon_{F}$ denotes the open $\varepsilon$-neighborhood of the reference simplex $\refsimplx{T}$ modified such that $\eta'(U^\varepsilon \cap \partial \Omega) = \refsimplx{F}^\varepsilon$ for a single simplex $\refsimplx{F} \in \partial \refsimplx{T}$. Given a partition of unity $\rho_i$ subordinate to the cover $U_i$ of $\Omega$, then for $\omega \in \Lambda^{m-1}(\Omega)$, Stokes' theorem implies
\eq{
	\sum\limits_{i \in A} \; \int\limits_{\Omega} \extd ( \rho_i\omega )   = \sum\limits_{i \in A'} \int\limits_{\partial \Omega}  \Tr{\partial \Omega }{} ( \rho_i \omega ) \,,
}
where $A'\subseteq A$ indexes boundary charts. We now wish to show that for every boundary chart the local statement
\eq{
	 \int\limits_{U_i^\varepsilon} \extd ( \rho_i\omega )   = \int\limits_{U_i^\varepsilon \cap \partial \Omega}  \Tr{U_i^\varepsilon \cap \partial \Omega}{} ( \rho_i \omega ) \,
}
holds. Regarding the closure $\mrm{cl}(U_i^\varepsilon)$ as a simply-connected manifold with boundary, by Stokes' theorem,
\eq{
	 \int\limits_{\mrm{cl}(U_i^\varepsilon)} \extd ( \rho_i\omega )   = \int\limits_{\partial \mrm{cl}(U_i^\varepsilon)}  \Tr{\partial \mrm{cl}(U_i^\varepsilon)}{} ( \rho_i \omega ) \,
}
holds. Now considering that $\rho_i$ is zero on the boundary of $\mrm{cl}(U_i^\varepsilon)$ everywhere except in regions $U_i^\varepsilon \cap \partial \Omega$ proves the assertion.
With the inclusion map $\imath : \refsimplx{F}^\varepsilon_i \rightarrow \refsimplx{T}^\varepsilon_F$ we obtain
\eq{
	\int\limits_{U^\varepsilon_i \cap \partial \Omega}  \Tr{U^\varepsilon_i \cap \partial \Omega}{} ( \rho_i \omega )
	= \int\limits_{\eta'(U^\varepsilon_i \cap \partial \Omega)}  \Tr{\eta'(U^\varepsilon_i \cap \partial \Omega)}{} (\bar{\eta}'_i)^* ( \rho_i \omega )
	= \int\limits_{\simplx{F}^\varepsilon_i}  \imath ^* (\bar{\eta}'_i)^* ( \rho_i \omega )
}
and therefore Stokes theorem in its local form on the reference simplex $\refsimplx{T}$ is
\eq{
	 \int\limits_{\refsimplx{T}^\varepsilon_F} \extd   (\bar{\eta}'_i)^*(\rho_i\omega)   =    \int\limits_{\refsimplx{F}^\varepsilon_i}   \imath^*  (\bar{\eta}'_i)^* (\rho_i\omega )  \,.
}

\subsubsection{Ensuring continuity of the assembled basis }\label{sec:basis_continuity}
A straightforward assembly procedure is obtained by choosing a globally-consistent ordering of the vertices when mapping to the reference element~\cite{Rognes2009}: for example, if each global vertex is assigned an integer index, then the local vertex ordering can be chosen according to ascending global indices. In addition, according to the construction outlined in section~\ref{sec:geometric_decomposition_simplex}, each subsimplex of the reference element is also ordered according to ascending indices of the reference element. This procedure has the advantage that for two adjacent $m$-simplices $\simplx{T}_1$ and $\simplx{T}_2$, all common subsimplices of $\simplx{T}_1$ and $\simplx{T}_2$ will agree on their local parametrization and therefore their relative orientation and origin, thus avoiding the need for additional basis function vertex permutations or reorientations. As will be demonstrated, such an assembly procedure does not require orientation information for inner-products over the domain volume. However, it does destroy a consistent orientation of the boundary. Hence, we import meshes with vertices numbered so that all simplices are positively-oriented. We then permute to ascending ordering and in the process may flip some positively-oriented simplices to negative orientation. Where this is the case, we store a negative orientation flag which we use for flipping to outward boundary orientation when required.

Suppose the collection of vertices defining a positively-oriented simplex of the mesh is not defined locally by ascending global indices, then the permutation $\pi$ of the vertices required to obtain the correct ordering will render the map $\eta_i$ from $U^\varepsilon_i$ to the positively-oriented reference simplex either orientation-preserving or -reversing depending on whether the parity $\mrm{sgn}(\pi)$ of $\pi$ is even or odd, respectively. To illustrate the relationship of these mappings with respect to inner products, we decompose $\eta_i$ into its orientation-preserving part $\eta'_i : U^\varepsilon_i \rightarrow \refsimplx{T}^{\varepsilon}$, $\mrm{deg}( \eta'_i ) = 1$,  and an automorphism $\eta^\pi_i$, mapping the $\varepsilon$-neighborhood of the reference simplex to the $\varepsilon$-neighborhood of the reference simplex with vertices permuted by $\pi$, such that $\eta_i = \eta^\pi_i \circ \eta'_i$ and $\bar{\eta}_i = \bar{\eta}'_i \circ (\eta^\pi_i)^{-1}$. Let $(\eta^\pi_i)^{-1} = \bar{\eta}^\pi_i$, then the pullback along $\bar{\eta}_i$ is given by $\bar{\eta}_i^* = (\bar{\eta}_i^\pi)^* \circ (\bar{\eta}'_i)^*$ and hence
\eq{
	\int\limits_{U^\varepsilon_i} \omega \wedge \star_{\bs{g}} \nu = \mrm{deg}(\eta_i) \int\limits_{\eta_i(U^\varepsilon_i)} (\bar{\eta}_i^\pi)^*\left[  (\bar{\eta}'_i)^*   \omega \wedge  \star_{(\bar{\eta}'_i)^*\bs{g}} (\bar{\eta}'_i)^* \nu \right]  \, .
}
Since $\eta^\pi_i$ is an automorphism $\mrm{deg}(\eta_i) = \mrm{sgn}\mrm{deg}(\eta_i^\pi) = \mrm{sgn}(\pi)$. Then, for a differential $k$-form, the Hodge star commutes with the pullback along $\bar{\eta}^\pi_i$ up to a change in sign, \ie ${(\bar{\eta}^\pi_i)^* ({\star_\bs{g}\omega})} = \mrm{sgn}(\pi) \,{\star_{(\bar{\eta}^\pi_i)^*\bs{g}} ((\bar{\eta}^\pi_i)^*\omega)}$, so that
\eq{
	\mrm{sgn}(\pi) \int\limits_{\eta_i(U^\varepsilon_i)} (\bar{\eta}_i^\pi)^*\left[  (\bar{\eta}'_i)^*   \omega \wedge  \star_{(\bar{\eta}'_i)^*\bs{g}} (\bar{\eta}'_i)^* \nu \right] = \int\limits_{\eta_i(U^\varepsilon_i)} (\bar{\eta}_i)^*   \omega \wedge  \star_{(\bar{\eta}_i)^*\bs{g}} (\bar{\eta}_i)^* \nu \, .
}
It may therefore be concluded that the orientation-reversal of a given $\eta_i$ must not explicitly be considered when integrating over the reference element $\refsimplx{T}$. The relative orientation of the sub-simplices in the mesh will be reflected in the sign of the associated coefficients. Thus, in the following, we identify terms of the form $(\bar{\eta}_i)^* \nu$ with the expansion $\nu_j \phi_j \in \Lambda^k(\eta_i(U^\varepsilon_i))$. In this case, the identity in Eq.~\ref{eq:refintident} remains valid when $\bar{\eta}'_i$ is replaced with $\bar{\eta}_i$.

\section{P-hierarchical refinement and a posteriori error-indicator of finite element differential forms}

\begin{figure}[bht]
	\centering
	\includegraphics[width=.65\linewidth]{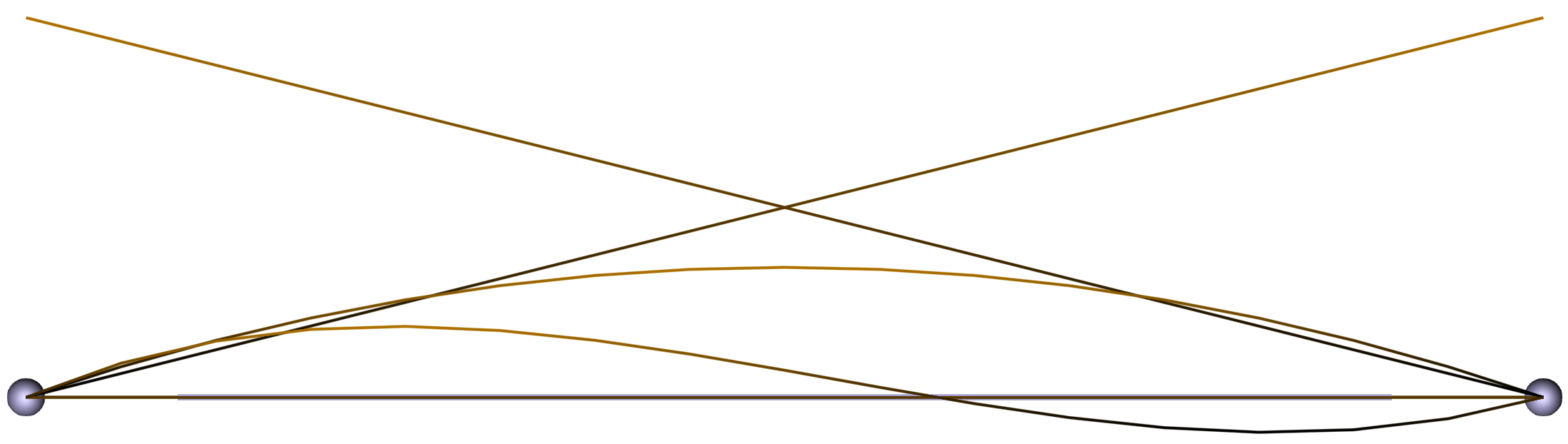}
	\caption{A $p$-hierarchical basis for $\mc{P}_3\Lambda^0$ in one dimension. Vertices are depicted as spheres and the 1-simplex is depicted as a black line.}
	\label{fig:hierarchicalbasis}
\end{figure}

It will now be the purpose to outline the assumptions placed upon a $p$-hierarchical finite element basis in the framework of the finite element exterior calculus. By $p$-hierarchical we mean a basis where polynomial degrees can be modified on subsimplices without having to recompute the basis for the adjacent simplices. Some of the spaces in the three-dimensional de Rham complex have been analyzed in the $p$-hierarchical setting~\cite{Ainsworth2003}, however we are not aware of an entirely general treatment. Aiming for this generality, our construction is entirely local to a reference simplex, greatly simplifying implementation in computer code by retaining only a single reference basis. An example of such a $p$-hierarchical basis is depicted in Fig.~\ref{fig:hierarchicalbasis}. Building upon these assumptions, we outline a simple algorithm for computing a $p$-hierarchical basis on the reference simplex given a valid sequence of non-hierarchical finite element exterior calculus bases of increasing polynomial order, such as those obtained for uniform polynomial order when solving for the basis using the degrees of freedom.

In order to obtain a globally valid polynomial finite element space from such a $p$-hierarchical finite element basis in the case where the polynomial order distribution is non-uniform in space, we recall the well-known minimum-rule. In literature, a basis satisfying the minimum rule is sometimes termed hierarchical~\cite{Licht2017}, although this is not the terminology we adopt in this work. In our terminology, valid finite element basis functions could be chosen on the complex which satisfy the minimum rule but which are not $p$-hierarchical. 

To conclude the section, we introduce an entirely simplex-local a posteriori error indicator based solely upon the examination of coefficients of the solution in a spectral basis with desirable properties. As a by-product of this a posteriori error-indicator, the decay rate of the coefficients in this basis enables us to determine whether to refine in $h$ or in $p$.

\subsection{Properties of finite element differential forms on simplices}

Theorems 4.15 and 4.21 of ~\cite{Arnold2006} establish that the spaces $\PLminus{r}{k}{T}$ and $\PLplus{r}{k}{T}$ admit a geometric decomposition into mutually disjoint subspaces $\fcspcminus{r}{k}{T}{f}{}$ and $\fcspcplus{r}{k}{T}{f}{}$ associated to $f \in \Delta(\simplx{T})$, respectively, such that for $k,r \ge 1$,
\eq{
\PLminus{r}{k}{T} = \bigoplus\limits_{f\in \Delta(T)} \fcspcminus{r}{k}{T}{f}{} \,,\quad \quad \quad
\PLplus{r}{k}{T} = \bigoplus\limits_{f\in \Delta(T)} \fcspcplus{r}{k}{T}{f}{} \,,
}
where $\fcspcminus{r}{k}{T}{f}{} = \fcspcplus{r}{k}{T}{f}{} = 0$ if $r+k \le \mrm{dim} f$ or $\mrm{dim} f  < k$. The elements of $\fcspcminus{r}{k}{T}{f}{}$ and $\fcspcplus{r}{k}{T}{f}{}$ are nonzero on all $\ssimplx{g}{\mu} \in \Delta\simplx{T}$ where $\ssimplx{f}{\sigma} \subseteq \tx{cl}(\ssimplx{g}{\mu})$ and vanish in the sense of traces where $\ssimplx{f}{\sigma} \not\subseteq \tx{cl}(\ssimplx{g}{\mu})$ respectively, establishing their limited support and therefore the geometric locality of the basis.

The inhomogeneous refinement and de-refinement of polynomial orders within the assembled finite element space is made possible by ensuring that such operations remain local. Transitions between areas of differing polynomial orders require the introduction of elements where each subsimplex in $\Delta(\simplx{T})$ may be assigned a polynomial order $r$ as required to retain a valid polynomial space. The construction of spaces which offer such choice is enabled by a $p$-hierarchical decomposition of the spaces $\fcspcminus{r}{k}{T}{f}{}$ and $\fcspcplus{r}{k}{T}{f}{}$ associated to each $f \in \Delta(\simplx{T})$. For this purpose, it is assumed that there exist extension operators $\extensminus{r,k}{f}{}{T}$ and $\extensplus{r,k}{f}{}{T}$,
such that
\eq{
\extensminus{r,k}{f}{}{T} \trspcminus{r}{k}{f}{} &\mapsto \fcspcminus{r}{k}{T}{f}{}  \,, \quad \quad \Tr{f}{} \extensminus{r,k}{f}{}{T} \trspcminus{r}{k}{f}{} \mapsto \trspcminus{r}{k}{f}{} \label{eq:extensminus} \,,\\
\extensplus{r,k}{f}{}{T} \trspcplus{r}{k}{f}{} &\mapsto \fcspcplus{r}{k}{T}{f}{} \,, \quad \quad \Tr{f}{} \extensplus{r,k}{f}{}{T} \trspcplus{r}{k}{f}{} \mapsto \trspcplus{r}{k}{f}{} \, . \label{eq:extensplus}
}
The spaces $\fcspcminus{r}{k}{T}{f}{}$ and $\fcspcplus{r}{k}{T}{f}{}$ contain polynomial differential forms of degree $k$ and order at least $p$, where $p \ge 1$, as well as of at most $q$, where $q \le r$, while $p$ and $q$ depend on the dimension of $f$. From the nesting relation $\PLminus{r}{k}{T} \subseteq \PLplus{r}{k}{T} \subseteq \PLminus{r+1}{k}{T}$, we obtain the facewise nesting: $\fcspcminus{r}{k}{T}{f}{} \subseteq \fcspcplus{r}{k}{T}{f}{}$ for $p \le r \le q$ and $\fcspcplus{r}{k}{T}{f}{} \subseteq \fcspcminus{r+1}{k}{T}{f}{}$ for $p \le r$ and $r+1 \le q$.

Returning to the goal of establishing properties of a hierarchical basis, the preferred hierarchical decomposition of the face spaces is constructed by first defining what the present work considers as being the most granular component of the facewise nesting, the transition spaces
\eq{
\Delta W^{k\pm}_r(f) =  \trspcplus{r}{k}{f}{}  \setminus \trspcminus{r}{k}{f}{} \,,\quad & \quad \Delta W^{k\mp}_{r}(f) =  W^{k-}_{r}(f)  \setminus W^{k+}_{r-1}(f) \,,\quad \, .
}
This allows for the definition of the hierarchical difference spaces as
\eq{
\Delta \trspcminus{r}{k}{f}{} =     \Delta W^{k\pm}_{r-1}(f) \oplus \Delta W^{k\mp}_{r}(f)  &\quad
\Delta  \trspcplus{r}{k}{f}{}=     \Delta W^{k\mp}_{r}(f) \oplus \Delta W^{k\pm}_{r}(f)
\, ,}
directly giving rise to the hierarchical decompositions
\eq{
\trspcminus{r}{k}{f}{} =    \trspcminus{1}{k}{f}{} \oplus  \Delta W^{k-}_{2}(f) \oplus \ldots \oplus \Delta W^{k-}_{r}(f) \,, \\
\trspcplus{r}{k}{f}{} =    \trspcminus{1}{k}{f}{} \oplus  \Delta W^{k\pm}_{1}(f) \oplus  \Delta W^{k+}_{2}(f) \oplus \ldots \oplus \Delta W^{k+}_{r}(f)  \,.
}
A proper extension from subsimplex $f$ to the simplex $\simplx{T}$ is straightforward as shown in Eqs.~\ref{eq:extensminus} and \ref{eq:extensplus}. Owing to the preferred structure of the hierarchical basis, we require that the extension operators $E^{k-}_{\ssimplx{f}{},\simplx{T}}$ and $E^{k-}_{\ssimplx{f}{},\simplx{T}}$ be defined as
\eq{
E^{k-}_{\ssimplx{f}{},\simplx{T}}\trspcminus{r}{k}{f}{} \mapsto   &\extensminus{1,k}{f}{}{T} \trspcminus{1}{k}{f}{} \oplus  \extensminus{2,k}{f}{}{T}\Delta W^{k-}_{2}(f) \oplus \ldots \oplus \extensminus{r,k}{f}{}{T}\Delta W^{k-}_{r}(f)  \,,\label{eq:hierarchicalextensionminus}\\
E^{k+}_{\ssimplx{f}{},\simplx{T}}\trspcplus{r}{k}{f}{} \mapsto    &\extensminus{1,k}{f}{}{T}\trspcminus{1}{k}{f}{} \oplus  \extensplus{1,k}{f}{}{T}\Delta W^{k\pm}_{1}(f) \notag\\
&\oplus  \extensplus{2,k}{f}{}{T}\Delta W^{k+}_{2}(f) \oplus \ldots \oplus \extensplus{r,k}{f}{}{T}\Delta W^{k+}_{r}(f) \,. \label{eq:hierarchicalextensionplus}
}
The analogue which defines the corresponding face spaces is
\eq{
\fcspcminus{r}{k}{T}{f}{} &=   \fcspcminus{1}{k}{T}{f}{} \oplus  \Delta \fcspcminus{2}{k}{T}{f}{} \oplus \ldots \oplus \Delta \fcspcminus{r}{k}{T}{f}{} \label{eq:hierarchicaldecompositionminus}\\
\fcspcplus{r}{k}{T}{f}{} &=    \fcspcminus{1}{k}{T}{f}{} \oplus  \Delta V^{k\pm}_{1}(\simplx{T}, f) \oplus  \Delta\fcspcplus{2}{k}{T}{f}{} \oplus \ldots \oplus \Delta\fcspcplus{r}{k}{T}{f}{} \,. \label{eq:hierarchicaldecompositionplus}
}
It is noted that the Whitney forms contained in $\fcspcminus{1}{k}{T}{f}{}$ are explicitly present in both polynomial order decompositions, a property that we require in establishing a hierarchy in discretization size as will be outlined in Section~\ref{sec:refeq}. For multi-field problems, it is necessary to ensure that no matter what polynomial order is requested of a subsimplex, a stable pairing over adjacent simplices is always present. This requires some coordination effort for an implementation, which however is entirely local to the subsimplex for which the polynomial order change is requested.

\subsection{An algorithm for computing a p-hierarchical basis}

An algorithm for computing a hierarchical basis from a valid sequence of given sets $\{\fcspcminus{q}{k}{T}{f}{}\}_{q=1}^r$ and $\{\fcspcplus{q}{k}{T}{f}{}\}_{q=1}^r$ decomposing finite element bases for $\PLminus{r}{k}{T}$ and $\PLplus{r}{k}{T}$ which respects~\ref{eq:hierarchicalextensionminus} and~\ref{eq:hierarchicalextensionplus} is then straightforward. For ${\mc{P}^-}\Lambda^{k}(\simplx{T})$, taking the trace of $\fcspcminus{q}{k}{T}{f}{}$ and $\fcspcminus{q-1}{k}{T}{f}{}$ determines $\trspcminus{q}{k}{f}{}$ and $\trspcminus{q-1}{k}{f}{}$. Compute $\Delta \trspcminus{q}{k}{f}{} = \trspcminus{q}{k}{f}{} \setminus \trspcminus{q-1}{k}{f}{}$ and note which elements are removed. Remove the corresponding elements of  $\fcspcminus{q}{k}{T}{f}{}$ to determine $\Delta\fcspcminus{q}{k}{T}{f}{}$. For $\PLplus{}{k}{T}$, satisfying equation~\ref{eq:hierarchicalextensionplus} requires a modification for $q = 1$ while otherwise using an analogous procedure. To acquire $\fcspcplus{1}{k}{T}{f}{}$ consistently, take the trace of $\fcspcplus{1}{k}{T}{f}{}$ and $\fcspcminus{1}{k}{T}{f}{}$ to determine $\trspcplus{1}{k}{f}{}$ and $\trspcminus{1}{k}{f}{}$. Compute $ \Delta W^{k\pm}_1(f) =  \trspcplus{1}{k}{f}{}\setminus \trspcminus{1}{k}{f}{}$ and note which elements are removed. Remove the corresponding elements in $\fcspcplus{1}{k}{T}{f}{}$ to determine $\Delta V^{k\pm}_1(f)$ and $\fcspcplus{1}{k}{T}{f}{} = \fcspcminus{1}{k}{T}{f}{} \oplus \Delta V^{k\pm}_1(f)$.

We use precisely this procedure to transform an arbitrarily generated finite element exterior calculus basis to its $p$-hierarchical form. In order to determine linear dependence, we employ GMRES orthogonalization. All the procedures in this algorithm are carried out in arbitrary-precision arithmetic to ensure its accuracy in the presence of numerical round-off. Following~\cite{Ainsworth2003}, in order to improve conditioning of the generated basis, we orthogonalize the volume bubble functions per polynomial order with respect to an appropriate inner product using GMRES. In our numerical studies, we do not observe ill-conditioned stiffness matrices.

\subsection{The minimum rule}
For $p$-adaptivity without $h$-refinement, let $\Omega_d$ represent the set of top-dimensional simplices in a given simplicial representation of $\Omega \subset \mathbb{R}^n$ with polynomial orders $r > 0$ assigned by the surjection $\mathfrak{O} : \Omega_d \rightarrow \mathbb{N}$. The minimum rule then states that to each subsimplex $\ssimplx{f}{} \in \Dksimplx{T}{k}$ of dimension $0 < k \le n$, the maximum allowed polynomial order assignable is determined by the map
\begin{align}
	 f \mapsto \min \left\{\mathfrak{O}(\simplx{T}) \, : \, \ssimplx{f}{} \subset \simplx{T}, \; \simplx{T} \in \Omega_d \right\} \, ,
	 \label{eq:minrule}
\end{align}
while subsimplices of dimension zero are assigned the order one. For a general choice of basis the minimum rule is required as, given a polynomial order distribution of the $n$-simplices, it is the only distribution over a set of $k < n$-subsimplices for which the approximation of a given function is invariant under differing choices of shape functions associated to them~\cite{solin2003higher}.

\subsection{A posteriori error and spectral decay indicator}

Inspired by the error indicator developed by Mavriplis (1989)~\cite{Mavriplis1989}, we develop an a posteriori error indicator by projection of solution polynomials onto eigenfunctions of the PDE problem
\eq{
	\mathcal{L} \phi = \omega \phi
}
on the $n$-simplex $\simplx{T}$.  To obtain symmetric orthogonal polynomials as eigenfunctions on the simplex, we employ a Legendre-like differential operator given by Braess and Schwab (2000)~\cite{Braess2000}, which in barycentric coordinates takes the form
\eq{
	\mathcal{L} = \sum\limits_{\ssimplx{f}{\sigma} \in \Dksimplx{T}{1}}  \lambda_{\sigma(1)} \lambda_{\sigma(2)} \left(\pd{}{\lambda_{\sigma(2)}} - \pd{}{\lambda_{\sigma(1)}}\right)^2 \, .
}
For cell complexes composed of non-simplicial cells, only the eigenvalue problem would be replaced. In the simplicial case, the generalized eigenvalue problem has the distinct eigenvalues $\omega_p = p(p+n)$ for the polynomial orders $p = 0, 1, 2, \ldots, r$. The eigenfunctions $\phi_i$ then fall in invariant subspaces $\Phi_p$ associated to these eigenvalues. Given an $n$-variate solution polynomial $a \in \mc{P}_r(\refsimplx{T})$ of order $r$ defined over the reference $n$-simplex with coordinates $\bs\xi$, the spectral coefficients $a_{pi}$ of $a$ may be obtained such that
\eq{
	a(\bs\xi) = \sum\limits_{p = 0, \ldots, r} \sum\limits_{i = 1,\ldots, |\Phi_p|} a_{pi} \phi_{pi}(\bs\xi) \, .
}
A simplified coefficient sequence in $p$ can be found by averaging over elements of the invariant subspaces:
\eq{
	a_p = \frac{1}{|\Phi_p|}\sum\limits_{i = 1,\ldots, |\Phi_p|}a_{pi} \, , \quad p = 0,\ldots,r \, .
}
In analogy to the error indicator developed by Mavriplis (1989)~\cite{Mavriplis1989}, we perform a fitting of the sequence $\{a_p\}_{p=0}^r$ to an exponential decay of the form
\eq{
	\hat{a}(p) = c \exp(-\sigma p) \, .
}
In the present case, for reasons of robustness, the sequence $\{a_p\}_{p=1}^r$ is fitted to the exponential decay by minimizing in the $L^1$-norm. As is standard procedure, we use the iteratively reweighted least squares (IRLS) method to carry out the minimization. 

Finally, we give a closed-form solution to the continuation problem given in~\cite{Mavriplis1989}, leading to the error indicator
\eq{
	\varepsilon(c, \sigma) = \sqrt{\frac{2 a_r^2}{2r+1} + c^2 \exp(\sigma) \, E_1((2r+3)\sigma)} \,
}
for the $L^2$-norm of the polynomial solution, where $E_1(z)$ is the generalized exponential integral. The exponent $\sigma$ gives an indication as to the rate at which the polynomial $a$ decays in the spectrum of $\mathcal{L}$, a vital tool in the selection of appropriate polynomial order distributions and in singularity detection.

For solution polynomial differential forms $\bs a \in \mc{P}_r\Lambda^k(\refsimplx{T})$, defined in terms of the reference element we require a definition of the differential form basis globally. We therefore apply the pushforward to the basis wedges
\eq{
a_{\bs\sigma}(\bs\xi) \, \cobasis{\xi}{\bs\sigma} \mapsto a_{\bs\sigma}(\bs\xi) \; (\varphi \circ \bar{\eta}_i)_*(\cobasis{\xi}{\bs\sigma})
}
and define $a^*_{\bs\nu} \in \mc{P}_r(\refsimplx{T})$ such that
\eq{
a_{\bs\sigma}(\bs\xi) \; (\varphi \circ \bar{\eta}_i)_*(\cobasis{\xi}{\bs\sigma}) = a^*_{\bs\nu}(\bs\xi) \, \cobasis{x}{\bs\nu} \, .
}
In order to ensure that $a^*_{\bs\nu}$ remains polynomial for general element mappings, the tangent map at the simplex centroid can be used in an approximation to the push-forward. It is then possible to apply the error and decay indicator outlined above to $a^*_{\bs\nu}$, for each $\bs{\nu}$ separately, aggregating the result appropriately into a single error-decay pair for the whole polynomial differential form. For vector-valued differential forms, each component of the vector of differential forms is assigned an error-decay pair. Aggregation can also be done here as is deemed appropriate.

\section{Hierarchical simplicial subdivision and hp-hierarchical basis of finite element differerential forms}

We will now combine the $p$-hierarchical basis outlined in the previous section with hierarchical $h$-refinement by building upon the works~\cite{Yserentant1985} and~\cite{Grinspun2003}. Apart from $hp$-adaptivity, the hierarchical basis can also be harnessed in gaining further advantages including multigrid convergence for the solution of linear systems~\cite{Yserentant1985}. Hierarchical methods were also considered briefly in~\cite{Stein2005} and for specific cases were applied in~\cite{Zander2015}. We  generalize these ideas and apply them to the case at hand. This will be made possible by introducing a regular subdivision of the simplex and extending this definition to an incomplete subdivision, \ie one which propagates a refinement of a subsimplex to adjacent simplices. By its design, this construction ensures that no so-called ``hanging'' mesh entities are created in the process. We then state the well-known refinement~\cite{Grinspun2003} or filtering~\cite{Stein2005} equation and adapt it to this construction. By a locality assumption, we resolve the linear dependence by the introduction of a linear constraint on the Whitney forms associated to the refined subsimplex. To ensure invariance of the solution under arbitrary choices of $p$-hierarchical bases, we extend the classical minimum-rule, as was outlined in the previous section, to the case of hierarchical $h$-refinement. While we restrict ourselves to simplices, the outlined methodology also applies in principle to general polyhedra.

\subsection{Incomplete Freudenthal subdivision}

The complete Freudenthal subdivision~\cite{Freudenthal1942, Edelsbrunner2000} $\freud{\simplx{T}}$ divides $\simplx{T}$ into $2^n$ smaller simplices  $\{ \simplx{U}_i\}_{i=1}^{2^n}$ with the same $n$-dimensional volume such that
\begin{align}
	\simplx{T} = \bigcup\limits_{i=1}^{2^n} \simplx{U}_{i} \, .
\end{align}
Following section \ref{sec:geometric_decomposition_simplex}, we define the set of all $k$-dimensional subsimplices $\Dkfreud{\simplx{T}}{k}$ and the set of all subsimplices $\Dfreud{\simplx{T}}$ contained in the complete Freudenthal subdivision $\freud{\simplx{T}}$, as
\begin{align}
	\Dkfreud{\simplx{T}}{k} = \bigcup\limits_{i=1}^{2^n} \Delta_k \simplx{U}_{i} \quad \tx{and} \quad \Dfreud{\simplx{T}} = \bigcup\limits_{k=0}^{n} \Delta_k \mc{F}(\simplx{T})\,,
\end{align}
consequently. Furthermore, a \textit{hierarchical parent-child-relationchip} between a parent-subsimplex $\ssimplx{f}{} \in \Dsimplx{T}$ and the $2^{\mrm{dim}\ssimplx{f}{}}$ children-subsimplices
\begin{align}
\hchildren{\ssimplx{f}{}} = \left\{ t : t \subset \ssimplx{f}{}{} , t\in \Dkfreud{\clssimplx{f}{}}{\mrm{dim}f} \right\} \,
\end{align}
is established. By the definition of the hierarchical parent-child relationchip, recursive refinements of a {\ital root}-subsimplex $\ssimplx{f}{}$ (one without a parent) can be tracked by a $2^{\text{dim}(\ssimplx{f}{})}$-ary tree. This hierarchy is distinct from the hierarchy of geometric ancestors and descendants. We will thus refer to the hierarchical analogues as {\ital hierarchical ancestors and descendants} if not already clear from context.
The subsimplex $\ssimplx{f}{}$ is termed a {\ital leaf}, if $\hchildren{\ssimplx{f}{}} = \emptyset$, identified by the indicator function
\begin{align}
\leaf{\ssimplx{f}{}} = \begin{cases} 1,\quad \text{if}\,\,\hchildren{\ssimplx{f}{}} = \emptyset,  \\ 0, \quad\text{otherwise.} \end{cases}
\end{align}

\newcommand\IMAGESWIDTHIF{.3}

\begin{figure}[htb]
	\centering
	\begin{minipage}[t]{\IMAGESWIDTHIF\linewidth}
        		\begin{tikzpicture}
			\node[inner sep=0pt]{\includegraphics[trim={29cm 11cm 30cm 3cm},clip,width=\linewidth]{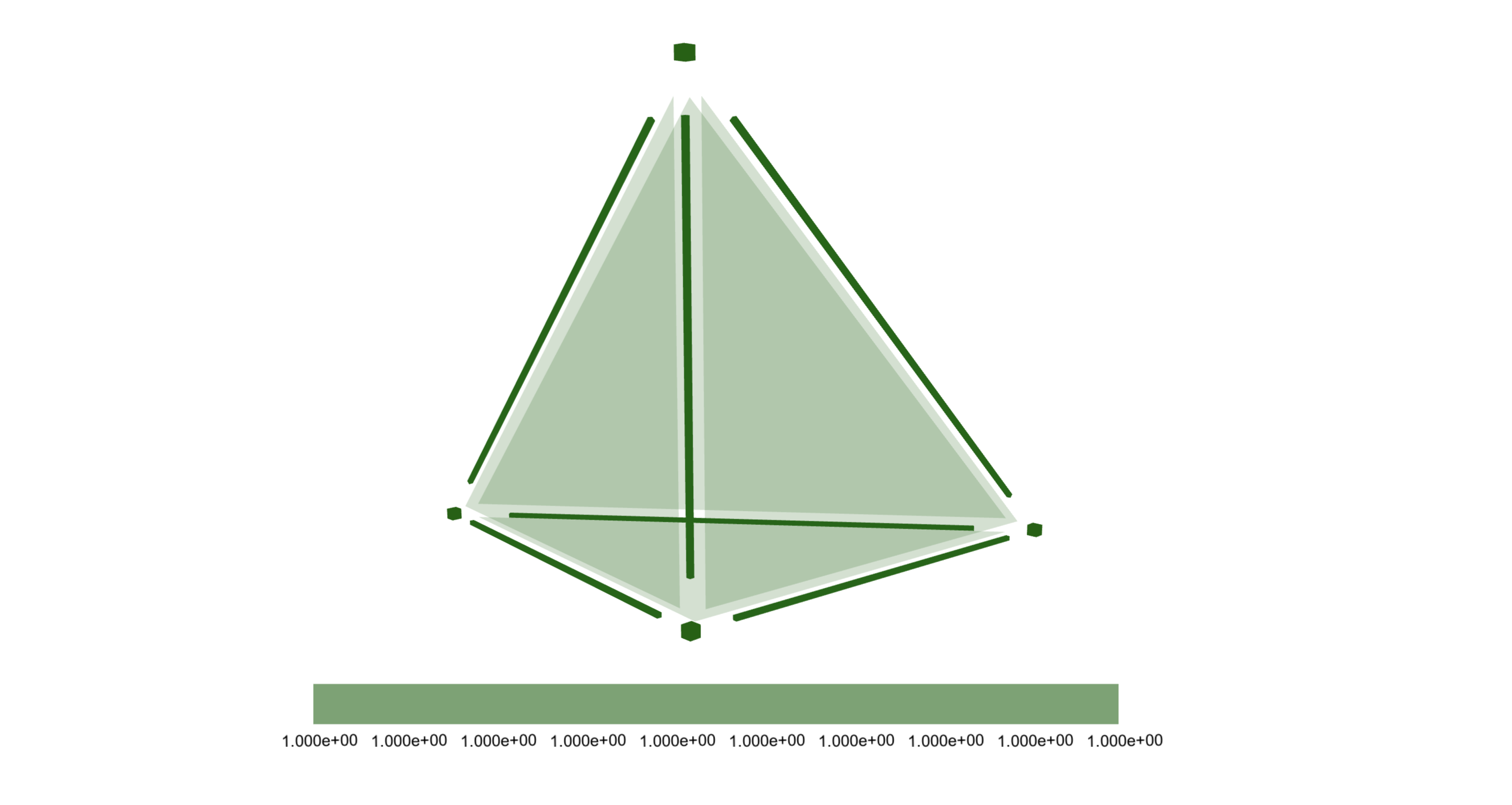}};
	 	\node at (0.025,-1.75) {$f_{(1)}$};
		\node at (2.3,-1.75) {$f_{(2)}$};
		\node at (-2.1,-1.7) {$f_{(3)}$};
		\node at (0.0,2.025) {$f_{(4)}$};
		\end{tikzpicture}
		\begin{center}
        		(a) unrefined subsimplices $\Dsimplx{T}$
		\end{center}
    	\end{minipage}%
    	\hfill
    	\begin{minipage}[t]{\IMAGESWIDTHIF\linewidth}
		\begin{tikzpicture}
			\node[inner sep=0pt]{\includegraphics[trim={29cm 11cm 30cm 3cm},clip,width=\linewidth]{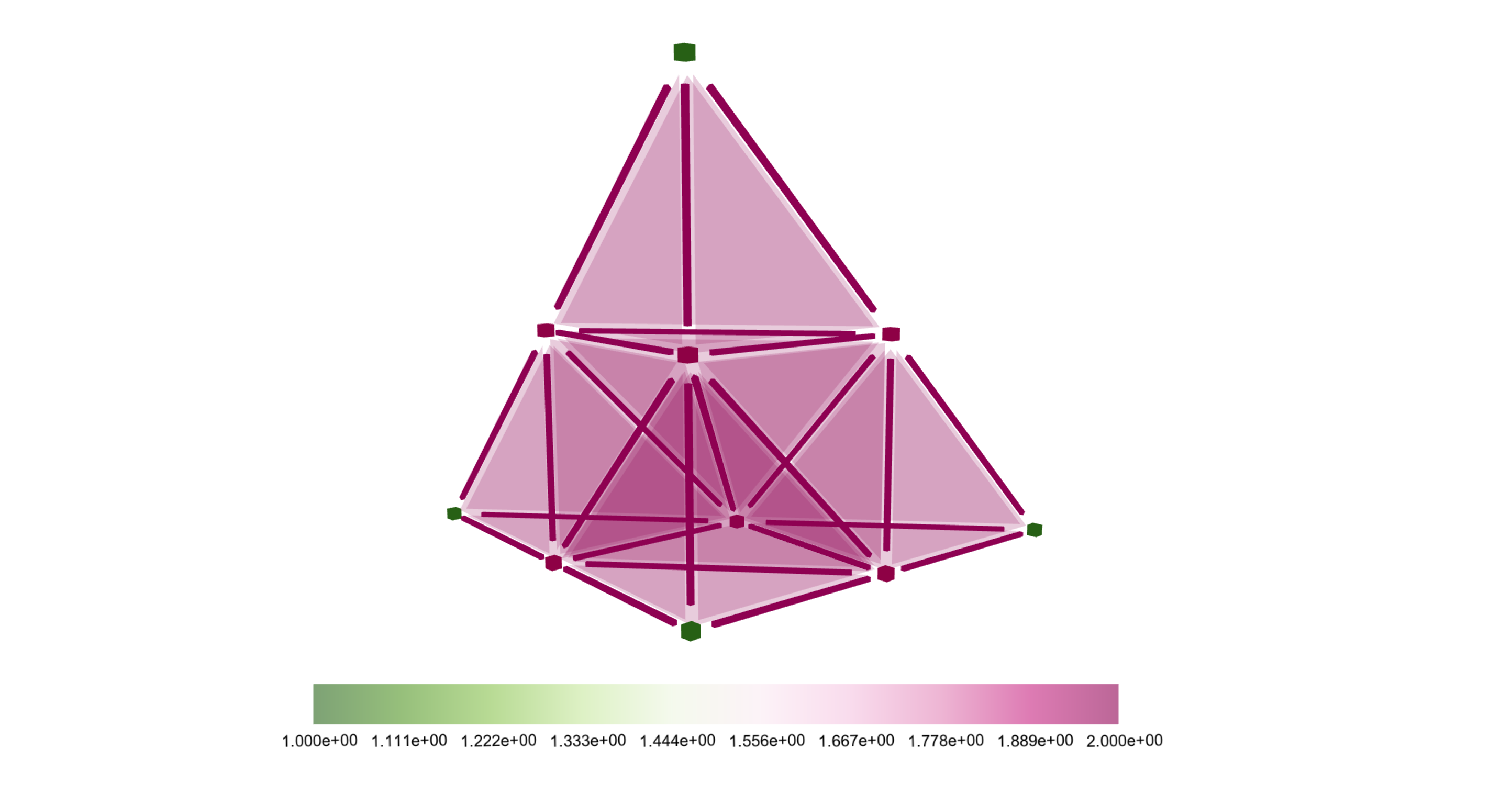}};
		\end{tikzpicture}
        		\begin{center}
        		(b) complete Freudenthal $\Dkfreud{\simplx{T}}{}$
		\end{center}
	\end{minipage}%
	\hfill
    	\begin{minipage}[t]{\IMAGESWIDTHIF\linewidth}
		\begin{tikzpicture}
			\node[inner sep=0pt]{\includegraphics[trim={29cm 11cm 30cm 1cm},clip,width=\linewidth]{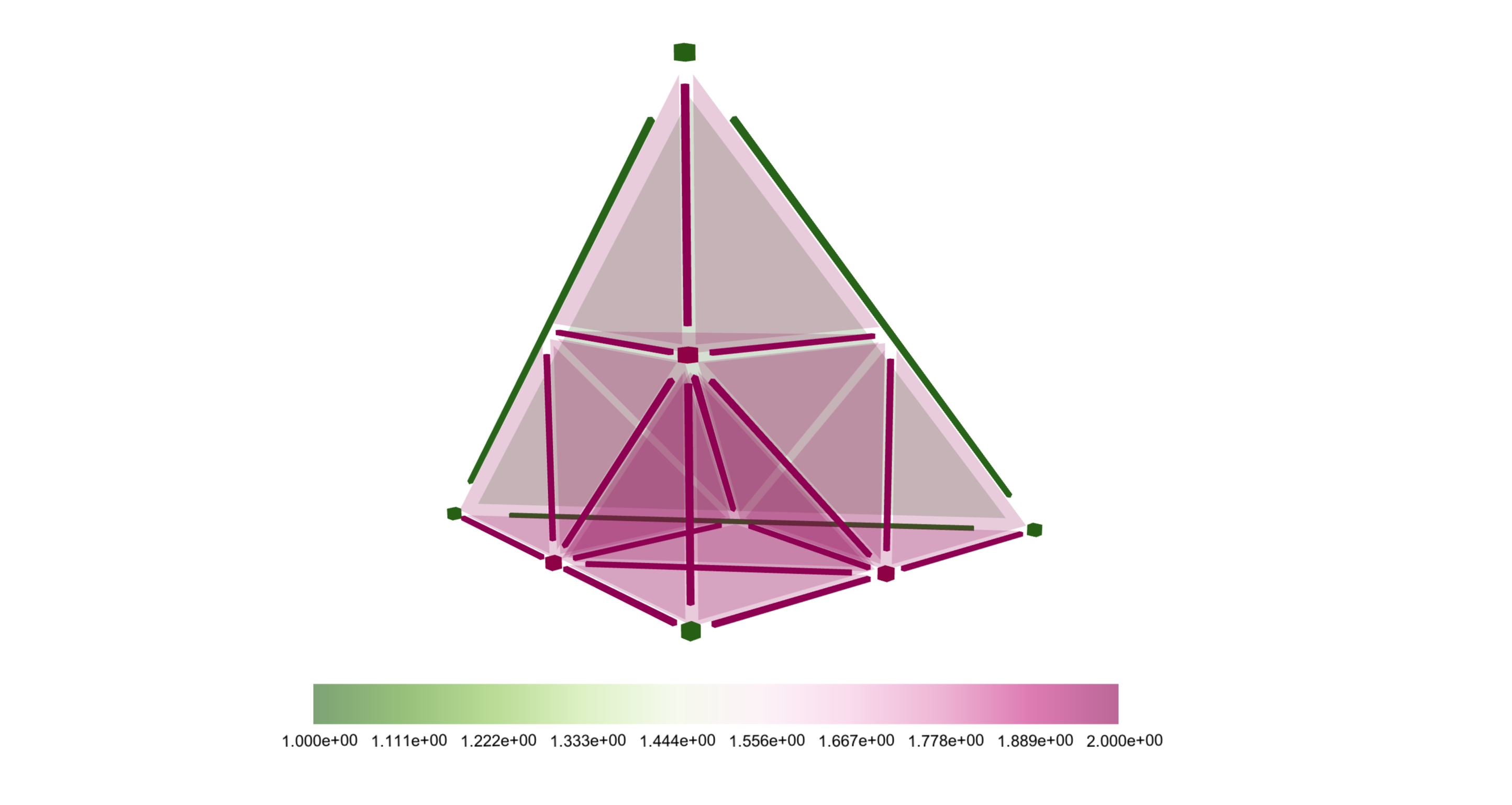}};
		\end{tikzpicture}
        		\begin{center}
        		(c) vertex refinement $\incfreud{\simplx{T},f_{(1)}}$
		\end{center}
	\end{minipage}%
	\hfill
    	\begin{minipage}[t]{\IMAGESWIDTHIF\linewidth}
		\begin{tikzpicture}
			\node[inner sep=0pt]{\includegraphics[trim={29cm 11cm 30cm 1cm},clip,width=\linewidth]{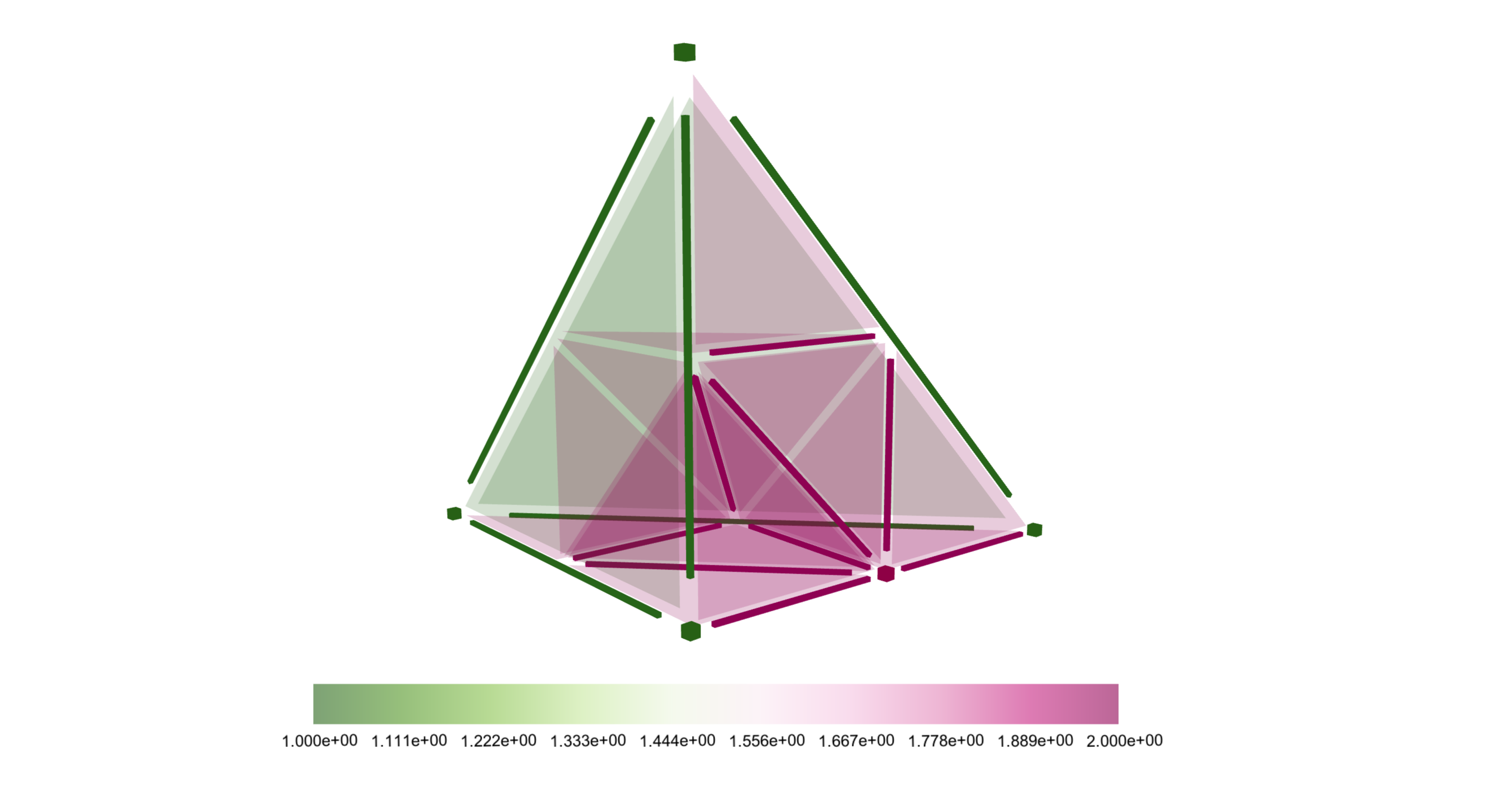}};
		\end{tikzpicture}
        		\begin{center}
        		(d) line refinement $\incfreud{\simplx{T},f_{(1,2)}}$
		\end{center}
	\end{minipage}%
	\hfill
    	\begin{minipage}[t]{\IMAGESWIDTHIF\linewidth}
		\begin{tikzpicture}
			\node[inner sep=0pt]{\includegraphics[trim={29cm 11cm 30cm 1cm},clip,width=\linewidth]{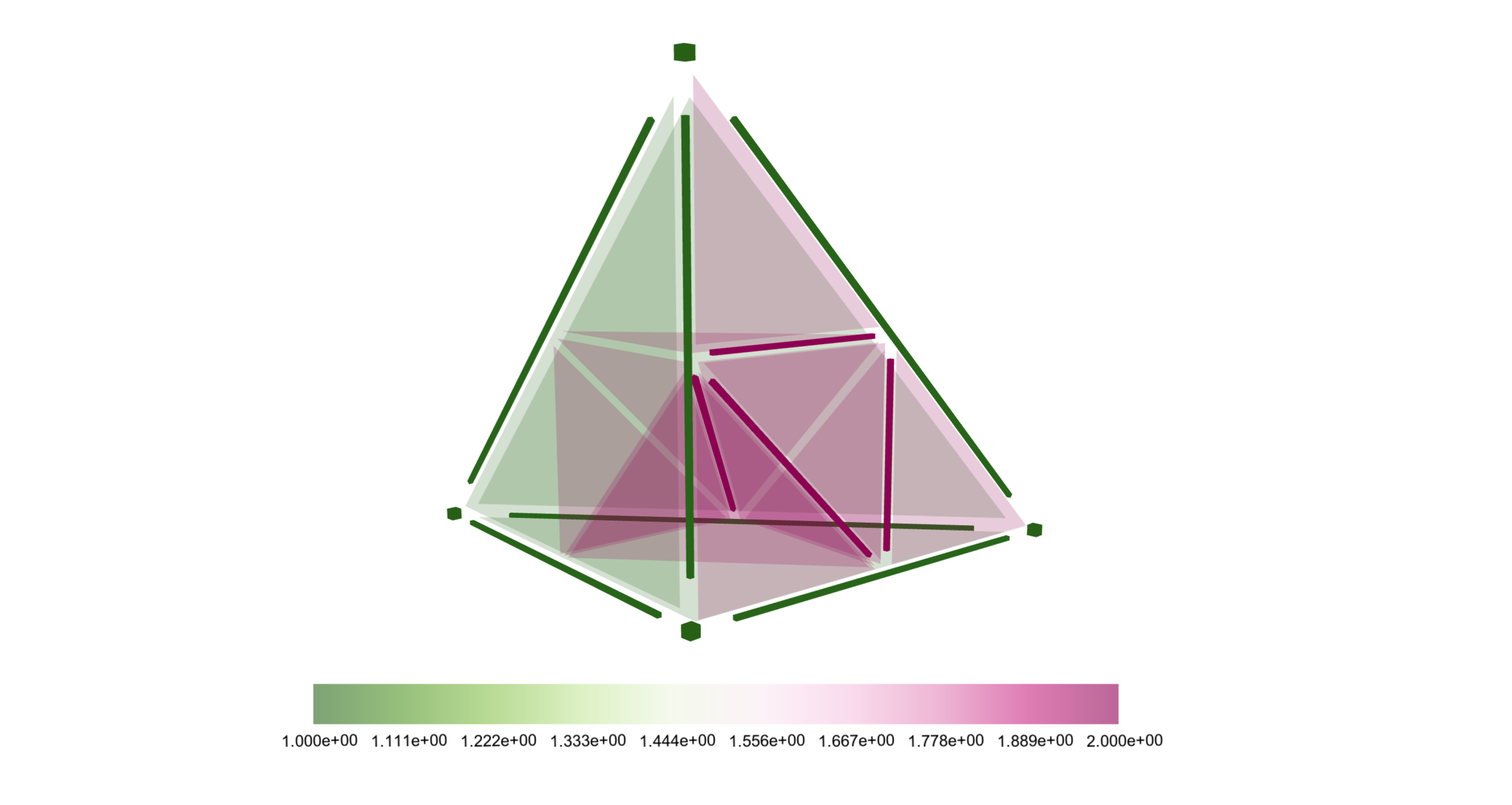}};
		\end{tikzpicture}
        		\begin{center}
        		(e) face refinement $\incfreud{\simplx{T},f_{(1,2,4)}}$
		\end{center}
	\end{minipage}%
	\hfill
    	\begin{minipage}[t]{\IMAGESWIDTHIF\linewidth}
		\begin{tikzpicture}
			\node[inner sep=0pt]{\includegraphics[trim={29cm 11cm 30cm 1cm},clip,width=\linewidth]{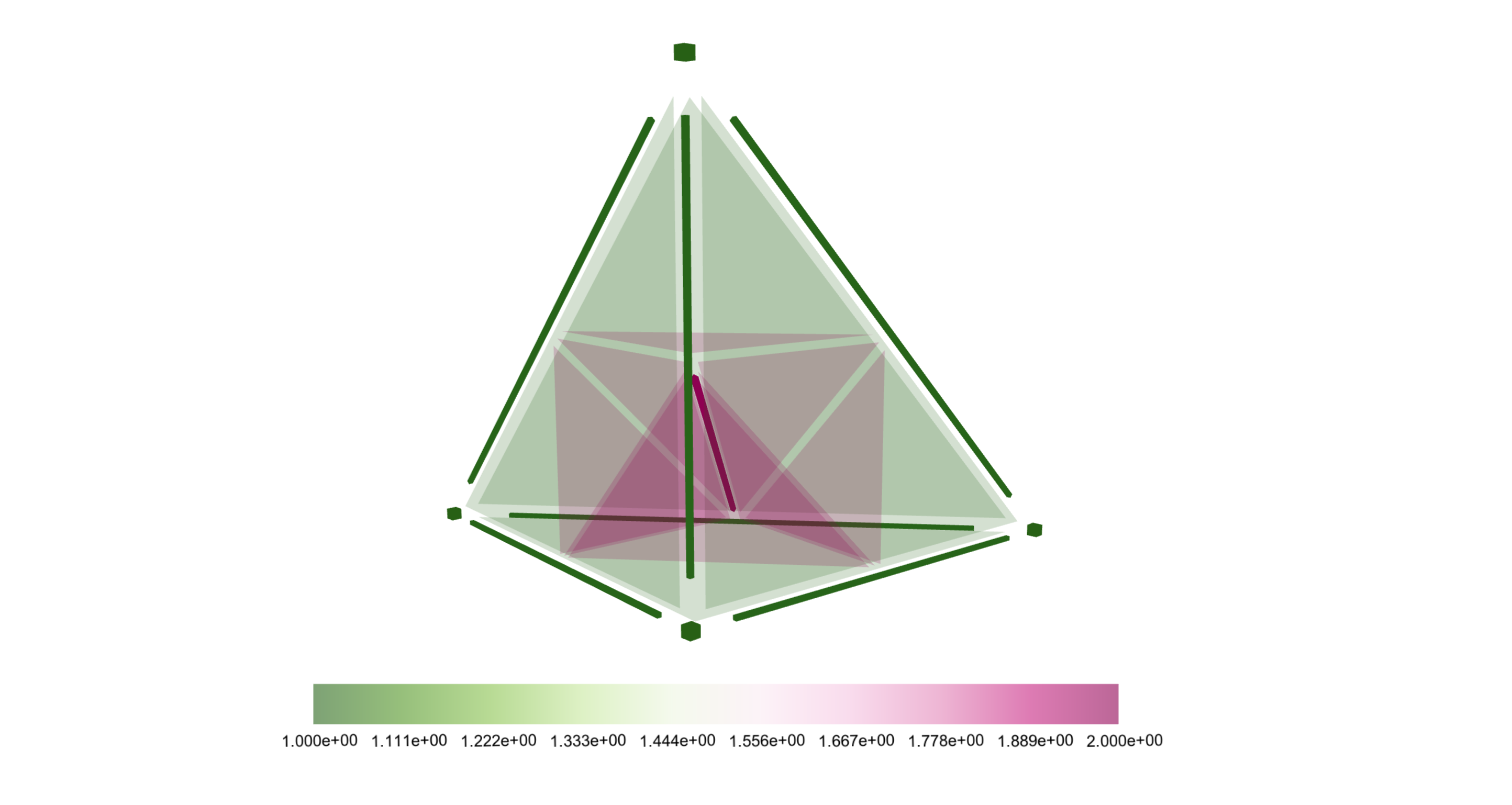}};
		\end{tikzpicture}
        		\begin{center}
        		(f) volume refinement $\incfreud{\simplx{T},f_{(1,2,3,4)}}$
		\end{center}
	\end{minipage}%
	\caption{Leaf vertices, edges and faces of the incomplete Freudenthal subdivision of a tetrahedron $\protect\simplx{T}$. The volumes are omitted in the graphics. They exist on the root level in (a) as well as on the refined level in the remaining depictions.}
\end{figure}
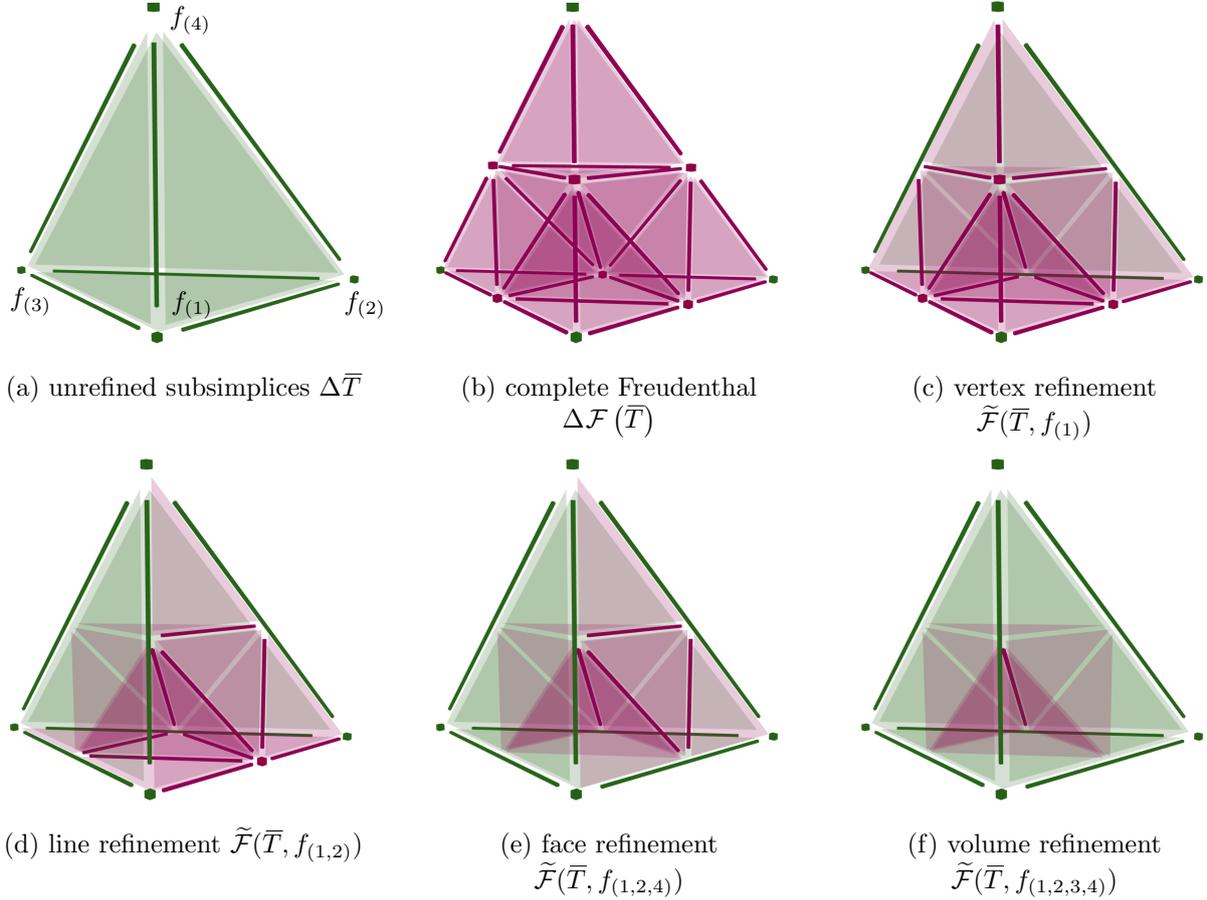

An incomplete Freudenthal subdivision is the set of subsimplices
\begin{align}
		\incfreud{\simplx{T},f} = \left\{ s \in \Delta \mc{F}(\simplx{T}) : s \subset \bigcup\limits_{r \in \upadjacstar{f}	} r \right\} \,.
\end{align}
referring to the subdivision of a single subsimplex $f \in \Dsimplx{T}$ and its compatible extension to $\simplx{T}$ such that the subdivision of all subsimplices of $\simplx{T}$ results in the complete Freudenthal subdivision of $\simplx{T}$.
Let $\Dkincfreud{\simplx{T},f}{k}$ denote the set of all $k$-subsimplices of the incomplete Freudenthal subdivision $\incfreud{\simplx{T},f}$.

In implementation, vertices in the refined complex carry labels such as $((v_0, \ldots, v_k), (m_0, \ldots, m_k))$ for a vertex lying on a $k$-dimensional root subsimplex with integer vertices $(v_1, \ldots, v_k)$. Such a labelling of points on a simplex can be traced to~\cite{Leung1979} and~\cite{Olmsted1986}. The multiplicity $(m_0, \ldots, m_k)$ corresponds to the location of the vertex in barycentric coordinates of the root subsimplex, \ie in its barycentric coordinates, the vertex is located at $\lambda_0 = m_0/(m_0 + \ldots + m_k)$,  $\lambda_1 = m_1/(m_0 + \ldots + m_k)$, , $\lambda_k = m_k/(m_0 + \ldots + m_k)$. In order to suppress duplicate vertices, the multipliers $(m_0, \ldots, m_k)$ may be normalized by $\mrm{gcd}(m_0, \ldots, m_k)$. Determining the hierarchical refinement depth of a given vertex then simply corresponds to computing
\eq{
	\mrm{log}_2\sum\limits_{i=0}^k m_i + 1 \, .
}

\subsection{The refinement equation}
\label{sec:refeq}

For hierarchical $h$-refinement, it is necessary to determine an equation relating coefficients across refinement scales, termed the refinement equation~\cite{Grinspun2003}. While elaborate higher-order or more global constructions are possible, giving rise to second-generation wavelets~\cite{Sweldens1998}, the present work will limit itself to an refinement equation built upon the Whitney forms $\fcspcminus{1}{k}{T}{f}{}$ which is locally restricted to the $k$-subsimplex $f$ being refined. Higher orders, \ie complements of $\fcspcminus{1}{k}{T}{f}{}$ in the sense of the decompositions given in Eqs.~\ref{eq:hierarchicaldecompositionminus} and~\ref{eq:hierarchicaldecompositionplus}, are simply treated separately by assuming that they reside only on the leaf subsimplices.

We thus begin by treating one level of refinement due to the splitting of a single $k$-subsimplex $f$ of a root simplex $\simplx{T}$. The finite element assembly operator as outlined in Sec.~\ref{sec:basis_continuity} is denoted by $\mathsf{A}$. It will be used to assemble local basis functions coming from multiple supporting simplices into global basis functions. The support of the global basis function associated with a $k$-subsimplex $f$ will be termed the {\ital star of $f$}. Now note that the subspace $\fcspcminus{1}{k}{T}{f}{}$ of the Whitney $k$-forms $\PLminus{1}{k}{T}$ associated with a $k$-subsimplex $\ssimplx{f}{}$ on the coarse level is contained in the direct sum \begin{align}
	\fcspcminus{1}{k}{T}{f}{} \subset \bigoplus\limits_{g \in \mc{S}_k\incfreud{\simplx{T},f}} \assembly\limits_{H \in  A(\simplx{T},\ssimplx{f}{},\ssimplx{g}{})} \fcspcminus{1}{k}{H}{g}{}\,
\end{align}
of appropriately assembled subspaces associated with the finer discretisation, specified by the incomplete Freudenthal $\incfreud{\simplx{T},\ssimplx{f}{}}$.
The support of a fine-level basis function $\phi_{}$ associated with a subsimplex $\ssimplx{g}{}\in \mc{S}_k\incfreud{\simplx{T},\ssimplx{f}{}}$ is given by the set of subsimplices
\begin{align}
A(\simplx{T},\ssimplx{f}{},\ssimplx{g}{}) =  \bigcup\limits_{H \in \mc{S}_n\incfreud{\simplx{T},f}} \upadjacstar{\simplx{H},g} \,.
\end{align}
This allows for the introduction of the refinement equation
\begin{align}
c_{\bs{\sigma}}\phi_{\bs{\sigma}} = \sum\limits_{\ssimplx{g}{\mu} \in \mc{S}_k\incfreud{\simplx{T},\ssimplx{f}{\sigma}}} c_{\bs\sigma\bs\mu} \, \phi_{\bs{\mu}}\,,
\quad \phi_{\bs{\mu}} \in  \assembly\limits_{H \in  A(\simplx{T},\ssimplx{f}{\sigma},\ssimplx{g}{\mu})} \fcspcminus{1}{k}{H}{{\ssimplx{g}{\mu}}}{} \,, \quad c_{\bs{\sigma}}, c_{\bs{\sigma}\bs{\mu}} \in \mathbb R\ ,
\end{align}
which represents the replacement of a coarse basis function $\phi_{\bs{\sigma}}$ of $\ssimplx{f}{\sigma}$ with a linear combination of finer basis functions, namely all basis functions associated with the incomplete Freudenthal $\incfreud{\simplx{T},f}$. 

Of course, the unknown coefficients are underdetermined by one in the refinement equation, leaving one degree of freedom to constrain. In seeking geometric locality of the refinement constraint we exploit the locality of the Whitney forms, \ie taking the trace of the refinement equation onto $f$
\begin{align}
\sum\limits_{\ssimplx{g}{\mu} \in \hchildren{\ssimplx{f}{\sigma}}}c_{\bs{\sigma\mu}}\mrm{Tr}_{\ssimplx{f}{\sigma}}\phi_{\bs\mu} &= c_{\bs{\sigma}}\mrm{Tr}_{\ssimplx{f}{\sigma}}\phi_{\bs{\sigma}}  \,, \quad \phi_{\bs{\mu}} \in  \assembly\limits_{H \in  A(\simplx{T},\ssimplx{f}{\sigma},\ssimplx{g}{\mu})} \fcspcminus{1}{k}{H}{{\ssimplx{g}{\mu}}}{}  \,, \quad c_{\bs{\sigma}}, c_{\bs{\sigma}\bs{\mu}} \in \mathbb R\
\end{align}
reduces $\ssimplx{g}{\mu} \in \mc{S}_k\incfreud{\simplx{T},\ssimplx{f}{\sigma}}$ to $\ssimplx{g}{\mu} \in \hchildren{\ssimplx{f}{\sigma}}$ and thereby eliminates the influence of basis functions involved in refinements of $t \in \Dsimplx{T}\setminus f$. Developing the constraint on the basis of the trace of the refinement equation ensures that the refinement equation is not only fulfilled on $\simplx{T}$, but also, in the sense of traces, on $\ssimplx{f}{}$. Thereby linear independence of the global system of equations is ensured for arbitrary repeated refinements. While we have considered the refinement on a single simplex only for notational convenience, the procedure carries over to simplicial complexes by simply accounting for the fact that the support of the coarse level basis functions now may lie on multiple simplices and then assembling accordingly. Carrying this over to the finer level of the discretization renders the procedure's generalization to repeated refinements straightforward.

The single constraint to be determined should be linear but is otherwise arbitrary. Two conventional choices exist: the {\ital hierarchical basis} in which a single constraint on the coefficients $c_{\bs{\sigma}\bs{\mu}}$ is introduced, while retaining the unknown coefficient $c_{\bs{\sigma}}$, and the {\ital quasi-hierarchical basis} in which the coefficient $c_{\bs{\sigma}}$ is assumed to vanish, leaving all $c_{\bs{\sigma}\bs{\mu}}$ undetermined~\cite{Grinspun2003}. The idea for the hierarchical basis of $\mc P^-_{1}\Lambda^{0}$ is to retain the coarse basis function associated with the vertex $v$ while removing the basis function associated with the vertex $\hchildren{v}$, instead. More generally, for a hierarchical basis of $\mc P^-_{1}\Lambda^{k}$, there exists more than one subsimplex in $\hchildren{f}$ for $k>0$, such that the choice in eliminating a single degree of freedom is not unique. A sufficient constraint is given by the zero-mean property
 \begin{align}
\sum\limits_{t_{\bs\mu} \in \mc C(\ssimplx{f}{\sigma})}c_{\bs{\sigma\mu}} = 0
\label{eq:zeroMean}
\end{align}
 of the coefficients associated with the subsimplices $\hchildren{f}$. The quasi-hierarchical basis is usually more efficiently implementable than the hierarchical basis. In practice, when working with ascending local ordering of global vertex numbers and orientation-reversing transformations to a positively-oriented reference simplex, attention must be paid to the relative orientation of $\ssimplx{t}{\mu}$ and $\ssimplx{f}{\sigma}$ and the coefficient sign corresponding to $\ssimplx{t}{\mu}$ in Eq.~\ref{eq:zeroMean} must be reversed accordingly.

\subsection{Derefinement of the simplicial hierarchy}
The process of removing refinement depth in previously-refined regions is a process which we refer to as derefinement. This terminology shall distinguish the process from coarsening. To be effective in implementation, we require each atomic operation in the process to be non-destructive apart from the entity being de-refined. In short, this amounts to tagging leaf subsimplices which meet such non-destructibility criteria and removing them. After this step, due to the change in the leaf set, the non-destructibility criteria are now met by another set of leaf subsimplices, and the process is repeated until a whole region of the hierarchical complex is derefined. The non-destructibility criterion is simple. Consider a set of leaf subsimplices covering its hierarchical parent which we now refer to as a {\ital subleaf}. Such a set of leafs is considered to be safely removeable if the corresponding subleaf has no geometric descendants which are refined. In fact, due to the hierarchical refinement process and the incomplete Freudenthal, this criterion reduces to corresponding subleafs without geometric children which have been refined.

\subsection{The minimum rule for hp-hierarchical finite element spaces}

Taking $h$-refinement into account, let $\Omega^1_d = \Omega_d$ represent the $n$-simplices at the root level of the refinement, while all $n$-simplices up to level $L$ are contained in $\cup_{l=1}^L \Omega^l_d$. It then must be assured that the minimum rule is enforced on all leaf subsimplices $\ssimplx{f}{}$ of dimension $0 < k \le n$. To this end, an $n$-simplex's order is defined as the order of its hierarchical parent. Thus by recurrence, the polynomial propagates upwards either from the root of the tree or from an {\ital order-defining} hierarchical ancestor, which remains to be identified appropriately. With this definition, it is sensible in the $h$-hierarchical context to assume that the polynomial order for leaf-type subsimplices of dimension $0 < k \le n$ must uniquely be defined by their geometric ancestor $n$-simplices without reference to the hierarchical descendants thereof. 

This assumption directly leads to two possible limits of choice in the definition of order-defining $n$-simplices. The lower limit is the case, shown in Eq.~\ref{eq:minrule}, when the order-defining volumes $\Omega_d$ are simply chosen as the roots of the hierarchical trees. The upper limit is the case where an $n$-simplex $\simplx{T} \in \cup_{l=1}^L \Omega^l_d$ is considered to be order-defining when
\begin{enumerate}[(a)]
	\item at least one of the $f \in \Delta\simplx{T}$ is of leaf type and
	\item none of its hierarchical ancestors fulfill condition (a).
\end{enumerate}
The latter case leads to the definition of the hierarchical minimum rule for a leaf subsimplex $\ssimplx{f}{}$ of dimension $0 < k \le n$
 \begin{align}
	 f \mapsto \min \left\{\mathfrak{O}^*(\simplx{T}) \, : \, \ssimplx{f}{} \subset \simplx{T}, \; \simplx{T} \in \Omega_d^* \right\} \, ,
\end{align}
where $\Omega_d^*$ is the set of all order-defining $n$-simplices with roots in $\Omega_d$ and $\mathfrak{O}^* : \Omega^*_d \rightarrow \mathbb{N}$. Again, subsimplices of dimension zero are assigned the order one. For multiple unknown solution fields, care is to be taken that the potentially required polynomial order offsets between the fields are preserved. Trivially, this can be ensured by simply adding or subtracting such offsets from the order obtained by the rules outlined above.

\section{Numerical studies}

For the following problems, we adopt a general refinement strategy built upon the following hypotheses.
\begin{enumerate}[(a)]
	\item In high-order methods, singularities in the solution not aligned with the mesh cause polynomial oscillations (Runge's phenomenon) which may pollute the solution even at great distance to the singularity.
	\item Runge's phenomenon may only be suppressed using lowest-order approximation where the solution is singular.
	\item The accuracy and trustworthiness of a posteriori error indicators depend on the solution and mesh quality. Only a trustworthy error indicator can be a guide to refinement.
\end{enumerate}
From these hypotheses, we deduce a refinement algorithm, which we sketch in the following.

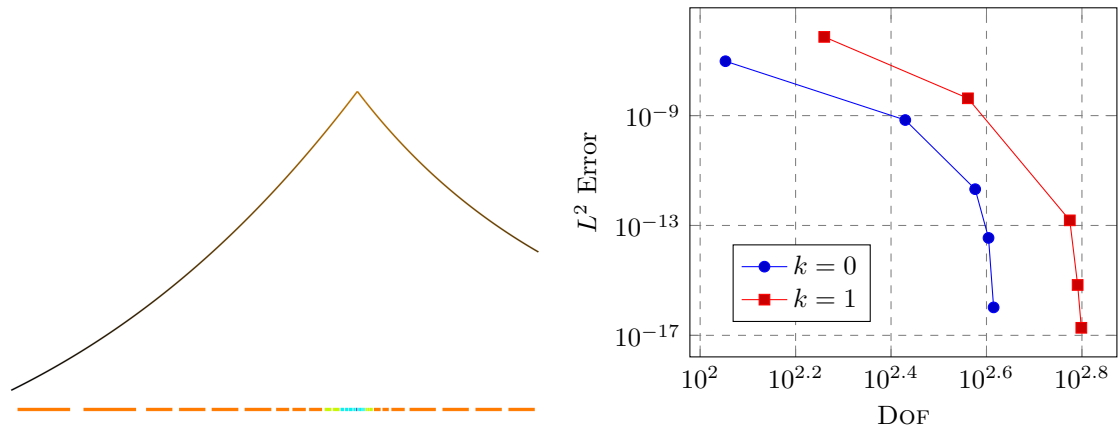
\begin{figure}[bht]
	\centering
	\includegraphics[width=.45\linewidth]{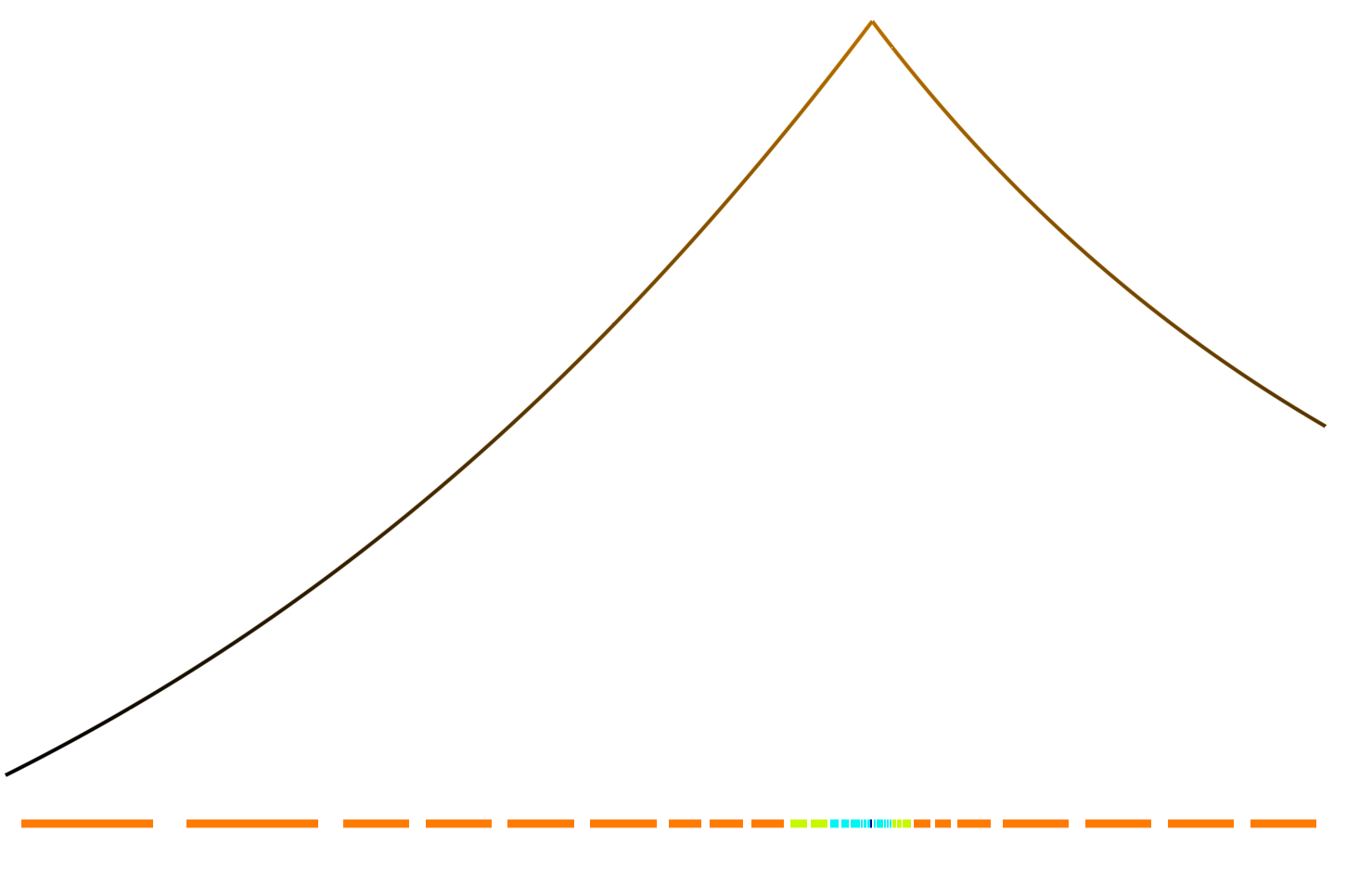}
\begin{tikzpicture}
    \begin{loglogaxis}[
        xlabel=\textsc{Dof},
        ylabel=$L^2$ Error,
        width=0.45\textwidth,
        grid style={dashed,gray},
        grid=both,
        minor tick num=2,
        legend style={at={(0.1,0.1)},anchor=south west}
    ]
 
           \addplot plot coordinates{
  (113,9.548154101647284e-8 )
 (269, 6.895516721315477e-10 )
 (377, 2.095832199415463e-12 )
 (402, 3.525016036646828e-14 )
 (412, 1.042044142571086e-16 )
    };

      \addplot plot coordinates {
 (182, 7.33627e-7) 
 (364, 4.26296e-9 )
 (595, 1.53706e-13)
 (618, 6.82813e-16)
 (629,1.88259e-17)
    };

      \legend{$k=0$,$k=1$}

    \end{loglogaxis}
\end{tikzpicture}
	
	\caption{Solution of the $L^2$ projection (k=1) and convergence (error computed using Monte Carlo integration)}
	\label{fig:L21D}
\end{figure}

First, we fix a maximum number of subdivision depths and initialize all polynomial orders to $r = 3$. This is done as the proposed a posteriori error indicator depends on a higher-order polynomial approximation, and $r=3$ was found to be sufficient for good performance. Then we solve the problem and subdivide the complex where the spectral decay indicates a loss of solution regularity. In regions of high regularity, we increase polynomial approximation order. We set to lowest-order approximation the cells with critically low solution regularity and solve the problem again, now with a more trustworthy a posteriori indicator. We then apply a de-refinement step to the the complex in regions where the solution shows high regularity. This process aims to ensure that at each stage, the solution is free of polynomial oscillations and that a valid a posteriori error-indicator can be computed. We repeat these steps until all singularities are localized to the maximum number of subdivisions and the solution is deemed regular everywhere in the complex with the exception of where lowest-order approximation is used. In a final sweep, we consider the areas of lowest approximation order to be finalized and iteratively refine and de-refine in $p$ and $h$ elsewhere as deemed appropriate.

We apply the algorithm to the solution of two simple $L^2$ projections onto the finite element basis for $k=0$ (piecewise continuous) and $k=1$ (piecewise discontinuous). In Fig.~\ref{fig:L21D} on the left, we depict the solution for $k=2$ together with the polynomial order distribution (orange: high, blue: low). On the right, we give the $L^2$ error computed using Monte Carlo integration for different numbers of degrees of freedom. We clearly observe exponential levels of convergence.

\subsection{Mixed finite element solution of the Hodge-Laplacian}
The Hodge Laplacian is the problem of finding $\bs{u} \in \Lambda^{k}(\Omega)$ such that
\eq{
(\extd\bs{\delta} + \bs{\delta}\extd)\bs{u} = \bs{f}-\bs{p}\,,  \quad \bs{p} \in \mathfrak{H}^k \,,
}
on a domain $\Omega \subset \field{R}^n$ where $\mathfrak{H}^k(\Omega) = \mrm{ker}\extd_k \setminus \mrm{im}\extd_{k-1}$ are the harmonic forms. Setting $k>0$, introducing $\bs\sigma \in \Lambda^{k-1}(\Omega)$ as an additional unknown, and enforcing orthogonality to the harmonic forms, the Hodge Laplacian problem reads
\eq{
	\bs{\sigma} = \bs{\delta} \bs{u} \,, \quad\quad
	\extd  \bs{\sigma} + \bs{\delta} \extd   \bs{u} +  \bs{p} = \bs{f} \,,\quad\quad
	P_{\mf{H}_k} \bs{u} = 0 \,.
}
The orthogonality to the harmonic forms is required as a consequence of the Hodge decomposition 
\eq{
	\Lambda^k \cong \mrm{im} \extd_{k-1} \oplus \mrm{im} \bs\delta_{k+1} \oplus \mathfrak{H}^k
}
associated with the de Rham complex
\eq{
 0 \rightarrow \Lambda^0(\Omega) \xrightarrow{\extd_0} \Lambda^1(\Omega) \xrightarrow{\extd_1} \ldots \xrightarrow{\extd_{n-1}} \Lambda^n(\Omega) \xrightarrow{\extd_n} 0 \, .
}%
Invoking the $L^2$ de Rham complex~\cite{Arnold2006}, the corresponding weak form, prior to integration by parts, is then the problem of finding $\bs{\sigma} \in H\Lambda^{k-1}(\Omega)$, $\bs{u} \in H\Lambda^{k}(\Omega)$, and $\bs{p} \in \mf{H}^k(\Omega)$ such that
\eq{
	\int\limits_\Omega\bs{\tau} \wedge \star \bs{\sigma}
	=  \int\limits_\Omega\bs{\tau} \wedge \star \bs{\delta} \bs{u}
	\quad\quad &\forall  \bs{\tau} \in H\Lambda^{k-1}\,, \\
	\int\limits_\Omega\bs{v} \wedge \star \extd  \bs{\sigma} + \int\limits_\Omega\bs{v} \wedge \star \bs{\delta} \extd   \bs{u} +  \int\limits_\Omega\bs{v} \wedge \star \bs{p}
	= \int\limits_\Omega\bs{v} \wedge \star \bs{f}  \quad\quad
	&\forall  \bs{v} \in H\Lambda^{k} \,,\\
	\int\limits_\Omega\bs{q} \wedge \star  \bs{u}
	=  0 \quad\quad
	&\forall  \bs{q} \in \mf{H}^k \, .
}
The boundary condition for $\bs\sigma$ can strongly be imposed by seeking its solution in a space of prescribed boundary trace.

\begin{figure}[htb]
\centering
\includegraphics[width=0.45\textwidth]{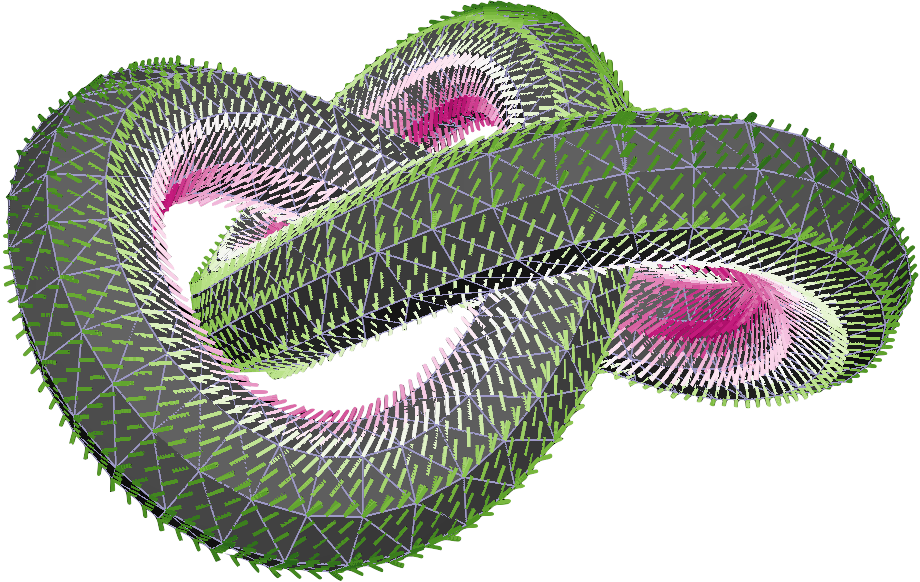}\hspace{0.7cm}
\includegraphics[width=0.45\textwidth]{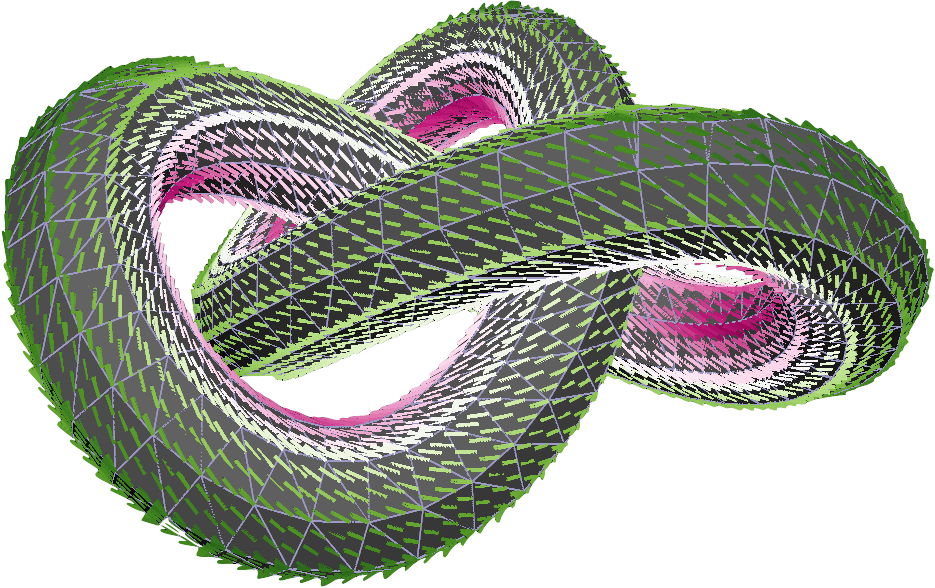}
\caption{Basis for the harmonic forms $\mathfrak{H}^1$ on a torus with Betti number $b_1 = 2$ embedded in $\mathbb{R}^3$}
\label{fig:harmonicforms}
\end{figure}

The terms from which the boundary conditions for $\star\bs{u}$ and $\star\extd\bs{u}$ can be found are $\bs{\tau} \wedge \star \bs{\delta} \bs{u}$ and $\bs{v} \wedge {\star \bs{\delta}} \extd  \bs{u}$, respectively. Using the Leibniz rule from Eq.~\ref{eq:leibnizrule}, these terms may be restated as follows:
\eq{
\bs{\tau} \wedge \star \bs{\delta} \bs{u}
&= \bs{\tau} \wedge   {\extd  \left[(-1)^k{\star \bs{u}} \right]}
=
{\extd } \bs{\tau} \wedge {\star \bs{u}}
- {\extd } \left( \bs{\tau} \wedge {\star \bs{u}}  \right)
 \\
\bs{v} \wedge {\star \bs{\delta}} \extd  \bs{u}
&= \bs{v} \wedge  {\extd  \left[ {(-1)^{k+1}\star \extd  \bs{u}}   \right] }
=  {\extd } \bs{v} \wedge  {\star \extd  \bs{u}}
- {\extd } \left( \bs{v} \wedge {\star \extd  \bs{u}}    \right) \,.
}
Application of Stokes' theorem then yields the natural boundary conditions such that the weak formulation becomes
\eq{
		\int\limits_{\Omega} \bs{\tau} \wedge \star \bs{\sigma}
		- \int\limits_{\Omega}{\extd } \bs{\tau} \wedge {\star \bs{u}}
=  -\int\limits_{\partial\Omega} \mrm{Tr} \left( \bs{\tau} \wedge {\star \bs{u}}  \right)
	\quad\quad &\forall  \bs{\tau} \in \Lambda^{k-1}\,, \\
	\int\limits_{\Omega} \bs{v} \wedge \star \extd  \bs{\sigma}
	 + \int\limits_{\Omega}  {\extd } \bs{v} \wedge  {\star \extd  \bs{u}}
	 + \int\limits_{\Omega} \bs{v} \wedge \star \bs{p}
	= \int\limits_{\Omega}\bs{v} \wedge \star \bs{f}
	+ \int\limits_{\partial\Omega} \mrm{Tr} \left( \bs{v} \wedge {\star \extd  \bs{u}}    \right)  \quad\quad
	&\forall  \bs{v} \in \Lambda^{k}\,, \\
	\int\limits_{\Omega} \bs{q} \wedge \star  \bs{u}
	=  0 \quad\quad
	&\forall  \bs{q} \in \mf{H}^k \,.
}

For verification of our implementation without adaptivity, we compute the kernel of the Hodge Laplacian problem with $k = n-1 = 1$ on the torus with Betti number $b_1 = 2$ embedded in $\mathbb{R}^3$. Using lowest-order elements of type $\bs\sigma \in \mc{P}_1^-\Lambda^{0}$ and $\bs u \in \mc{P}_1^-\Lambda^{1}$ we thus obtain an approximation of the basis for the harmonic forms $\mathfrak{H}^1$, which is shown in Fig.~\ref{fig:harmonicforms}. It can clearly be seen that the covector fields obtained circulate around the 1-holes of the torus.

\begin{figure}[htb]
\begin{minipage}{0.49\textwidth}

\begin{tikzpicture}[every node/.style={inner sep=2,outer sep=0}]
\node (label) at (0,0)[]{
	\includegraphics[width=\textwidth]{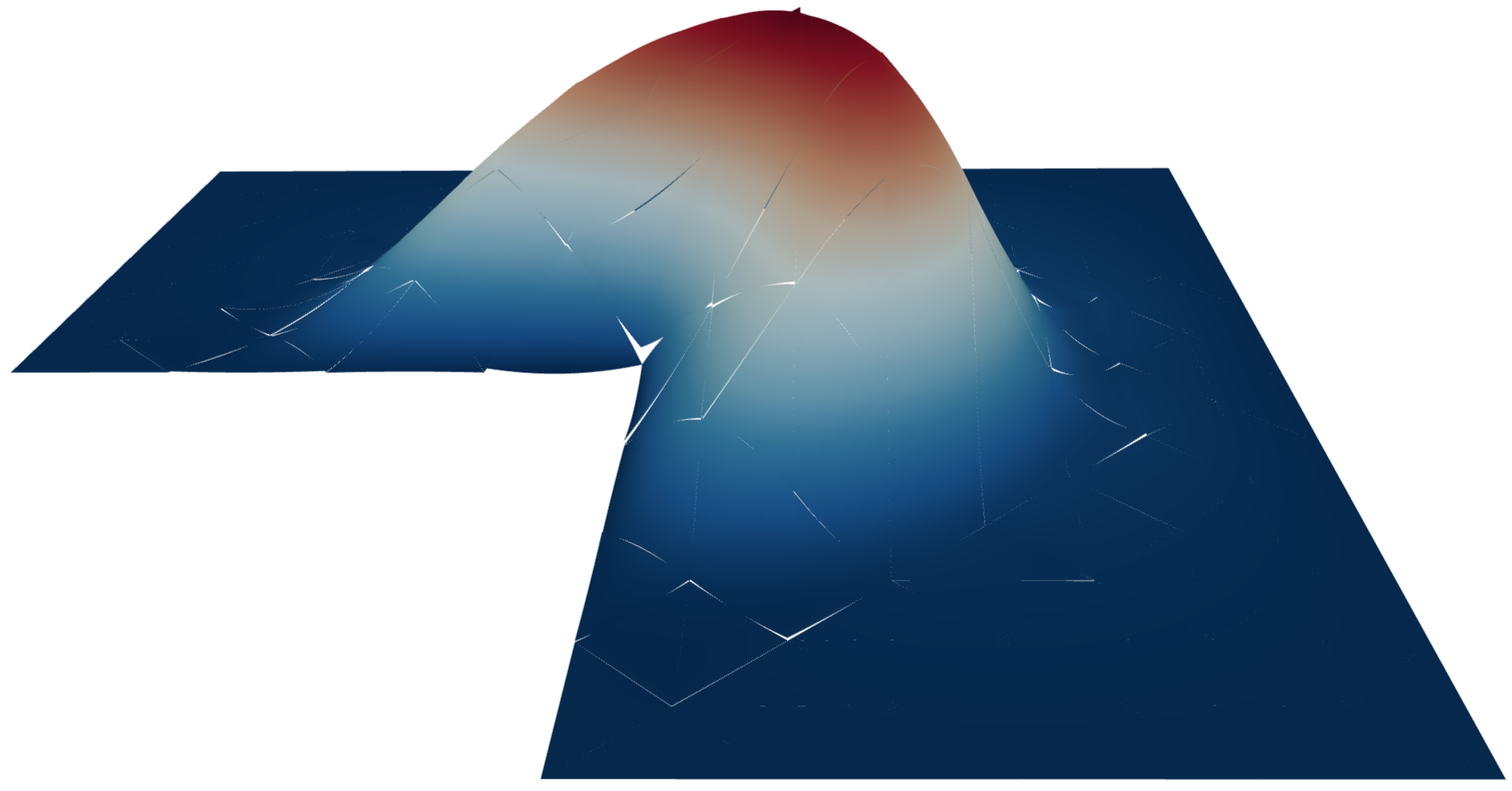}
      };
      \node (label) at (0,-2.5)[]{
\pgfplotscolorbardrawstandalone[ 
    colormap name = BuRd,
    colorbar horizontal,
    point meta min=-4.058e-3,
    point meta max=0.5049,
    colorbar style={
        width=\textwidth-1.5cm,
        height=.25cm,
        }
        ]};
\end{tikzpicture}
\end{minipage}
\hfill
\begin{minipage}{0.49\textwidth}
	\begin{tikzpicture}[every node/.style={inner sep=2,outer sep=0}]
\node (label) at (0,0)[]{
	\includegraphics[width=\textwidth]{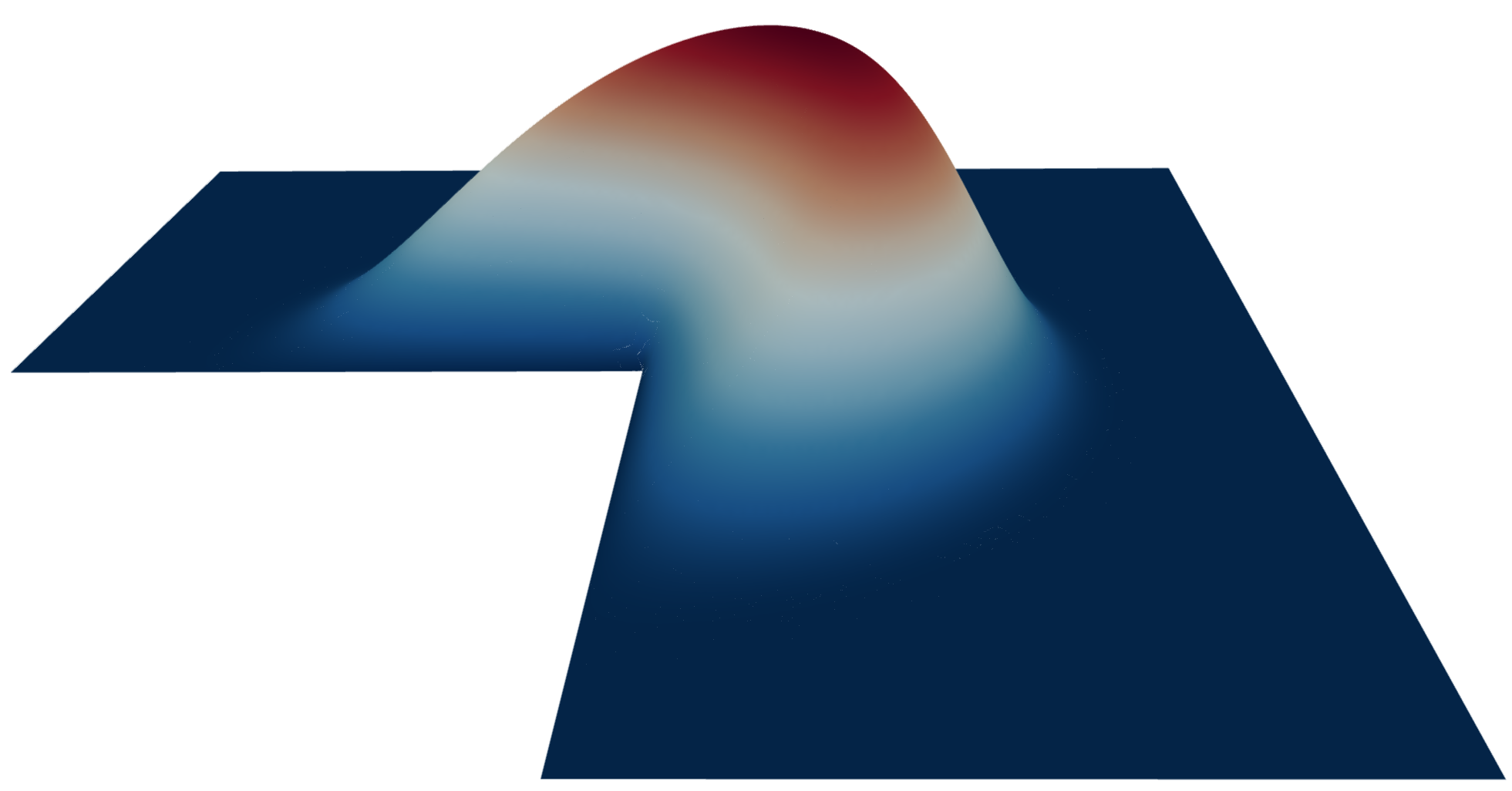}
      };
      \node (label) at (0,-2.5)[]{
\pgfplotscolorbardrawstandalone[ 
    colormap name = BuRd,
    colorbar horizontal,
    point meta min=-4.056e-5,
    point meta max=0.4682,
    colorbar style={
        width=\textwidth-1.5cm,
        height=.25cm,
        }
        ]};
\end{tikzpicture}
\end{minipage} \caption{Solutions $\bs{u} \in \Lambda^2$ for the initial (left) and the refined problem (right)}
\label{fig:solveclap}
\end{figure}

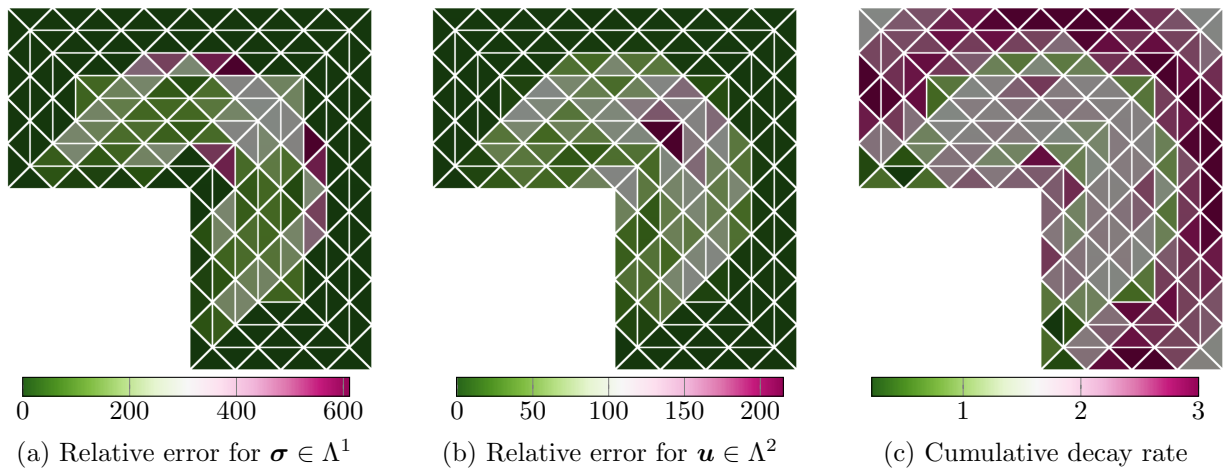
\begin{figure}[htb]
\begin{minipage}{0.3\textwidth}
\centering
\begin{tikzpicture}[every node/.style={inner sep=2,outer sep=0}]
\node (label) at (0,0)[]{
	\includegraphics[width=\textwidth]{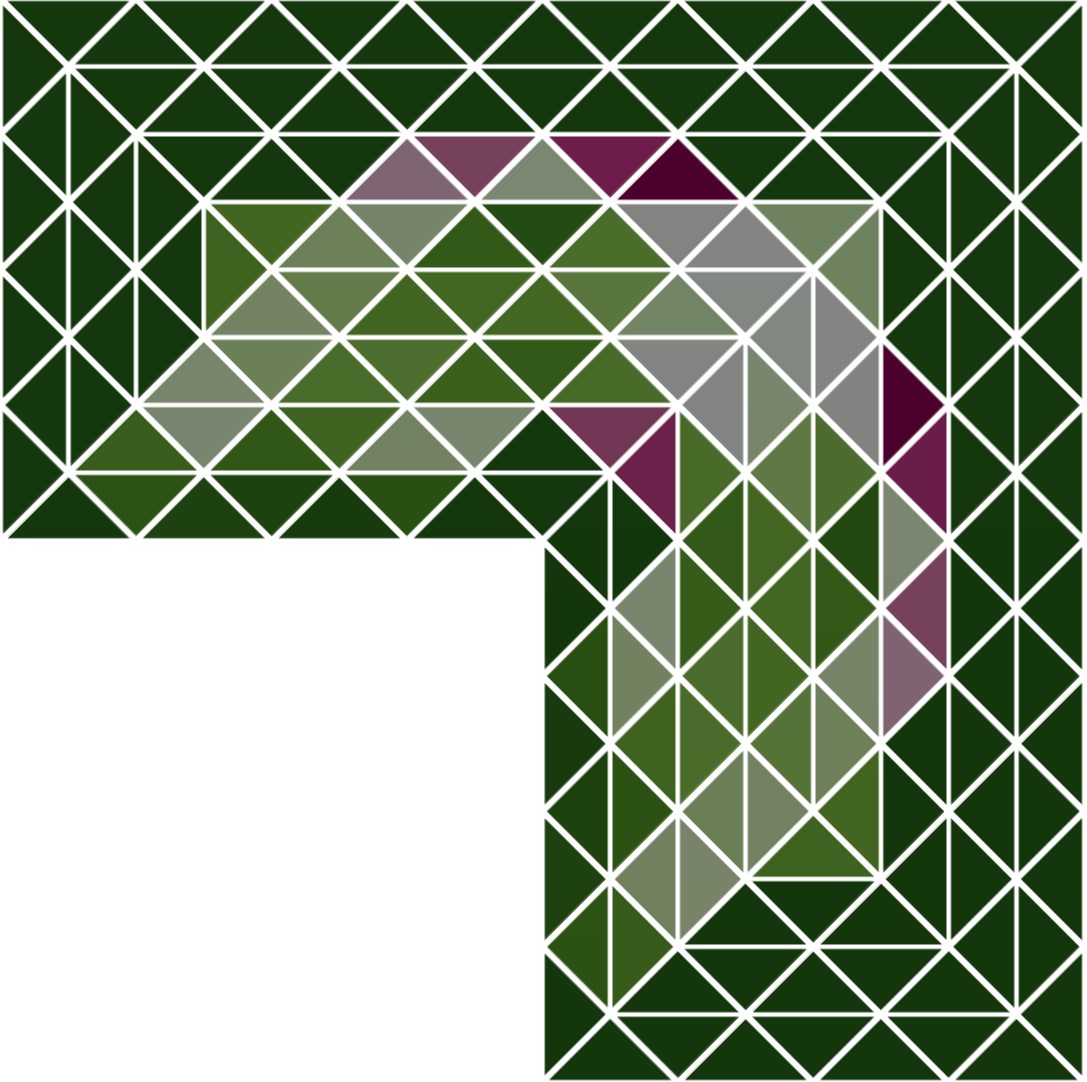}
      };
      \node (label) at (0,-2.8)[]{
\pgfplotscolorbardrawstandalone[ 
    colormap name = arnold,
    colorbar horizontal,
    point meta min=0.0,
    point meta max=611.7,
    colorbar style={
        width=\textwidth-.5cm,
        height=.25cm,
        }
        ]};
\end{tikzpicture}\\
(a) Relative error for $\bs{\sigma} \in {\Lambda^1}$
\end{minipage}
\hfill
\begin{minipage}{0.3\textwidth}
\centering
\begin{tikzpicture}[every node/.style={inner sep=2,outer sep=0}]
\node (label) at (0,0)[]{
	\includegraphics[width=\textwidth]{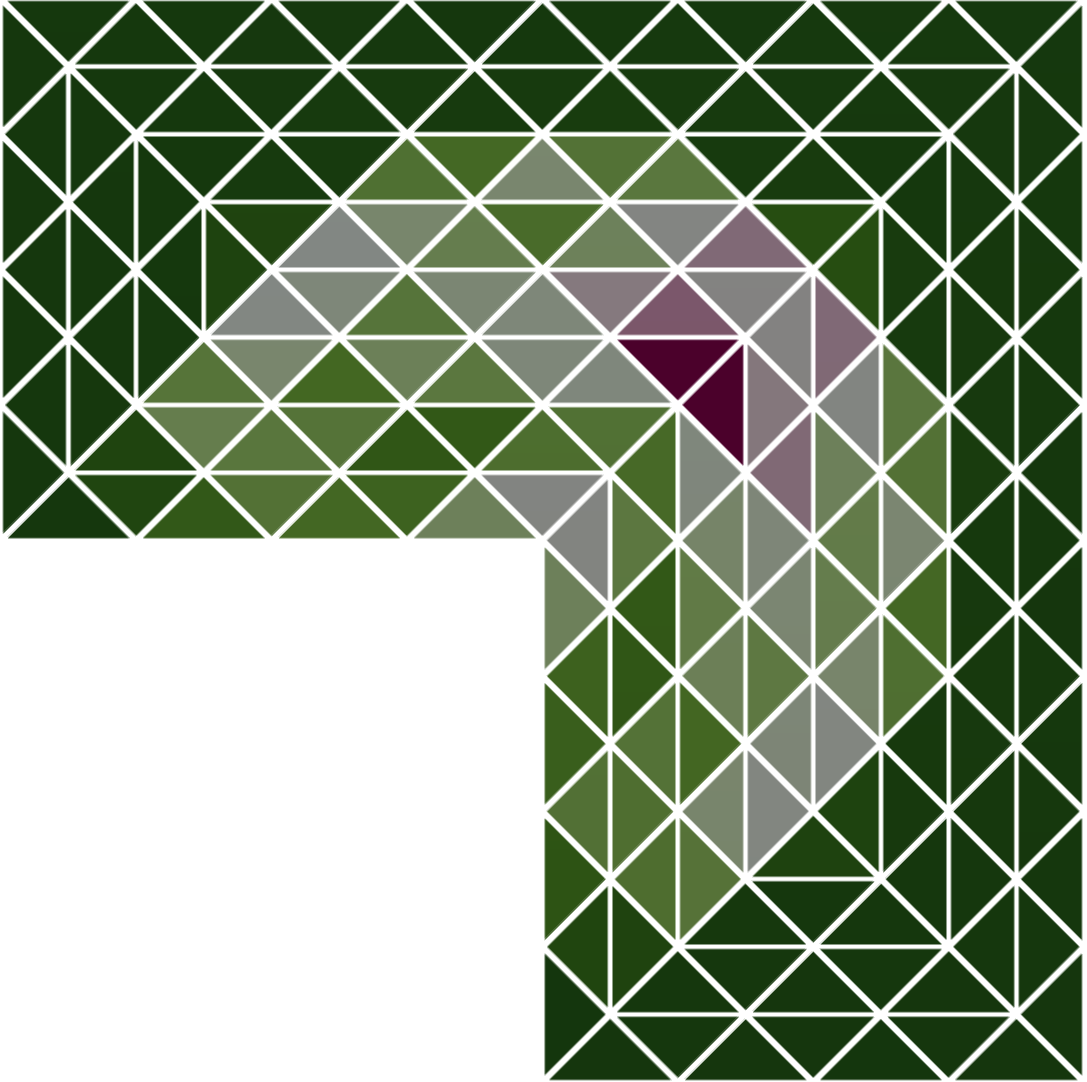}
      };
      \node (label) at (0,-2.8)[]{
\pgfplotscolorbardrawstandalone[ 
    colormap name = arnold,
    colorbar horizontal,
    point meta min=-3.060e-1,
    point meta max=2.152e2,
    colorbar style={
        width=\textwidth-.5cm,
        height=.25cm,
        }
        ]};
\end{tikzpicture}\\
(b) Relative error for $\bs{u} \in {\Lambda^2}$
\end{minipage}
\hfill
\begin{minipage}{0.3\textwidth}
\centering
\begin{tikzpicture}[every node/.style={inner sep=2,outer sep=0}]
\node (label) at (0,0)[]{
	\includegraphics[width=\textwidth]{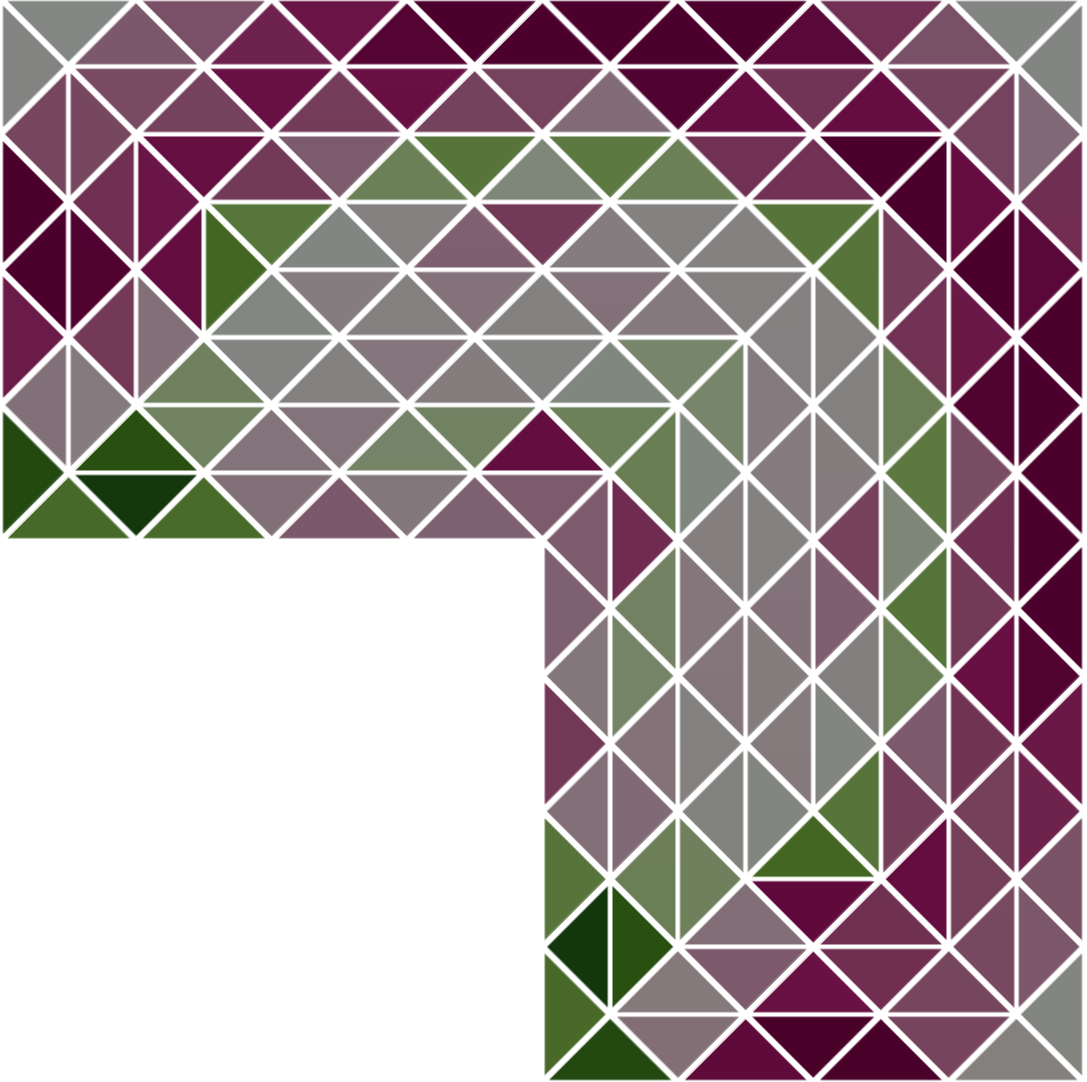}
      };
      \node (label) at (0,-2.8)[]{
\pgfplotscolorbardrawstandalone[ 
    colormap name = arnold,
    colorbar horizontal,
    point meta min=2.24e-1,
    point meta max=3,
    colorbar style={
        width=\textwidth-.5cm,
        height=.25cm,
        }
        ]};
\end{tikzpicture}\\
(c) Cumulative decay rate
\end{minipage}
 \caption{Errors and decay rate for the initial Laplacian problem}
\label{fig:errdecayveclapinit}
\end{figure}

\begin{figure}[htb!]
\begin{minipage}{0.45\textwidth}
\centering
\begin{tikzpicture}[every node/.style={inner sep=2,outer sep=0}]
\node (label) at (0,0)[]{
	\includegraphics[width=\textwidth]{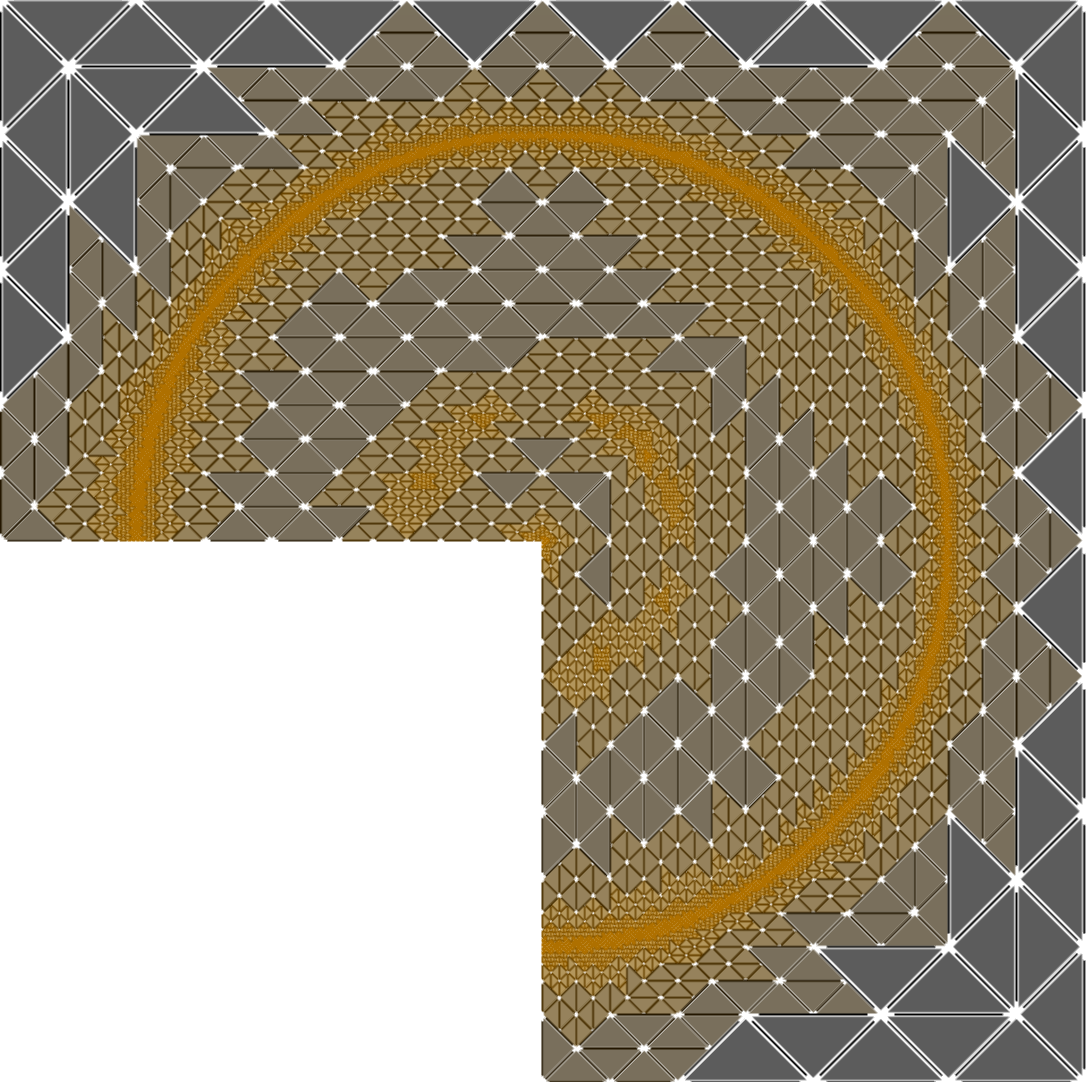}
      };
      \node (label) at (0,-4.0)[]{
\pgfplotscolorbardrawstandalone[ 
    colormap name = BkOr,
    colorbar horizontal,
    point meta min=1,
    point meta max=6,
    colorbar style={
        width=\textwidth-.5cm,
        height=.25cm,
        }
        ]};
\end{tikzpicture}\\
(a) Hierarchical subdivision depth
\end{minipage}
\hfill
\begin{minipage}{0.45\textwidth}
\centering
\begin{tikzpicture}[every node/.style={inner sep=2,outer sep=0}]
\node (label) at (0,0)[]{
\includegraphics[width=\textwidth]{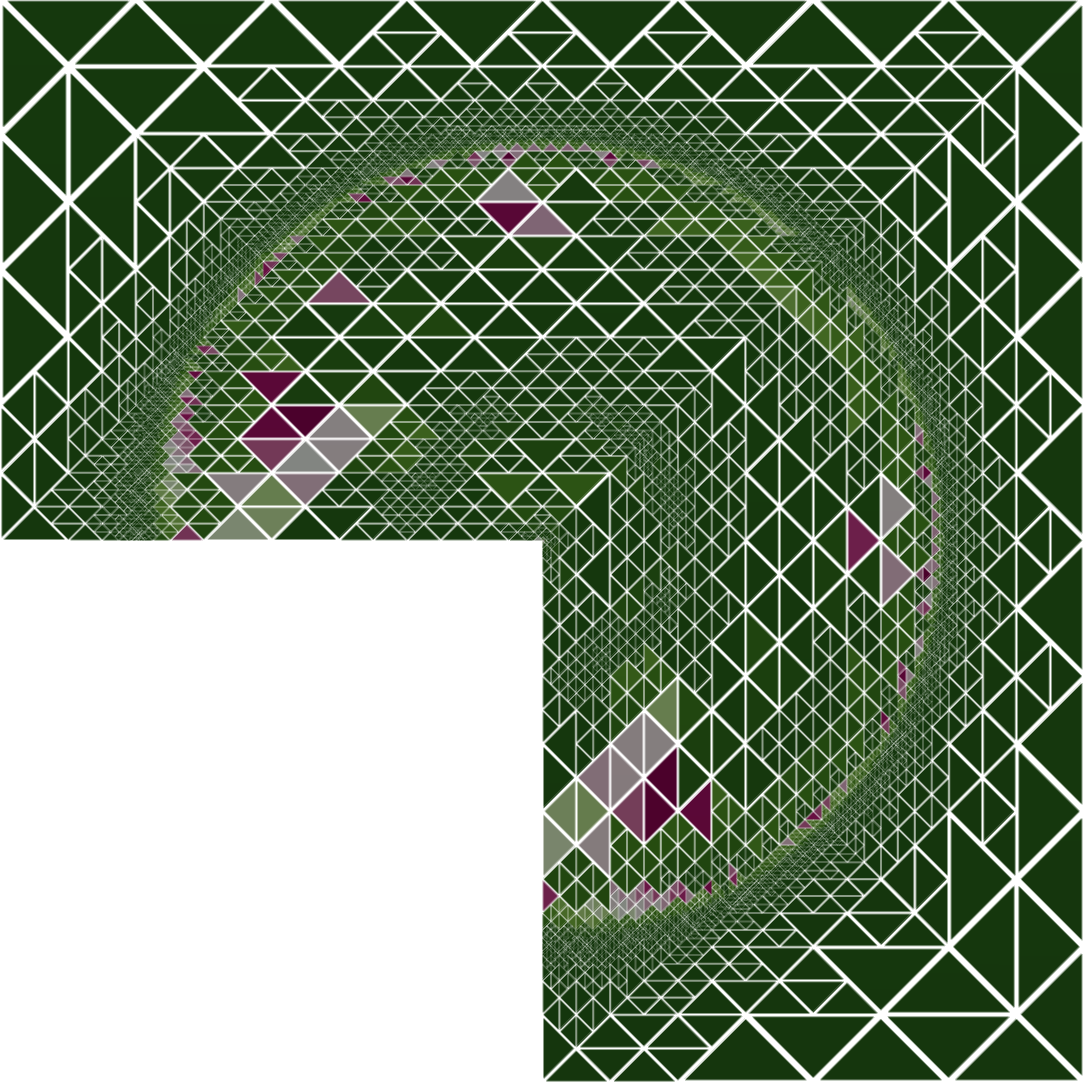}
      };
      \node (label) at (0,-4.0)[]{
\pgfplotscolorbardrawstandalone[ 
    colormap name = arnold,
    colorbar horizontal,
    point meta min=0.0,
    point meta max=1.659,
    colorbar style={
        width=\textwidth-.5cm,
        height=.25cm,
        }
        ]};
\end{tikzpicture}\\
(b) Relative error for $\bs{u} \in \mc{P}^-\Lambda^2$
\end{minipage}\vspace{.3cm}\\
\begin{minipage}{0.45\textwidth}
\centering
\begin{tikzpicture}[every node/.style={inner sep=2,outer sep=0}]
\node (label) at (0,0)[]{
\includegraphics[width=\textwidth]{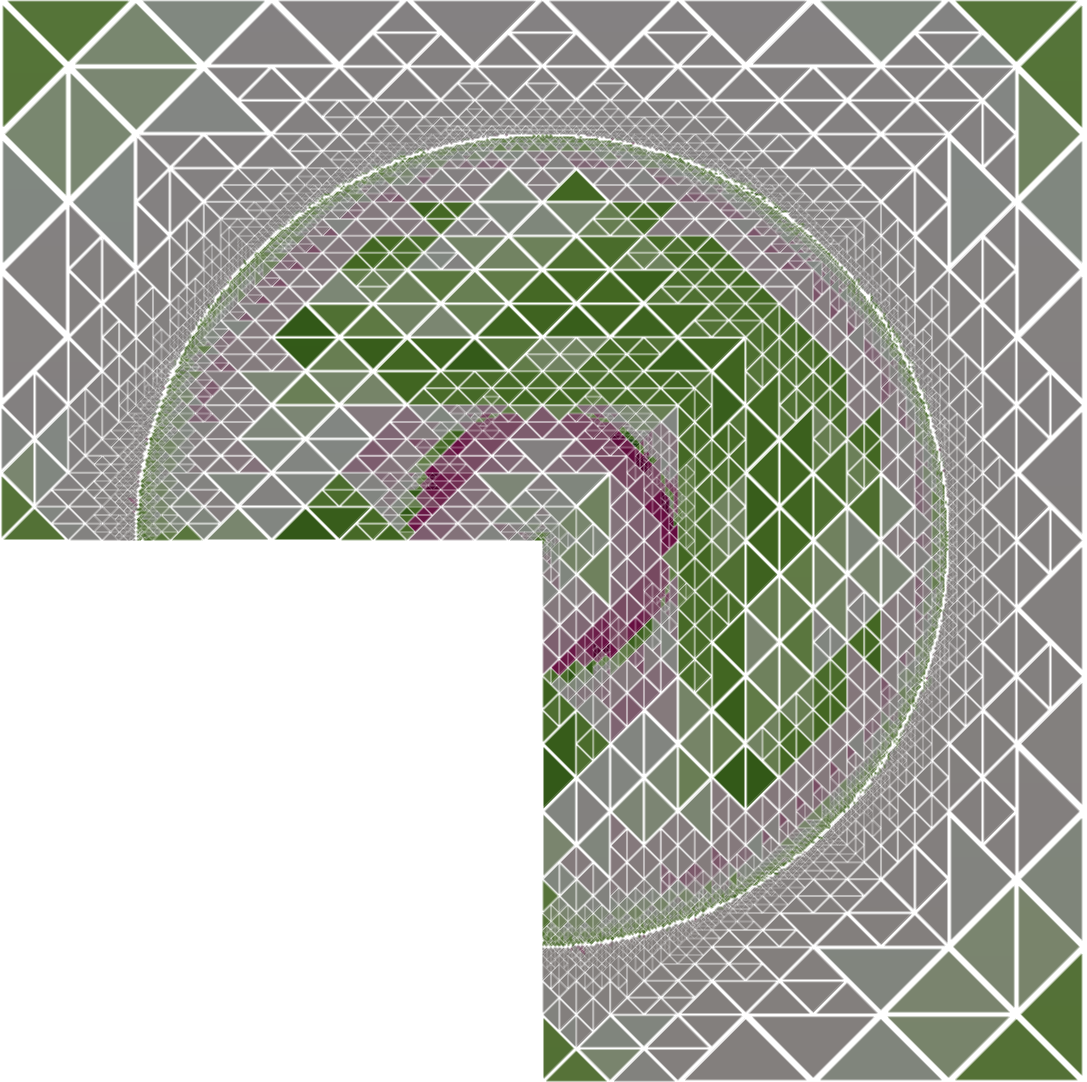}
      };
      \node (label) at (0,-4.0)[]{
\pgfplotscolorbardrawstandalone[ 
    colormap name = arnold,
    colorbar horizontal,
    point meta min=0.5278,
    point meta max=5,
    colorbar style={
        width=\textwidth-.5cm,
        height=.25cm,
        }
        ]};
\end{tikzpicture}\\
(c) Cumulative decay rate
\end{minipage}
\hfill
\begin{minipage}{0.45\textwidth}
\centering
\begin{tikzpicture}[every node/.style={inner sep=2,outer sep=0}]
\node (label) at (0,0)[]{
\includegraphics[width=\textwidth]{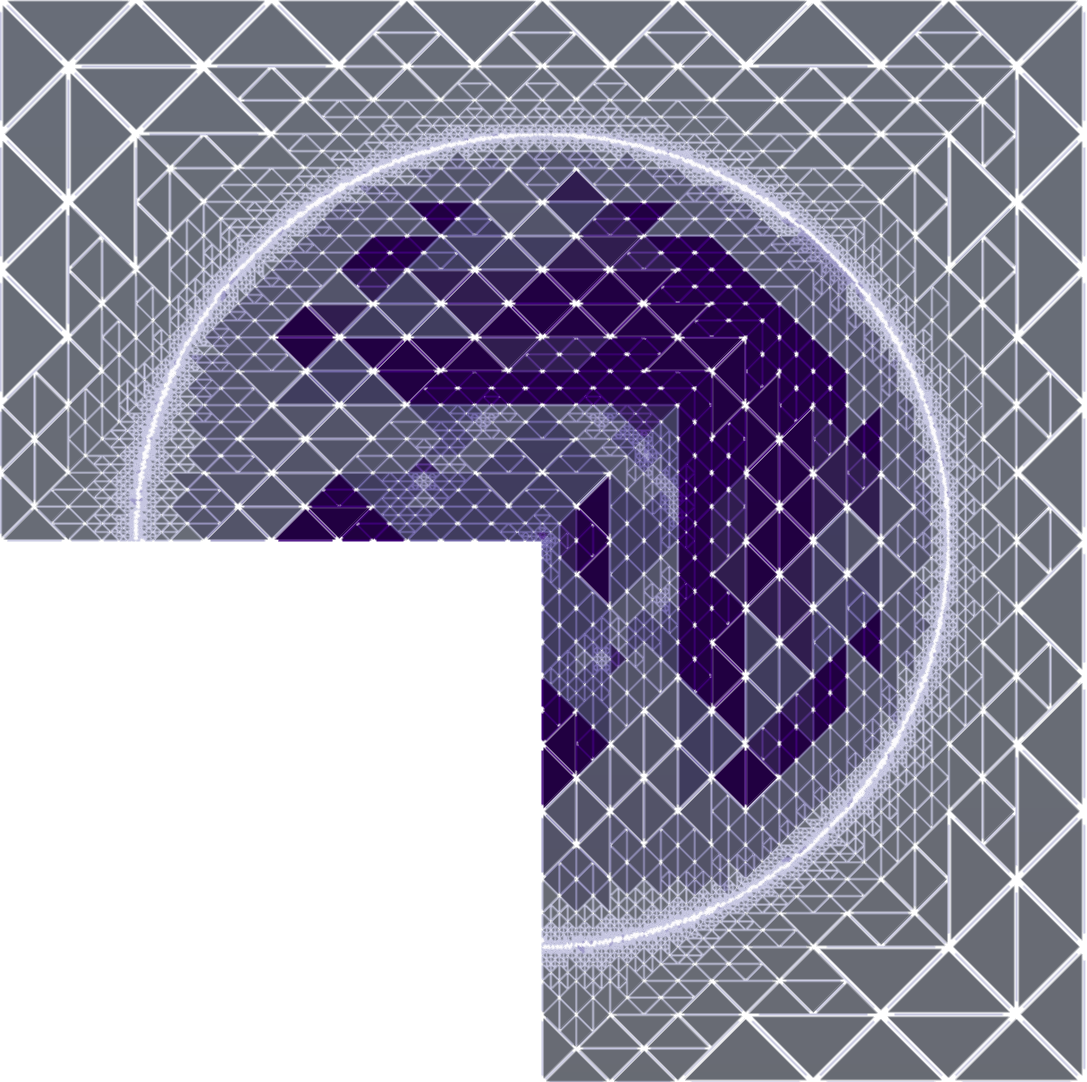}
      };
      \node (label) at (0,-4.0)[  ]{
\pgfplotscolorbardrawstandalone[ 
    colormap name = Purples,
    colorbar horizontal,
    point meta min=1,
    point meta max=7,
    colorbar style={
        width=\textwidth-.5cm,
        height=.25cm,
        }
        ]};
\end{tikzpicture}\\
(d) Polynomial order distribution
\end{minipage}
 \caption{Error, decay and polynomial order distribution for the refined Laplacian problem}
\label{fig:errdecayveclapref}
\end{figure}
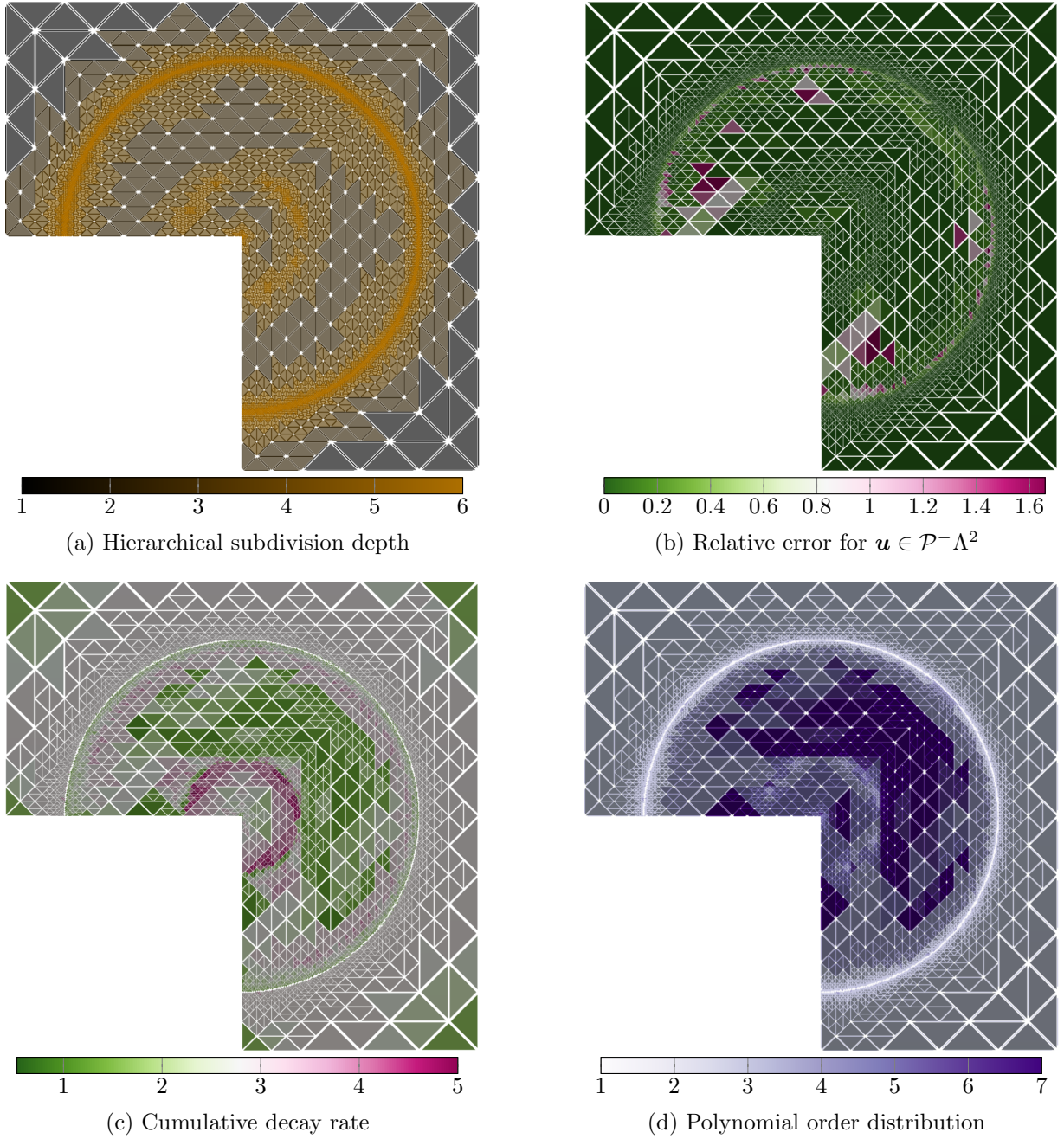

To demonstrate the potential of adaptivity, we solve the Hodge Laplacian problem for $k = n = 2$, \ie the mixed formulation of the scalar Laplacian, on an L-shaped domain. By Arnold, Falk and Winther (2006)~\cite{Arnold2006}, we choose as a stable discretization to the problem the finite element space $\bs\sigma \in \mc{P}_r^-\Lambda^{1}$, $\bs u \in \mc{P}_r^-\Lambda^{2}$. As a right-hand-side we choose a function singular on two concentric rings. The results are depicted in Figures~\ref{fig:solveclap}, \ref{fig:errdecayveclapinit} and \ref{fig:errdecayveclapref}. Figure~\ref{fig:solveclap} shows a comparison of the solutions for the initial and the final refined problem.\footnote{While the discontinuities in the left rendering are to be interpreted as solution errors, we note that at sufficient magnification, discontinuities also seem to appear in the rendering of the refined solution on the right. This is simply a consequence of our post-processing routine which only coarsely samples the polynomials on simplices at high refinement depths.} Figure~\ref{fig:errdecayveclapref} depicts the error, decay and polynomial order distribution for the refined problem, clearly indicating the concentric singularities propagated to the solution and the low or lowest-order approximation within the singularities.

\begin{figure}[htb]
\centering
\includegraphics[width=0.45\textwidth]{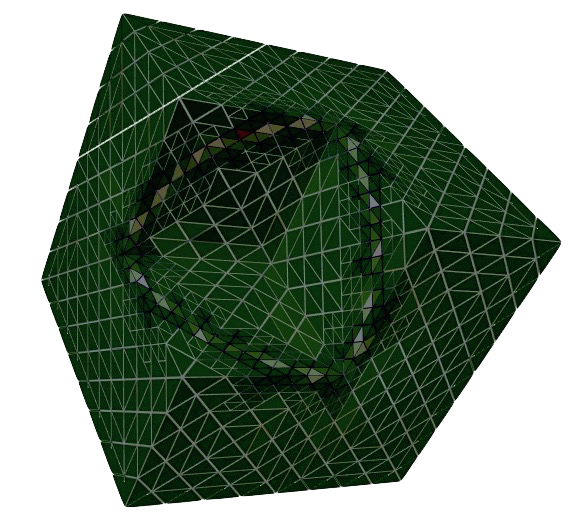}\hspace{0.7cm}
\includegraphics[width=0.45\textwidth]{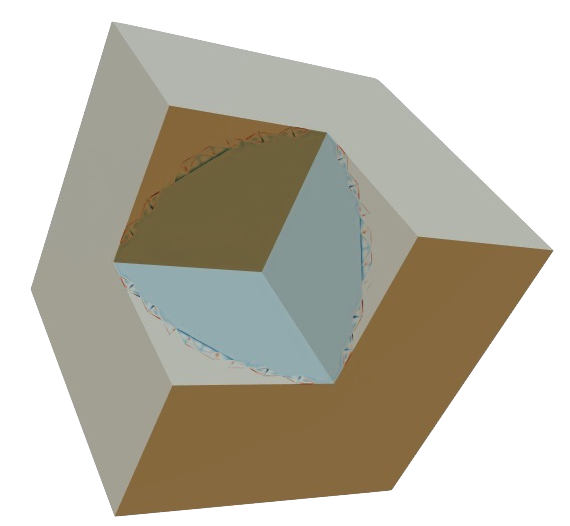}
\caption{Fichera's corner problem of the Hodge Laplacian for $k=n=3$ with a single level of refinement depth. Error indicator on the left, solution for $\bs{u}$ on the right.}
\label{fig:ficheracorner}
\end{figure}

As a proof of the dimension-independent nature of the method outlined herein, we also solve a small example in three dimensions for $k = n = 3$ of Fichera's corner problem. The finite element space employed is $\bs\sigma \in \mc{P}_r^-\Lambda^{2}$, $\bs u \in \mc{P}_r^-\Lambda^{3}$. We use a manufactured solution $\bs{u}$ corresponding to the indicator function of a ball with radius one and impose the right-hand-side accordingly. Then the mesh is subdivided only once without setting to lowest order the elements lying on the singularity. The results are depicted in Fig.~\ref{fig:ficheracorner}, where it can clearly be seen that the polynomial oscillations in the solution $\bs{u}$ are confined to the elements impacted by the singularity without propagating throughout the mesh. This result is due to the discontinuous approximation $\bs{u}\in \mc{P}_r^-\Lambda^{3}$.

Although the examples given for the Hodge Laplacian only employ the pairing $\bs\sigma \in \mc{P}_r^-\Lambda^{k-1}$ and $\bs u \in \mc{P}_r^-\Lambda^{k}$, we verified that the method works for all stable pairings of the Hodge Laplacian as given in~\cite{Arnold2006}. A more elaborate example which demonstrates more complex pairings is given in the following section in the context of linearized elasticity.

\subsection{Mixed finite element solution of incompressible linearized elasticity}

In the linearization of the general elasticity problem of a reference body $\Omega$ embedded in an ambient space $\mc{S}$, we are left with computing an infinitesimal displacement field $\star\bs{u} \in \varphi^*\tanbndl{S}$ which arises from the linearization of the motion $\varphi$ with respect to an admissible $\epsilon$-variation. The linearization of the deformation gradient then is the covariant derivative of $\star\bs{u}$ with respect to the induced connection ${{}^{\varphi}\nabla}$. The linearization of the right Cauchy Green tensor, in turn involves only the metric, the deformation gradient and its linearization. Linearization about a null-deformation, gives $\mrm{symm}{\nabla}(\star\bs{u})$ as the linearization of the right Cauchy Green tensor. We will now assume $\mc{S} = \field{R}^n$ and call all such tensor fields the space of linear strains $\mc{E} =  \Lambda^1(\Omega; T^*\Omega)$. The space of stresses considered as vector-valued $n-1$-forms is denoted as $\Sigma =  \Lambda^{n-1}(\Omega; \field{R}^n)$ and the constitutive relation is given by $\star_E : \Sigma \rightarrow \mc{E}$. Using the divergence in order to relate stress to force density, we arrive at the following diagram:
\eq{
	 \Lambda^0(\Omega; \field{R}^n )  \xrightarrow{\mrm{symm}{\nabla}} \mc{E} \xrightarrow{{\star^{-1}_E}} \Sigma \xrightarrow{{\mrm{div}}}   \Lambda^n(\Omega;\field{R}^n) \, .
}
If we restrict to three dimensions, we have $\sigma_{ij} = \epsilon_{ikm}\epsilon_{jln} \Phi_{kl,mn}$ where $\Phi_{kl} \in \mc{E}$ is the Beltrami stress tensor. Then the diagram given above is a vector-valued reduced version of the de Rham complex which is termed the {\ital complex of linear elasticity}~\cite{Arnold2006, Eastwood2000} and thus is an exact sequence. As shown in~\cite{Eastwood2000}, it inherits all properties of the de Rham complex and can be related to the Calabi complex. 

In the finite element exterior calculus~\cite{Arnold2006}, a stable and consistent discretization of the elasticity complex with weak symmetry is given. The following example of Cook's membrane will demonstrate its performance in the $hp$-hierarchical context. From a mechanical standpoint, it has been shown~\cite{Stein1990} that discrete stability in linearized plane elasticity is essentially the satisfaction of three conditions: the stress jump condition at inter-element boundaries, the absence of spurious modes on a two-element-patch, and the absence of zero-energy stresses on an element. We will begin with the nature of stress in the formulation, and then derive the constitutive relation as well as the weak formulation of the mixed problem. Finally, we discuss stress boundary conditions and present a numerical solution.

\subsubsection{The definition of linear momentum flux} In the context of differential forms we define the linear momentum flux $\mbb{t}$ as a vector-valued $n-1$-form, related to the Cauchy stress tensor by ${\star\mbb{t}} = \bs{\sigma}$. In the present case, it is the goal that each component of the linear momentum flux be expanded in a basis of finite element differential forms of degree $n-1$ defined on the reference element. It is therefore sensible to expand the stress vectors for each plane in the standard basis of the Euclidean surrounding space with basis vectors $\{\bs{e}^I\}_{I=1}^n$, \ie $\bs{t} = t_I \bs{e}^I$ while defining the area normals in the basis of the reference element with coordinates $\{\xi^i\}_{i=1}^n$. Its mapping to Euclidean ambient space is $\psi = \varphi\circ\bar{\eta}_i$ which induces a metric $\psi^*\bs{\mrm{Id}}$ on the reference element, in the following denoted by $\bs{g}$. The Cauchy stress tensor, considered in implementation point-wisely as an element of $\Sigma = \Lambda^{n-1}(T^*_P \simplx{T}, \field{R}^n)$, then is expanded as $\bs{\sigma} = \sigma_{Ij} \bs{e}^I \otimes \extd \xi^j$. While $\bs{\sigma}$ conceptually is symmetric, the components $\sigma_{Ij}$ in this choice of basis are, in general, not symmetric. This choice of basis for the Cauchy stress can be considered as natural from a physical standpoint, as stress vectors can be defined sensibly in a global setting whereas the planes upon which they act are more conveniently described relative to the material, \eg crystallographic planes which need not be aligned globally. An analogue choice of basis can be made for the linear momentum flux such that $\mbb{t} = \mbb{t}_{\bs\sigma I} \, \bs{e}^I \otimes \extd \xi^\bs{\sigma}$.

\subsubsection{Elasticity Hodge star}
Initially, consider the compliance tensor, \ie the metric associated with the strain energy inner product $\Sigma \times \Sigma \rightarrow \mathbb{R}^+$:
\eq{
	\mbb{D}(\bs{\tau}, \bs{\sigma}) = \frac{1}{2\mu} \left( \bs{\tau} : \bs{\sigma} - \frac{\lambda}{2\mu + n\lambda} \mrm{tr} \bs{\tau} \cdot \mrm{tr}\bs{\sigma} \right) \, .
}
In the incompressible limit, as $\lambda \rightarrow \infty$, the factor $c = \frac{\lambda}{2\mu + n\lambda}$ becomes $1/n$. In two dimensions, the latter definition corresponds to the plane strain case. For a compressible material under plane stress, replace $\lambda$ in $c$ with $\lambda^* = 2\lambda\mu / (\lambda + 2\mu)$, while for the incompressible case, $c$ becomes $1/(1+n)$. Assuming the previously introduced local representation of the Cauchy stress, the trace is to be interpreted as
\eq{
	\mrm{tr} \left(\sigma_{Ij} \bs{e}^I \otimes \extd x^j \right) = \underbrace{\left(\frac{\bs{\partial}}{\bs{\partial} \xi^i}\right)^I}_{\tensor{J}{^I_i}}\sigma_{Ik} g^{ki} \, .
}
Then, the local representation of $\mbb{D}$ is given by
\eq{
	\frac{1}{2\mu} \left(
	\delta^{KI}g^{jl}    - \frac{\lambda}{2\mu + n\lambda}
	 \tensor{J}{^I_\alpha}  g^{j\alpha}
	\tensor{J}{^K_\beta}  g^{l\beta} \right) \;\bs{e}_I \otimes \frac{\bs{\partial}}{\bs{\partial} \xi^j} \otimes \bs{e}_K \otimes \frac{\bs{\partial}}{\bs{\partial} \xi^l}
}
Denoting the components in this representation as $D^{IjKl}$, the physical notion of strain arising from a constitutive law (as opposed to geometric linear strain, \ie the symmetric gradient of the displacements) is given by the vector-valued one-form
\eq{
	\star_{E, \bs{g}} \mbb{t} = \mbb{e} = \delta_{IJ} D^{IjKl} g_{jm} \sigma_{Kl} \, \bs{e}^J \otimes \extd \xi^m = (\mbb{D} : \star_\bs{g} \mbb{t})^\flat.
}
The map $\mbb{t} \mapsto (\mbb{D} : \star_\bs{g} \mbb{t})^\flat$ is termed the {\ital elasticity Hodge star}, $\star_{E, \bs{g}} : \Sigma \rightarrow \mc{E}$, which satisfies
\eq{
	\sigma_{Ij} D^{IjKl} \sigma_{Kl} \, {\star_\bs{g} 1} = \mbb{t} \overset{\cdot}{\wedge} \star_{E, \bs{g}} \mbb{t} \, .
}
A similar Hodge star was considered in terms of a DEC formulation in~\cite{Yavari2008}.

\subsubsection{Mixed formulation of linearized elasticity }

\begin{figure}[htb!]
\begin{minipage}{0.45\textwidth}
\begin{center}
\begin{tikzpicture}[every node/.style={inner sep=2,outer sep=0}]
\node (label) at (0,0)[]{
	\includegraphics[width=.8\textwidth]{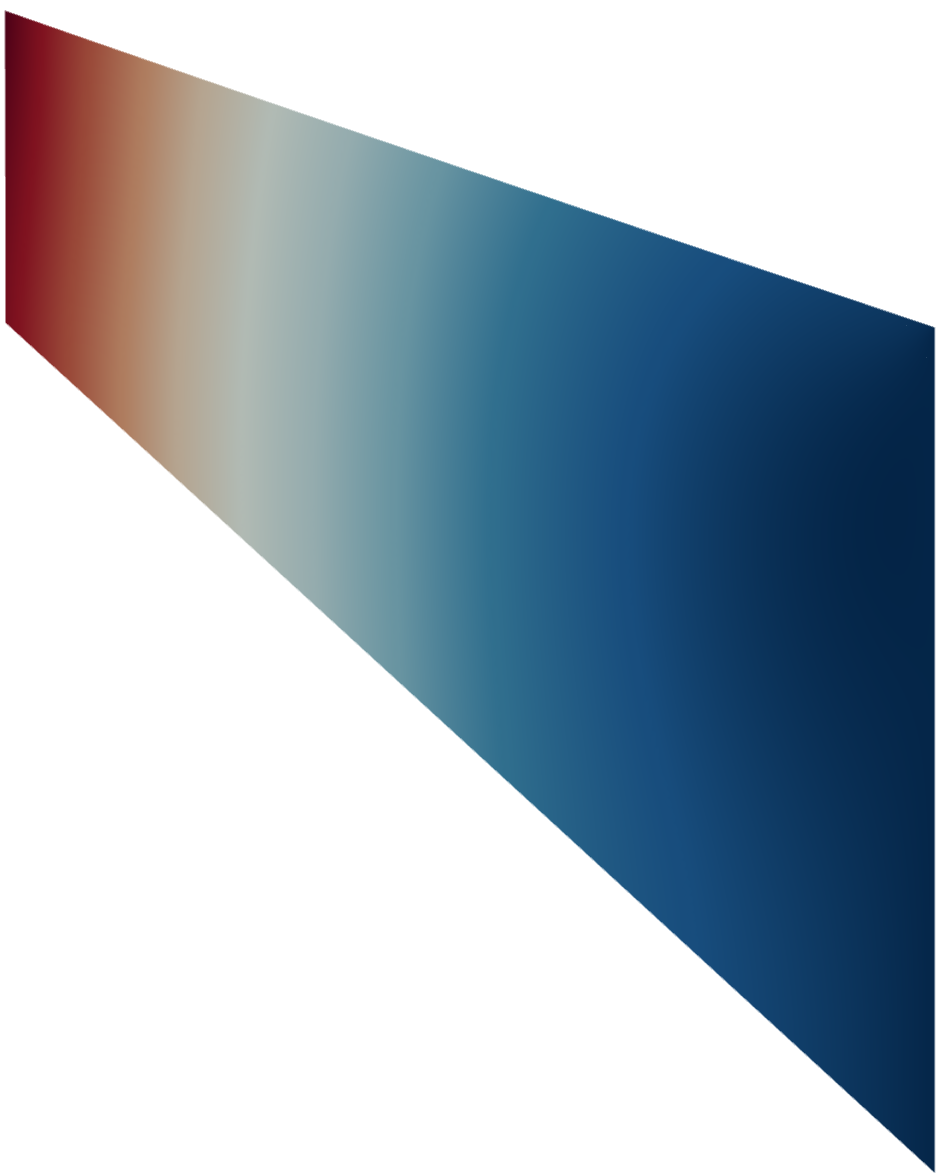}
      };
      \node (label) at (0.05,-4.0)[]{
\pgfplotscolorbardrawstandalone[ 
    colormap name = BuRd,
    colorbar horizontal,
    point meta min=-3.088e-2,
    point meta max=7.771,
    colorbar style={
        width=\textwidth-1.5cm,
        height=.25cm,
        }
        ]};
      \draw[-stealth, thick] (-2,-2.5)--(-2,-1.5) node[above] {$x_2$};  
      \draw[-stealth, thick] (-2,-2.5)--(-1,-2.5) node[right] {$x_1$};  
\end{tikzpicture}\\
(a) Vertical displacement
\vspace{.25cm}
\begin{tikzpicture}[every node/.style={inner sep=2,outer sep=0}]
\node (label) at (0,0)[]{
\includegraphics[width=\textwidth]{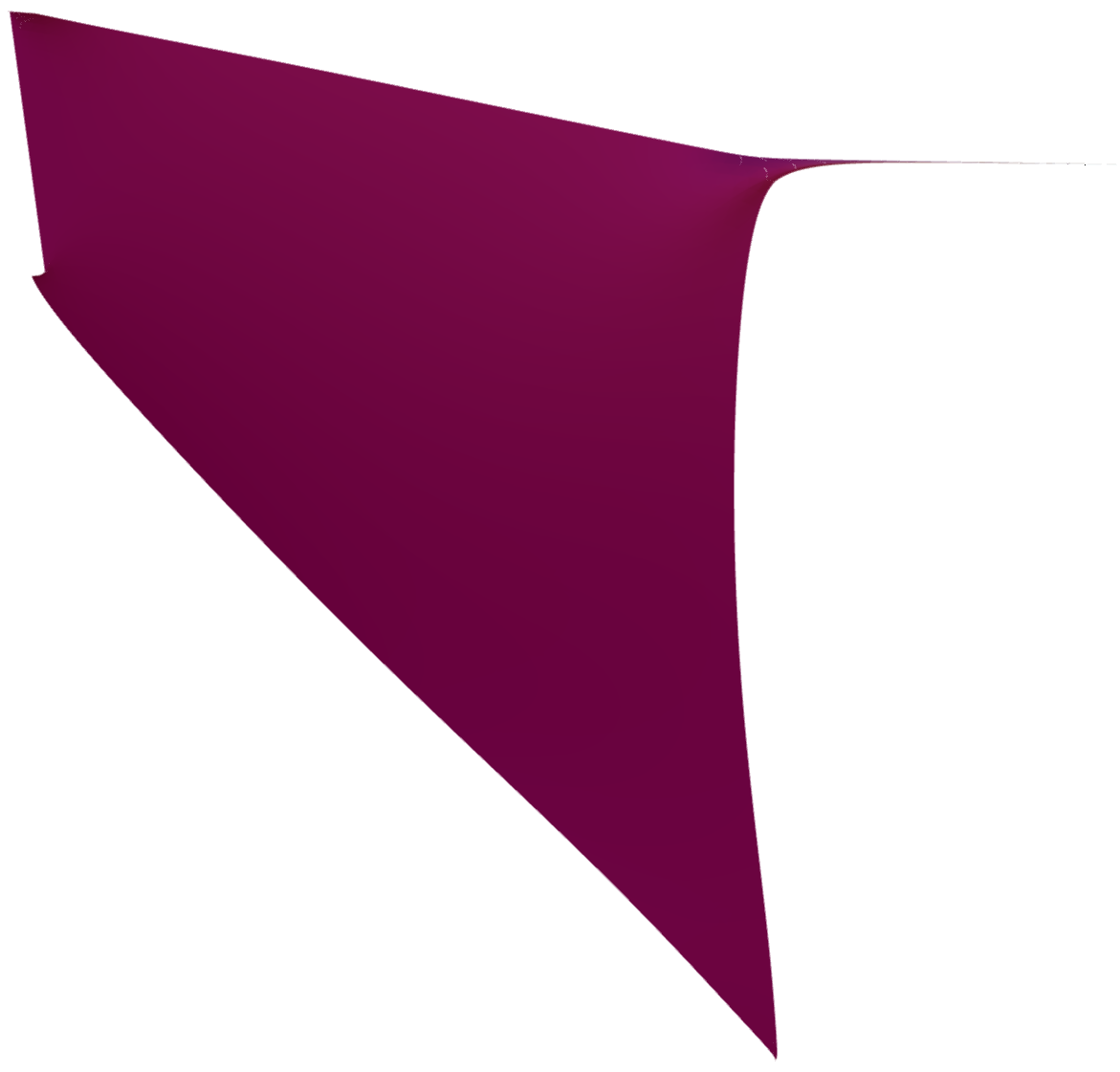}
      };
      \node (label) at (0.05,-3.9)[]{
\pgfplotscolorbardrawstandalone[ 
    colormap name = arnold,
    colorbar horizontal,
    point meta min=-519.1,
    point meta max=13.84,
    colorbar style={
        width=\textwidth-1.5cm,
        height=.25cm,
        }
        ]};
\draw[-stealth, thick] (2.25,-1.5)--(2.95,-1.7) node[right] {$x_1$};  
      \draw[-stealth, thick] (2.25,-1.5)--(2.17,-0.75) node[above] {$x_2$};  
            \draw[-stealth, thick] (2.25,-1.5)--(1.65,-1.65) node[below] {$x_3$};  
\end{tikzpicture}\\
(c) Hydrostatic stress $\frac{1}{3}\mrm{tr}(\star_{\bs{g}}\mbb{t})$ 
\end{center}
\end{minipage}
\hfill
\begin{minipage}{0.45\textwidth}
\centering
\begin{center}
\begin{tikzpicture}[every node/.style={inner sep=2,outer sep=0}]
\node (label) at (0,0)[]{
\includegraphics[width=.8\textwidth]{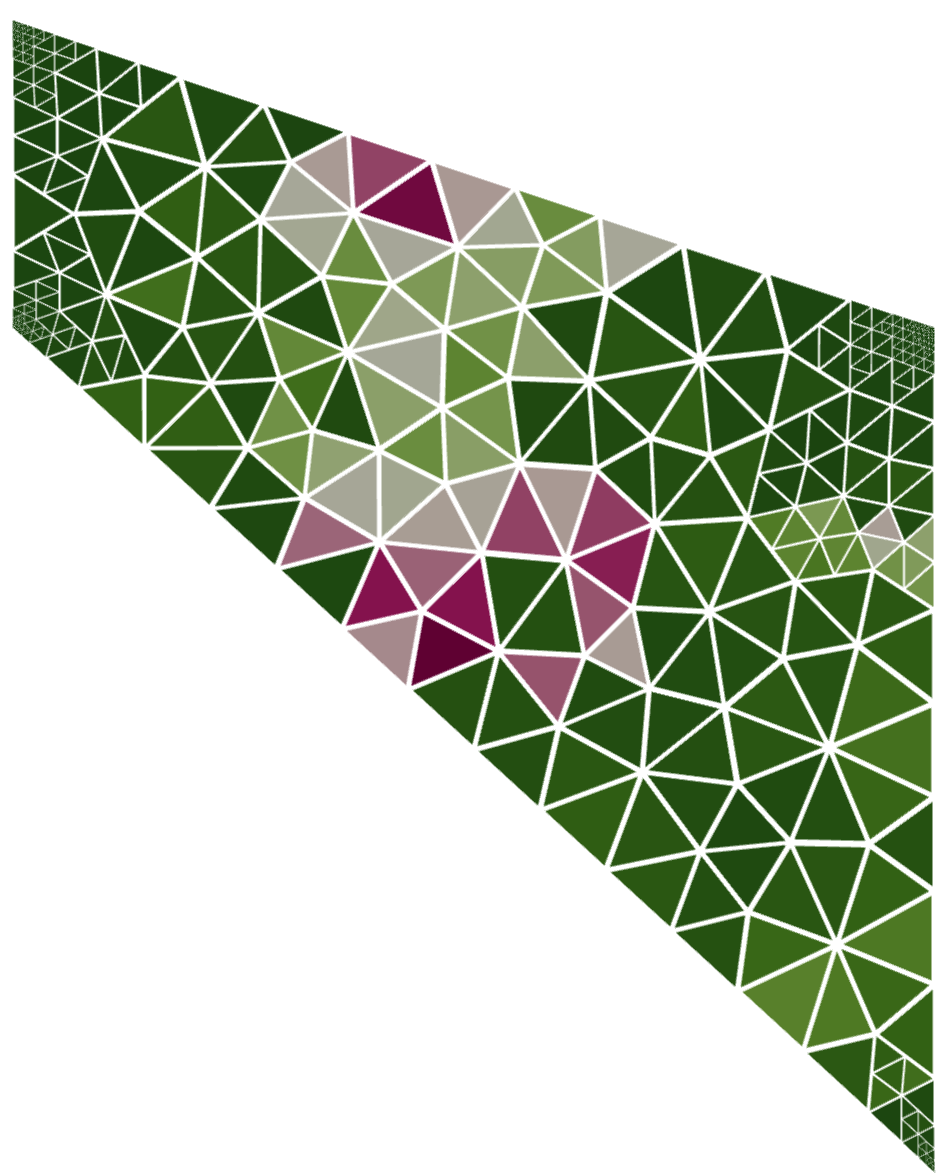}
      };
      \node (label) at (0.05,-4.0)[]{
\pgfplotscolorbardrawstandalone[ 
    colormap name = arnold,
    colorbar horizontal,
    point meta min=-0,
    point meta max=1.693,
    colorbar style={
        width=\textwidth-1.5cm,
        height=.25cm,
        }
        ]};
     \draw[-stealth, thick] (-2,-2.5)--(-2,-1.5) node[above] {$x_2$};  
      \draw[-stealth, thick] (-2,-2.5)--(-1,-2.5) node[right] {$x_1$};  
\end{tikzpicture}\\
(b) Relative error
\vspace{.25cm}
\begin{tikzpicture}[every node/.style={inner sep=2,outer sep=0}]
\node (label) at (0,0)[]{
\includegraphics[width=\textwidth]{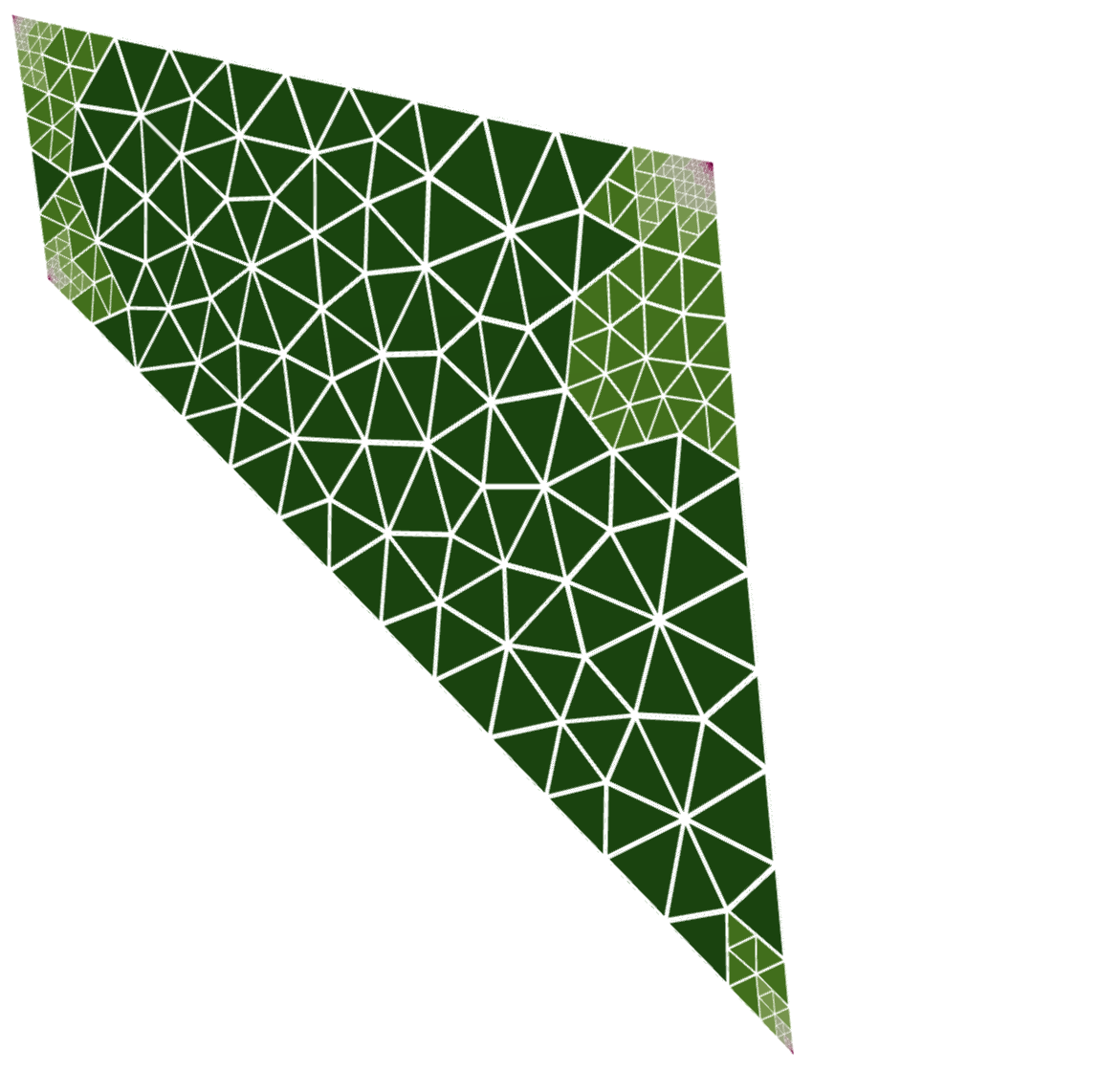}
      };
      \node (label) at (0.05,-3.9)[]{
\pgfplotscolorbardrawstandalone[ 
    colormap name = arnold,
    colorbar horizontal,
    point meta min=1,
    point meta max=8,
    colorbar style={
        width=\textwidth-1.5cm,
        height=.25cm,
        }
        ]};
      \draw[-stealth, thick] (2.25,-1.5)--(2.95,-1.7) node[right] {$x_1$};  
      \draw[-stealth, thick] (2.25,-1.5)--(2.17,-0.75) node[above] {$x_2$};  
\end{tikzpicture}\\
(d) Hierarchical subdivision depth
\end{center}
\end{minipage} \caption{Results for the plane strain linearized elasticity problem on Cook's membrane}
\label{fig:elasticityresult}
\end{figure}
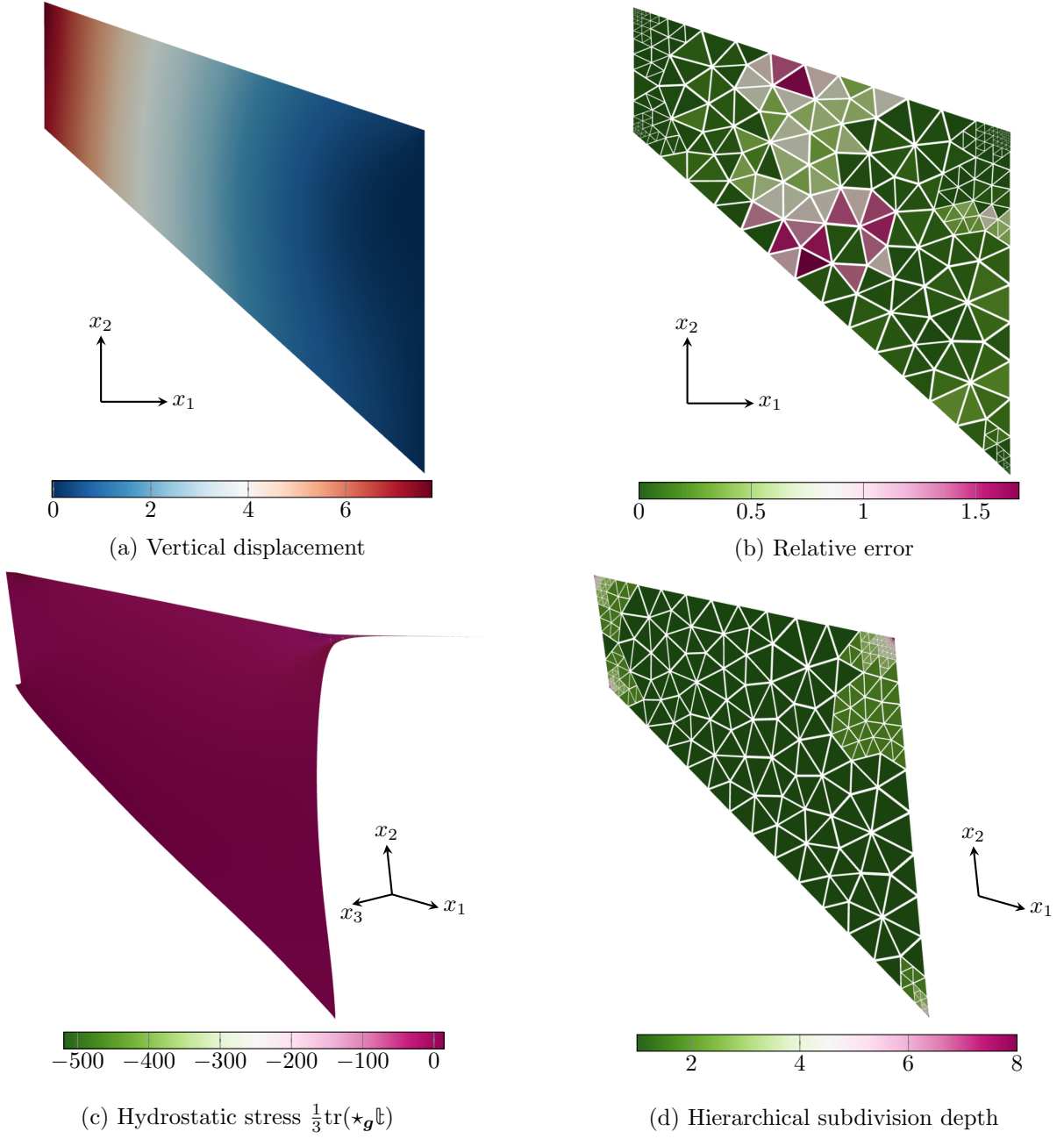

\begin{figure}[htb!]
\begin{center}
\begin{tikzpicture}
\node[inner sep=0pt] (sc1) at (-3.5,0)
    {\includegraphics[height=5.5cm]{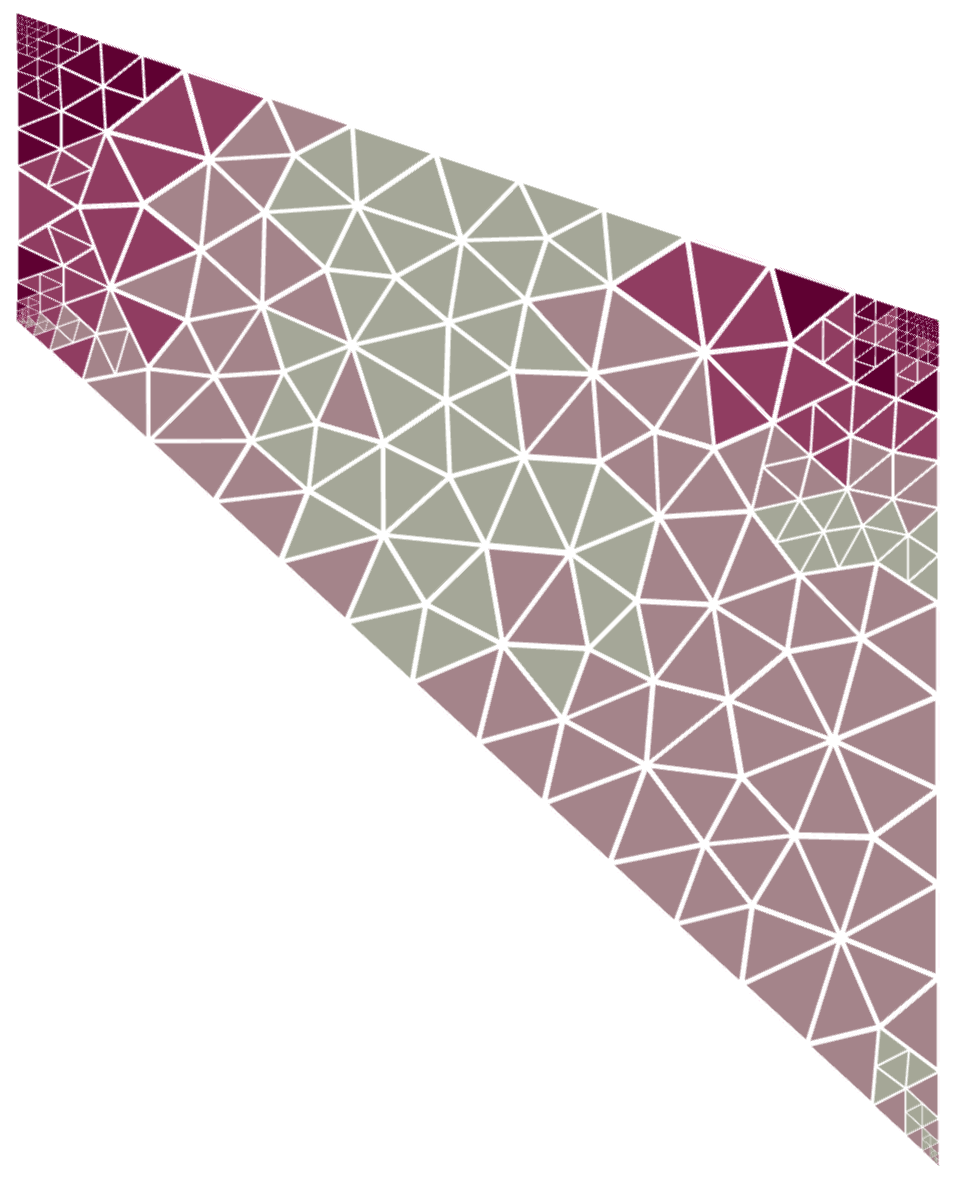}};
\node[draw,rectangle,inner sep=0pt,line width=2] (sc2) at (1.,0)
    {\includegraphics[height=4cm]{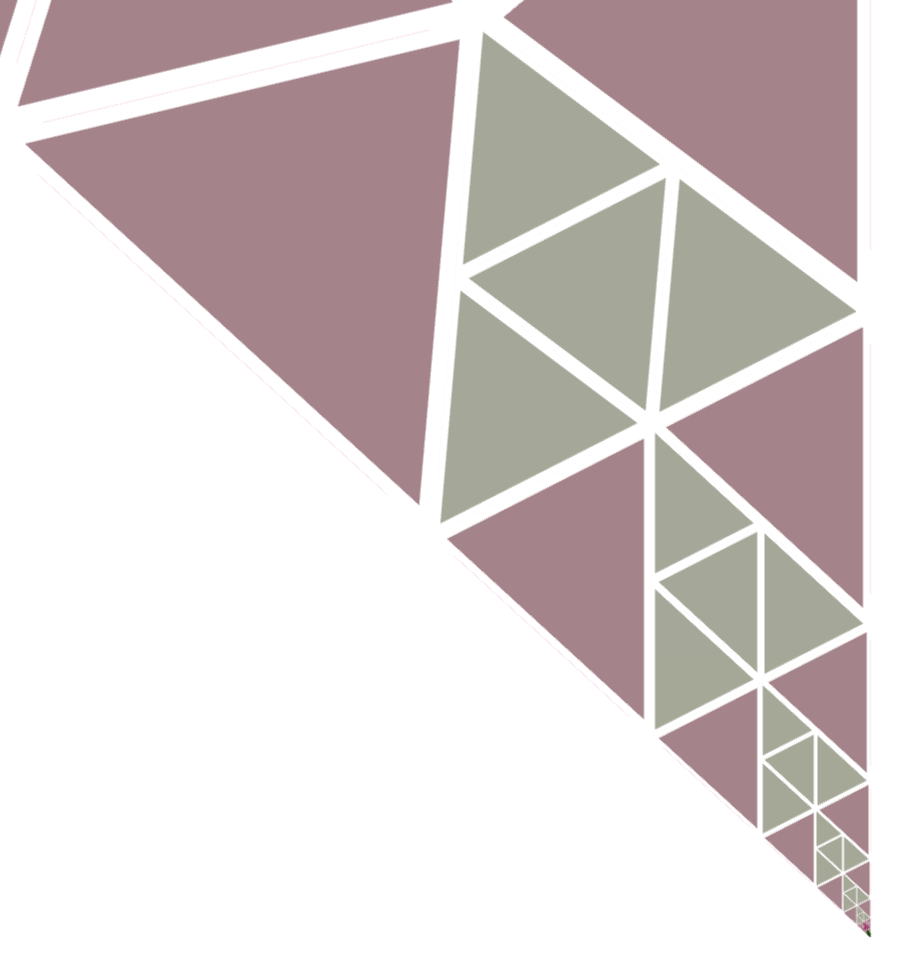}};
    \node[draw,rectangle,inner sep=0pt,line width=2] (sc3) at (4.85,0)
    {\includegraphics[height=4cm]{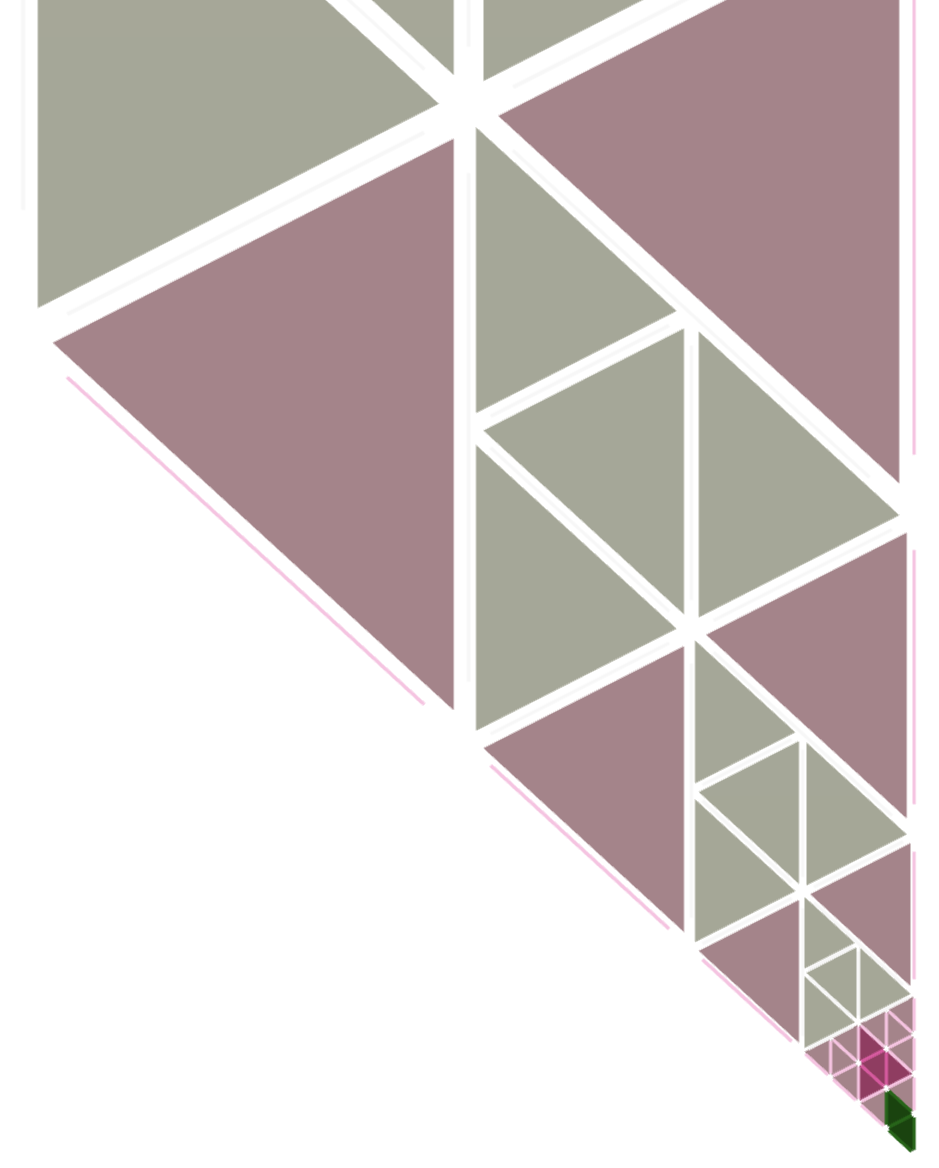}};
	\node[draw,rectangle,minimum height=.8cm,minimum width=.8cm,line width=1] (scb1) at (-1.65,-2.35){};
	\draw[-stealth,line width = 1] (scb1) --  (-0.25,-1.5);	
	\node[draw,rectangle,minimum height=1.25cm,minimum width=1.0cm,line width=1] (scb2) at (2.325,-1.325){};
	\draw[-stealth,line width = 1] (scb2) --  (3.75,-1.25);
	  \node (label) at (.35,-3.5)[]{
\pgfplotscolorbardrawstandalone[ 
    colormap name = arnold,
    colorbar horizontal,
    point meta min=1,
    point meta max=7,
    colorbar style={
	xtick={1,2,...,7},
        height=.25cm,
        width=5cm,
        }
        ]};	
\end{tikzpicture}
\end{center}
 \caption{Polynomial order distribution for the elasticity problem on Cook's membrane}
\label{fig:elasticitypolydist}
\end{figure}
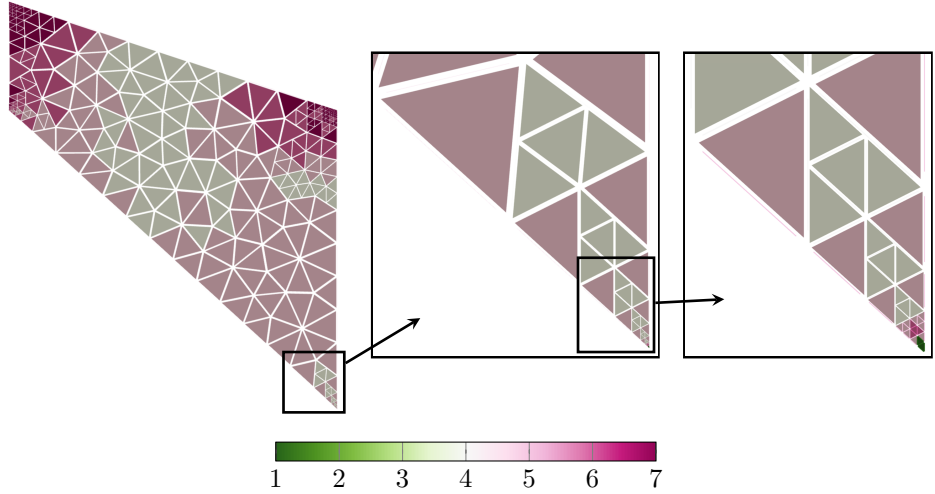

The linearized elasticity equations in terms of the linear momentum flux $\mbb{t} \in \Lambda^{n-1}(\Omega; \field{R}^n)$ and the displacement $\bs{u} \in \Lambda^n(\Omega; \field{R}^n)$, can be stated as
\eq{
	(-1)^{n-1}\extd \mbb{t} + \bs{b} = 0 \quad \quad {\star_{E,\bs{g}} \mbb{t}} = \extd{\star_\bs{g}\bs{u}} \quad \quad {\star_{\bs{g}} \mbb{t}}  \in \mrm{symm}
}
on $\Omega$. Since $\mc{S} = \mbb{R}^n$, the corresponding displacement vector field is written simply as
\eq{
	\star_\bs{g}\bs{u} = \bs{e}^I \otimes u_I \, .
}
A weakly-symmetric version is obtained by replacing the angular momentum balance equation by
\eq{
	 ({\star_{\bs{g}} \mbb{t}}) : \bs{q} = 0 \quad  \forall  \bs{q} \in \mrm{skw} \, .
}
The algebraic operator $S : \Lambda^{n-1}(\Omega; \field{R}^n) \rightarrow \Lambda^{n-1}(\Omega; \field{R}^n \wedge \field{R}^n)$ is introduced according to~\cite{Arnold2006}. Then, by taking inner products of the balance of linear momentum and the infinitesimal constitutive relation as well as by introducing a rotation-valued Lagrange multiplier term in the constitutive equation, the formulation
\eq{
	(-1)^{n-1}\bs{v} \overset{\cdot}{\wedge} \star_\bs{g} \extd \mbb{t} + \bs{v} \overset{\cdot}{\wedge} \star_\bs{g} \bs{b} = 0 \quad &\forall \bs{v} \in \mrm{\Lambda}^{n}(\Omega; \field{R}^n)\\
	\mbb{s} \overset{\cdot}{\wedge} \star_{E, \bs{g}} \mbb{t} = \mbb{s} \overset{\cdot}{\wedge}\extd{\star_\bs{g}\bs{u}} + S\mbb{s} \overset{\cdot}{\wedge} \star_\bs{g}\bs{p} \quad &\forall \mbb{s} \in \Lambda^{n-1}(\Omega; \field{R}^n) \\
	\bs{q} \overset{\cdot}{\wedge} \star_\bs{g}S\mbb{t} = 0 \quad &\forall\bs{q} \in \Lambda^{n}(\Omega; \field{R}^n \wedge \field{R}^n)
}
is obtained. Note that for the wedge-product of the bivector-valued differential forms, we mean 
\eq{
	\bs v \overset{\cdot}{\wedge} \bs w \mapsto (\bs v_{\bs{\sigma}} \wedge \star\bs w_{\bs{\nu}}) \cobasis{x}{\bs{\sigma}} \wedge  \cobasis{x}{\bs{\nu}} 
}
where the Hodge star is taken to be related to the ambient space $\mc S = \field{R}^n$.
Applying Leibniz' rule to the weak constitutive equation,
\eq{
	\mbb{s} \overset{\cdot}{\wedge} \star_{E, \bs{g}} \mbb{t} = (-1)^{n-1} \extd(\mbb{s} \overset{\cdot}{\wedge}{\star_\bs{g}\bs{u}}) + (-1)^n \extd \mbb{s} \overset{\cdot}{\wedge} {\star_\bs{g}\bs{u}} + S\mbb{s} \overset{\cdot}{\wedge} \star_\bs{g}\bs{p}
}
follows.
By Stokes' theorem applied to the left-most term on the right-hand-side of the latter equation,
\eq{
	\int\limits_{\Omega}\extd(\mbb{s} \overset{\cdot}{\wedge}{\star_\bs{g}\bs{u}}) = \int\limits_{\partial \Omega}\mbb{s} \overset{\cdot}{\wedge}{\star_\bs{g}\bs{u}} \,,
}
the displacement boundary condition is found. The final weak form then reads: find $\mbb{t} \in H\Lambda^{n-1}(\Omega; \mbb{R}^n)$, $\bs{u} \in L^2\Lambda^n(\Omega; \mbb{R}^n)$ and $\bs{p} \in L^2\Lambda^n(\Omega, \Lambda^2\mbb{R}^n)$ such that
\eq{
	(-1)^{n}\int\limits_{\Omega}\bs{v} \overset{\cdot}{\wedge} \star_\bs{g} \extd \mbb{t}   = \int\limits_{\Omega}\bs{v} \overset{\cdot}{\wedge} \star_\bs{g} \bs{b}  \quad &\forall \bs{v} \in L^2\mrm{\Lambda}^{n}(\Omega; \field{R}^n)\\
	\int\limits_{\Omega}\mbb{s} \overset{\cdot}{\wedge} \star_{E, \bs{g}} \mbb{t} + (-1)^{n-1} \int\limits_{\Omega}\extd \mbb{s} \overset{\cdot}{\wedge} {\star_\bs{g}\bs{u}} + \int\limits_{\Omega}S\mbb{s} \overset{\cdot}{\wedge} \star_\bs{g}\bs{p} = (-1)^{n-1} \int\limits_{\partial \Omega}\mbb{s} \overset{\cdot}{\wedge}{\star_\bs{g}\bs{u}} \quad &\forall \mbb{s} \in H\Lambda^{n-1}(\Omega; \field{R}^n) \\
	\int\limits_{\Omega}  \bs{q}\overset{\cdot}{\wedge} \star_\bs{g}S\mbb{t} = 0 \quad &\forall\bs{q} \in L^2\Lambda^{n}(\Omega; \Lambda^2\field{R}^n)\,,
}
where the stress-boundary condition must be enforced strongly, as outlined in the following section. We discretize the problem using the spaces $\mc{P}_{r}\Lambda^{n-1}(\Omega_\mrm{d}; \field{R}^n)$ for stress, $\mc{P}_{r-1}\Lambda^{n}(\Omega_\mrm{d}; \field{R}^n)$ for displacement and $\mc{P}_{r-1}\Lambda^{n}(\Omega_\mrm{d}; \Lambda^2\field{R}^n)$ for rotation~\cite{Arnold2006}. Note that if desired, the infinitesimal rotation can serve as an error indicator, however this is not pursued further.

\subsubsection{Strong imposition of stress boundary conditions }

For enforcing the stress boundary condition $\bs{\sigma}\cdot\bs{n} = \bar{\bs{t}}$ on $\partial_\sigma \Omega$, the required relation between Cauchy stress and linear momentum flux for each boundary chart indexed by $i$ is
\eq{
	\sigma_{Ij} n^j \, \bs{e}^I = ({\star_\bs{g}\mbb{t}})(\bs{n}) \quad  \tx{on} \quad \eta_i(\partial_\sigma \Omega)\,,
}
where $n^j$ are the components of the vector normal to $\eta_i(\partial_\sigma \Omega)$ obtained by pushing the unit normal vector of $\partial_\sigma \Omega$ forward along the map $\eta_i$. Imposing the boundary condition leads to the constraint equation
\eq{
	({\star_\bs{g}\mbb{t}})(\bs{n})   - \bar{\bs{t}} = 0 \quad \tx{on} \quad \eta_i(\partial_\sigma \Omega) \, .
}
The stress can be expanded in the reference element basis as $\mbb{t} = \hat{\sigma}_{\alpha,I}  \bs{e}^I \otimes \phi^\sigma_\alpha$ with $\phi^\sigma_\alpha \in \mc{P}_r\Lambda^{n-1}(\simplx{T})$, where the basis functions $\phi^\sigma$ indexed by $\alpha \in \mc{I}(\eta_i(\partial_\sigma \Omega))$ are assumed to have an influence at the boundary $\eta_i(\partial_\sigma \Omega)$ such that $(\star_\bs{g} \phi^\sigma_\alpha)(\bs{n}) = (\star_\bs{g}\phi_\alpha^\sigma)_j n^j \neq 0$ on $\eta_i(\partial_\sigma \Omega)$.
With $B^\sigma_{j\alpha} = (\star_\bs{g} \phi^\sigma_\alpha)_j$, the constraint becomes
\eq{
	\hat{\sigma}_{\alpha,I}  \bs{e}^I  B^\sigma_{j\alpha} n^j   - \bar{\bs{t}} = 0 \quad \tx{on} \quad \eta_i(\partial_\sigma \Omega) \, .
}
To ensure the existence of solution coefficients $\hat{\bs{\sigma}}_{\alpha}$ satisfying this constraint, the imposed stress vector field $\bar{\bs{t}}$ must be exactly representable in the space of the basis functions employed for stress approximation at the boundary. As such a restriction on the stress vector field is impractical, we solve for the unknown coefficients $\hat{\bs{\sigma}}_{\alpha}$ by seeking the Galerkin projection of $\bar{\bs{t}}$ onto the subspace spanned by the boundary basis functions. For each spatial direction $J \in (1\ldots n)$ and boundary chart indexed by $i$, this results in solving the linear system
\eq{
	\int\limits_{\eta_i(\partial_\sigma \Omega)}   n^i n^j   \, \mrm{Tr} (B^\sigma_{i\alpha} B^\sigma_{j\beta}  ) \;  \hat{\sigma}_{\beta,J}  = \int\limits_{\eta_i(\partial_\sigma \Omega)}  n^i  \bar{t}_J \mrm{Tr} ( B^\sigma_{i\alpha}   )   \,, \quad \alpha, \beta \in \mc{I}(\eta_i(\partial_\sigma \Omega)) \,.
}
The coefficients representing the imposed stress field at the boundary can then be incorporated into the global system of equations for the elasticity problem, effectively reducing its size.

\subsubsection{Cook's membrane in plane strain}

For the incompressible Cook's membrane problem in plane strain, we choose Young's modulus $E = 250$ and Poisson's ratio $\nu = 0.5$. As boundary conditions we prescribe $\star \bs{u} = 0$ (right boundary), $(\star\mathbb{t})(\bs{n}) = 0$ (top and bottom boundaries) and $(\star\mathbb{t})(\bs{n}) = (0, 100/16)^\mrm{T}$ (units of force per length squared, left boundary). The problem is 48 by 60 units of length, the length of the left boundary is 16, and the length of the right boundary is 44.

Results of the computation are displayed in Fig.~\ref{fig:elasticityresult}. The maximum vertical displacement in $x_2$-direction is approximately 7.771 units of length. The response is softer than the result given in~\cite{Nakshatrala2008} which is plausible considering the high accuracy due to the local refinement. A perspective projection is used for the hydrostatic stress in order to show the singularities more clearly. As can be seen, the compressive singularity at the re-entrant corner is well approximated. This is also the case for the other singularities. The minimum of the hydrostatic stress, which takes on a value of $-519.1$ units of force per length squared, is located at the re-entrant corner.\footnote{We note that at sufficient magnification, discontinuities seem to appear in the hydrostatic stress rendering. This is not a deficiency of the solution but simply a consequence of our post-processing routine which only coarsely samples the polynomials on simplices at high refinement depths.} As a reference, the subdivided mesh is rendered using the same projection. In Fig.~\ref{fig:elasticitypolydist}, the polynomial order distribution is depicted in order to demonstrate the lowest-order approximation within the singularity. Although not visible in print, the refinement algorithm imposes lowest-order approximation in all singularities present. 

The FEEC discretization of elasticity is sometimes criticized for its large number of degrees of freedom per element, even the lowest-order case. While this is obviously true, the remarkable accuracy obtained hints that the effort may well be worth while in complex problems amenable to adaptivity. Owing to an efficient implementation, computational costs were moderate and computations were carried out on a standard workstation.

\section{Conclusion and outlook}

A simple method for $hp$-hierarchical discretization of the finite element exterior calculus on piecewise Riemannian manifolds was outlined and a problem-independent error indicator was presented. The method as well as the algebraic concepts were implemented in the Julia language and verified using the Hodge Laplacian and incompressible linearized elasticity. The method is easily extended to relocation refinements such as those arising from variational ALE methods~\cite{Thoutireddy2003}. Further natural applications for an adaptive FEEC as outlined herein include spacetime~\cite{Salamon2014}, finite elasticity~\cite{Angoshtari2016} and elastodynamics problems~\cite{Arnold2014}. The framework for adaptivity in higher dimensions developed herein is especially key in the four-dimensional spacetime case, as it allows for spatially inhomogeneous time-subdivisions tailored to the solution, leading to efficient methods. While the finite element exterior calculus appears to be ideally-suited to multi-physics problems, it remains to be seen what practical gain can be achieved in problems such as those arising in magnetohydrodynamics.

\clearpage

\bibliography{hpfeec}{}
\bibliographystyle{wileyj}
\end{document}